\documentclass[11pt,twoside]{article}

\usepackage{amsfonts,epsfig,graphicx}
\usepackage{amsmath,amssymb,amsthm}
\usepackage{fullpage}

\usepackage{epsf}
\usepackage{fancyheadings}
\usepackage{graphics}
\usepackage{amsfonts}
\usepackage{amsmath}
\usepackage{psfrag}
\usepackage{mathrsfs}

\usepackage{color}

\usepackage{multirow,tabularx}
\usepackage{hhline}

\newcolumntype{L}[1]{>{\raggedright\let\newline\\\arraybackslash\hspace{0pt}}m{#1}}
\newcolumntype{C}[1]{>{\centering\let\newline\\\arraybackslash\hspace{0pt}}m{#1}}
\newcolumntype{R}[1]{>{\raggedleft\let\newline\\\arraybackslash\hspace{0pt}}m{#1}}

\theoremstyle{plain}

\newtheorem{theo}{Theorem}[section]

\newtheorem{lem}{Lemma}[section]
\newtheorem{prop}{Proposition}[section]
\newtheorem{cor}{Corollary}[section]

\theoremstyle{definition} 

\newtheorem{nota}{Notation}[section]
\newtheorem{de}{Definition}[section]
\newtheorem{exa}{Example}[section]
\newtheorem{as}{Assumption}[section]
\newtheorem{alg}{Algorithm}[section]

\newcommand{\btheo}{\begin{theo}}
\newcommand{\bde}{\begin{de}}
\newcommand{\ble}{\begin{lem}}
\newcommand{\bpr}{\begin{prop}}
\newcommand{\bno}{\begin{nota}}
\newcommand{\bex}{\begin{exa}}
\newcommand{\bcor}{\begin{cor}}
\newcommand{\spro}{\begin{proof}}
\newcommand{\bas}{\begin{as}}
\newcommand{\balg}{\begin{alg}}

\newcommand{\etheo}{\end{theo}}
\newcommand{\ede}{\end{de}}
\newcommand{\ele}{\end{lem}}
\newcommand{\epr}{\end{prop}}
\newcommand{\eno}{\end{nota}}
\newcommand{\eex}{\end{exa}}
\newcommand{\ecor}{\end{cor}}
\newcommand{\fpro}{\end{proof}}
\newcommand{\eas}{\end{as}}
\newcommand{\ealg}{\end{alg}}

\theoremstyle{plain}

\newtheorem{theos}{Theorem}
\newtheorem{props}{Proposition}
\newtheorem{lems}{Lemma}
\newtheorem{cors}{Corollary}

\theoremstyle{definition}
\newtheorem{exas}{Example}
\newtheorem{algs}{Algorithm}
\newtheorem{asss}{Assumption}
\newtheorem{defns}{Definition}

\newcommand{\btheos}{\begin{theos}}
\newcommand{\etheos}{\end{theos}}
\newcommand{\bprops}{\begin{props}}
\newcommand{\eprops}{\end{props}}
\newcommand{\bdes}{\begin{defns}}
\newcommand{\edes}{\end{defns}}
\newcommand{\blems}{\begin{lems}}
\newcommand{\elems}{\end{lems}}
\newcommand{\bcors}{\begin{cors}}
\newcommand{\ecors}{\end{cors}}
\newcommand{\bexs}{\begin{exas}}
\newcommand{\eexs}{\end{exas}}
\newcommand{\balgs}{\begin{algs}}
\newcommand{\ealgs}{\end{algs}}
\newcommand{\bass}{\begin{asss}}
\newcommand{\eass}{\end{asss}}

\newcommand{\numobs}{\ensuremath{n}}
\newcommand{\posterior}{\ensuremath{\pi_n}}
\newcommand{\covdim}{\ensuremath{p}}
\newcommand{\real}{\ensuremath{\mathbb{R}}}
\newcommand{\defn}{\ensuremath{: \, =}}
\newcommand{\inprod}[2]{\ensuremath{\langle #1 , \, #2 \rangle}}

\newcommand{\Gfun}{\ensuremath{\mathcal{G}}}
\newcommand{\betastar}{\ensuremath{\beta^\ast}}
\newcommand{\MarkovChain}{\ensuremath{\mathcal{C}}}

\newcommand{\Cm}{\ensuremath{\nu}}

\newcommand{\yobs}{\ensuremath{Y}}
\newcommand{\Sset}{\ensuremath{S}}
\newcommand{\regu}{\ensuremath{\lambda}}
\newcommand{\maxsize}{\ensuremath{s_0}}
\newcommand{\gammastar}{\ensuremath{{\gamma^\ast}}}
\newcommand{\Gammastar}{\ensuremath{\Gamma^\ast}}
\newcommand{\sstar}{\ensuremath{s^\ast}}
\newcommand{\hyperpara}{\ensuremath{g}}

\newcommand{\design}{\ensuremath{X}}
\newcommand{\Xmat}{\ensuremath{\design}}
\newcommand{\wnoise}{\ensuremath{w}}
\newcommand{\wtil}{\ensuremath{\tilde{w}}}
\newcommand{\Proj}{\ensuremath{\Phi}}

\newcommand{\Exs}{\ensuremath{\mathbb{E}}}

\newlength{\widebarargwidth}
\newlength{\widebarargheight}
\newlength{\widebarargdepth}

\newcommand{\Prob}{\ensuremath{\mathbb{P}}}

\newcommand{\Ind}{\ensuremath{\mathbb{I}}}

\newcommand{\widgraph}[2]{\includegraphics[keepaspectratio,width=#1]{#2}}

\newcommand{\UNICON}{\ensuremath{c}}

\newcommand{\matsnorm}[2]{|\!|\!| #1 | \! | \!|_{{#2}}}

\newcommand{\opnorm}[1]{\ensuremath{\matsnorm{#1}{\tiny{\mbox{op}}}}}

\newcommand{\Aevent}{\ensuremath{\mathcal{A}}}

\newcommand{\mprob}{\ensuremath{\mathbb{P}}}
\newcommand{\neigh}{\ensuremath{\mathcal{N}}}
\newcommand{\usedim}{\covdim}

\newcommand{\CB}{C_\beta}

\newcommand{\MODEL}{\ensuremath{\mathbb{M}}}
\newcommand{\Mspace}{\ensuremath{\mathscr{M}}}

\long\def\comment#1{}

\newcommand{\plprior}{\ensuremath{\pi}}

\newcommand{\hamm}{\ensuremath{d_H}}

\newcommand{\CONTWO}{\ensuremath{c'}}

\newcommand{\sigmazero}{\ensuremath{\sigma_0}}

\newcommand{\order}{\ensuremath{\mathcal{O}}}

\newcommand{\betamin}{\ensuremath{\beta_{\mbox{\tiny{min}}}}}

\newcommand{\Ltil}{\ensuremath{\widetilde{L}}}

\newcommand{\betahat}{\ensuremath{\widehat{\beta}}}
\newcommand{\gammatil}{\ensuremath{\widetilde{\gamma}}}

\newcommand{\lammin}{\ensuremath{\lambda_{\mbox{\tiny{min}}}}}
\newcommand{\lammax}{\ensuremath{\lambda_{\mbox{\tiny{max}}}}}
\newcommand{\Sigmabad}{\ensuremath{\Sigma_{\mbox{\tiny{bad}}}}}

\newcommand{\rcon}{\ensuremath{C_0}}

\newcommand{\calpha}{\ensuremath{C_1}}

\newcommand{\Path}{\ensuremath{T}}
\newcommand{\PathEns}{\ensuremath{\mathcal{T}}}

\newcommand{\GAP}{\ensuremath{\mbox{Gap}}}

\newcommand{\PMAT}{\ensuremath{\mathbf{P}}}
\newcommand{\QMAT}{\ensuremath{\mathbf{Q}}}
\newcommand{\SMAT}{\ensuremath{\mathbf{S}}}
\newcommand{\ProbMH}{\ensuremath{\PMAT_{\tiny{\mbox{MH}}}}}

\newcommand{\ACCEPT}{\ensuremath{\mathbf{R}}}

\newcommand{\gammabar}{\ensuremath{\bar{\gamma}}}

\newcommand{\PREC}{\ensuremath{\Lambda}}
\newcommand{\Stat}{\ensuremath{\pi}}

\newcommand{\Bevent}{\ensuremath{\mathcal{B}}}
\newcommand{\Cevent}{\ensuremath{\mathcal{C}}}
\newcommand{\Devent}{\ensuremath{\mathcal{D}}}

\newcommand{\spindex}{\ensuremath{s}}
\newcommand{\RE}{\ensuremath{\operatorname{RE}}}
\newcommand{\IC}{\ensuremath{\operatorname{SI}}}
\newcommand{\HACKER}{sparse projection~}
\newcommand{\rowdesign}{\ensuremath{x}}
\newcommand{\SNR}{\ensuremath{\mbox{SNR}}}

\newcommand{\MspaceMax}{\ensuremath{\Mspace(\maxsize)}}

\makeatletter
\long\def\@makecaption#1#2{
        \vskip 0.8ex
        \setbox\@tempboxa\hbox{\small {\bf #1:} #2}
        \parindent 1.5em  %% How can we use the global value of this???
        \dimen0=\hsize
        \advance\dimen0 by -3em
        \ifdim \wd\@tempboxa >\dimen0
                \hbox to \hsize{
                        \parindent 0em
                        \hfil 
                        \parbox{\dimen0}{\def\baselinestretch{0.96}\small
                                {\bf #1.} #2
                                } 
                        \hfil}
        \else \hbox to \hsize{\hfil \box\@tempboxa \hfil}
        \fi
        }
\makeatother

\begin{document}

\begin{center} {\LARGE{\bf{On the Computational Complexity of 
High-Dimensional Bayesian Variable Selection}}} \\
  \vspace{1cm}

\begin{center}
\begin{tabular}{ccccc}
Yun Yang$^1$ & &  Martin J. Wainwright$^{1,2}$ & & Michael I. Jordan$^{1,2}$
\end{tabular}
\end{center}

  \vspace{1cm}
  {\large University of California, Berkeley} \\
  \vspace{.15cm} $^1$Department of Electrical Engineering and Computer
  Science ~~~~ $^2$Department of Statistics

\vspace*{.2in}

\today

\vspace*{.5in}

\begin{abstract}
We study the computational complexity of Markov chain Monte Carlo
(MCMC) methods for high-dimensional Bayesian linear regression under
sparsity constraints. We first show that a Bayesian approach can
achieve variable-selection consistency under relatively mild
conditions on the design matrix.  We then demonstrate that the
statistical criterion of posterior concentration need not imply the
computational desideratum of rapid mixing of the MCMC algorithm.  By
introducing a truncated sparsity prior for variable selection, we
provide a set of conditions that guarantee both variable-selection
consistency and rapid mixing of a particular Metropolis-Hastings
algorithm. The mixing time is linear in the number of covariates up to
a logarithmic factor.  Our proof controls the spectral gap of the
Markov chain by constructing a canonical path ensemble that is
inspired by the steps taken by greedy algorithms for variable
selection.
\end{abstract}

\end{center}

%%%%%%%%%%%%%%%%%%%%%%%%%%%%%%%%%%%%%%%%%%%%%%%%%%%%%%%%%%%%%%%%%%%%%%%

\section{Introduction}

In many areas of science and engineering, it is common to collect a
very large number of covariates $\design_1,\ldots,\design_\covdim$ in
order predict a response variable $\yobs$.  We are thus led to
instances of high-dimensional regression, in which the number of
covariates $\covdim$ exceed the sample size $\numobs$.  A large
literature has emerged to address problems in the regime $\covdim 
\gg \numobs$, where the ill-posed nature of the problem is addressed
by imposing sparsity conditions---namely, that the response $\yobs$ 
depends only on a small subset of the covariates.  Much of this 
literature is based on optimization methods, where penalty terms 
are incorporated that yield both convex~\cite{Tibshirani96} and 
nonconvex~\cite{FanLi91,Zha12} optimization problems.  Theoretical 
analysis is based on general properties of the design matrix and the 
penalty function.

Alternatively, one can take a Bayesian point of view on high-dimensional
regression, placing a prior on the model space and performing the 
necessary integration so as to obtain a posterior distribution~\cite{George1993,
Ishwaran2005,Anirban2015}.  Obtaining such a posterior allows one to 
report a subset of possible models along with their posterior probabilities 
as opposed to a single model.  One can also report the marginal posterior 
probability of including each covariate.  Some recent work has provided some
theoretical understanding of the performance of Bayesian approaches to 
variable selection.  In the moderate-dimension scenario (in which
$\covdim$ is allowed to grow with $\numobs$ but $\covdim\leq
\numobs$), Shang and Clayton~\cite{Shang2011} establish posterior
consistency for variable selection in a Bayesian linear model, meaning
that the posterior probability of the true model that contains all
influential covariates tends to one as $\numobs$ grows. Narisetty and
He~\cite{Narisetty2014} consider a high-dimensional scenario in which
$\covdim$ can grow nearly exponentially with $\numobs$; in this
setting, they show the Bayesian spike-and-slab variable-selection
method achieves variable-selection consistency. Since this particular
Bayesian method resembles a randomized version of $\ell_0$-penalized
methods, it could have better performance than $\ell_1$-penalized
methods for variable selection under high-dimensional
settings~\cite{Narisetty2014,Shen2012}.  Empirical evidence for this
conjecture is provided by Guan et al.~\cite{Guan2011} for SNP
selection in genome-wide association studies, but it has not been
confirmed theoretically.

The most widely used tool for fitting Bayesian models are sampling
techniques based on Markov chain Monte Carlo (MCMC), in which a Markov
chain is designed over the parameter space so that its stationary
distribution matches the posterior distribution.  Despite its
popularity, the theoretical analysis of the computational efficiency
of MCMC algorithms lags that of optimization-based methods. In
particular, the central object of interest is the \emph{mixing time}
of the Markov chain, which characterizes the number of iterations
required to converge to an $\epsilon$-distance of stationary
distribution from any initial configuration.  In order for MCMC
algorithms to be controlled approximations, one must provide
meaningful bounds on the mixing time as a function of problem
parameters such as the number of observations and the dimensionality.
Of particular interest is determining whether the chain is
\emph{rapidly mixing}---meaning that the mixing time grows at most
polynomially in the problem parameters---or \emph{slowly mixing}
meaning that the mixing time grows exponentially in the problem
parameters.  In the latter case, one cannot hope to obtain approximate
samples from the posterior in any reasonable amount of time for large
models.

Unfortunately, theoretical analysis of mixing time is comparatively
rare in the Bayesian literature, with a larger number of negative
results.  On the positive side, Jones and Hobert~\cite{Jones2004}
consider a Bayesian hierarchical version of the one-way random effects
model, and obtain upper bounds on the mixing time of Gibbs and block
Gibbs samplers as a function of the initial values, data and
hyperparameters.  Belloni and Chernozhukov~\cite{Belloni2009} show
that a Metropolis random walk is rapidly mixing in the dimension for
regular parametric models in which the posterior converges to a normal
limit.  It is more common to find negative results in the literature.
Examples include Mossel and Vigoda~\cite{Mossel2006}, who show that
the MCMC algorithm for Bayesian phylogenetics takes exponentially long
to reach the stationary distribution as data accumulates, and Woodard
and Rosenthal~\cite{Woodard2013}, who analyze a Gibbs sampler used for
genomic motif discovery and show that the mixing time increases
exponentially as a function of the length of the DNA sequence.

The goal of the current paper is to study the computational complexity
of Metropolis-Hastings procedures for high-dimensional Bayesian
variable selection.  For concreteness, we focus our analysis on a
specific hierarchical Bayesian model for sparse linear regression, and
an associated Metropolis-Hastings random walk, but these choices
should be viewed as representative of a broader family of methods.  In
particular, we study the well-known Zellner $\hyperpara$-prior for
linear regression~\cite{Zellner1986}. The main advantage of this prior
is the simple expression that it yields for the marginal likelihood,
which is convenient in our theoretical investigations. As in past
analyses~\cite{Narisetty2014}, we consider the marginal probability of
including each covariate into the model as being on the order of
$\covdim^{-\mathcal{O}(1)}$.  Moreover, we restrict the support of the
prior to rule out unrealistically large models. As a specific computational methodology, 
we focus on an iterative, local-move and neighborhood-based procedure for sampling from the model space, which is motivated by the shotgun stochastic search~\cite{Hans2007}.

Our main contribution is to provide
conditions under which Bayesian posterior consistency holds, and
moreover, the mixing time grows linearly in $\covdim$ (up to
logarithmic factor), implying that the chain is rapidly mixing. As a
by-product, we provide conditions on the hyper-parameter $\hyperpara$
to achieve model-selection consistency.  We also provide a
counter-example to illustrate that although ruling out unrealistically
large models is not necessary for achieving variable-selection
consistency, it is necessary in order that the Metropolis-Hastings
random walk is rapidly mixing.  To be clear, while our analysis
applies to a fully Bayesian procedure for variable selection, it is
based on a frequentist point of view in assuming that the data are
generated according to a true model.

There are a number of challenges associated with characterizing the
computational complexity of Markov chain methods for Bayesian models.
First, the posterior distribution of a Bayesian model is usually a
much more complex object than the highly structured distributions in
statistical physics for which meaningful bounds on the Markov chain
mixing times are often obtained (e.g.\!  \cite{Borgs99},
\cite{Martinelli2012}, \cite{Levin2010}). Second, the transition
probabilities of the Markov chain are themselves stochastic, since
they depend on the underlying data-generating process.  In order to
address these challenges, our analysis exploits asymptotic properties
of the Bayesian model to characterize the typical behavior of the
Markov chain. We show that under conditions leading to Bayesian
variable-selection consistency, the Markov chain over the model space
has a global tendency of moving towards the true data-generating
model, even though the posterior distribution can be highly
irregular. In order to bound the mixing time, we make use of the
canonical path technique developed by
Sinclair~\cite{Sinclair1988,Sinclair1992} and Diaconis and Stroock
\cite{Diaconis1991}. More precisely, the particular canonical path
construction used in our proof is motivated by examining the solution
path of stepwise regression procedures for linear model selection
(e.g., ~\cite{Zhang11,An08}), where a greedy criterion is used to
decide at each step whether a covariate is to be included or deleted
from the curent model.

Overall, our results reveal that there is a delicate interplay between
the statistical and computational properties of Bayesian models for
variable selection.  On the one hand, we show that concentration of
the posterior is not only useful in guaranteeing desirable statistical
properties such as parameter estimation or model-selection
consistency, but they also have algorithmic benefits in certifying the
rapid mixing of the Markov chain methods designed to draw samples from
the posterior.  On the other hand, we show that posterior consistency
on its own is \emph{not} sufficient for rapid mixing, so that
algorithmic efficiency requires somewhat stronger conditions.

The remainder of this paper is organized as follows.
Section~\ref{SecBackground} provides background on the Bayesian
approach to variable selection, as well as Markov chain algorithms for
sampling and techniques for analysis of mixing times.  In
Section~\ref{SecMain}, we state our two main results
(Theorems~\ref{ThmBVSconsistency} and~\ref{ThmMain}) for a class of
Bayesian models for variable selection, along with simulations that
illustrate the predictions of our theory.  Section~\ref{SecProofMain}
is devoted to the proofs of our results, with many of the technical
details deferred to the appendices.  We conclude in
Section~\ref{SecDiscussion} with a discussion.

%%%%%%%%%%%%%%%%%%%%%%%%%%%%%%%%%%%%%%%%%%%%%%%%%%%%%%%%%%%%%%%%%%%%%%%%%%%%

\section{Background and problem formulation}
\label{SecBackground}
In this section, we introduce some background on the Bayesian approach
to variable selection, as well some background on Markov chain
algorithms for sampling, and techniques for analyzing their mixing
times.

\subsection{Variable selection in the Bayesian setting}
\label{SectionBVS}

Consider a response vector $\yobs \in \real^\numobs$ and a design
matrix $\design \in \real^{\numobs \times \covdim}$ that are linked by
the standard linear model
\begin{align}
\label{EqnStandardLinear}
\yobs = \design \betastar + \wnoise,\quad \mbox{where $\wnoise \sim
  \mathcal{N}(0,\sigma^2 I_\numobs)$,}
\end{align}
and $\betastar \in \real^\covdim$ is the unknown regression vector.
Based on observing the pair $(\yobs, \design)$, our goal is to recover
the support set of $\betastar$---that is, to select the subset of
covariates with non-zero regression weights, or more generally, a
subset of covariates with absolute regression weights above some
threshold.

In generic terms, a Bayesian approach to variable selection is based
on first imposing a prior over the set of binary indicator vectors,
and then using the induced posterior (denoted by $\pi( \gamma \, \mid
\, Y)$) to perform variable selection. Here each binary vector $\gamma
\in \{0,1\}^\covdim$ should be thought of as indexing the model which
involves only the covariates indexed by $\gamma$.  We make use of the
shorthand $|\gamma| = \sum_{j=1}^\covdim \gamma_j$ corresponding to
the number of non-zero entries in $\gamma$, or the number of active
covariates in the associated model.  It will be convenient to adopt a
dualistic view of $\gamma$ as both a binary indicator vector, and as a
subset of $\{1,\ldots, \covdim\}$.  Under this identification, the
expression $\gamma\subset \gamma'$ for a pair of inclusion vectors
$(\gamma,\gamma')$ can be understood as that the subset of variables
selected by $\gamma$ is contained in the subset of variables selected
by $\gamma'$. Similarly, it will be legitimate to use set operators on
those indicator vectors, such as $\gamma \cap \gamma'$, $\gamma \cup
\gamma'$ and $\gamma \setminus \gamma'$.  Using the set
interpretation, we let $\design_{\gamma} \in \real^{\numobs \times
  |\gamma|}$ denote the submatrix formed of the columns indexed by
$\gamma$, and we define the subvector $\beta_\gamma \in
\real^{|\gamma|}$ in an analogous manner.  We make use of this
notation in defining the specific hierarchical Bayesian model analyzed
in this paper, defined precisely in Section~\ref{SecHierarchical} to
follow.

%%%%%%%%%%%%%%%%%%%%%%%%%%%%%%%%%%%%%%%%%%%%%%%%%%%%%%%%%%%%%%%%%%%%%%%%%%

\subsection{MCMC algorithms for Bayesian variable selection}

Past work on MCMC algorithms for Bayesian variable selection can be
divided into two main classes---Gibbs samplers
(e.g.,~\cite{George1993,Ishwaran2005,Narisetty2014}) and
Metropolis-Hastings random walks (e.g.~\cite{Hans2007,Guan2011}).  In
this paper, we focus on a particular form of Metropolis-Hastings
updates.

In general terms, a Metropolis-Hastings random walk is an iterative
and local-move based procedure involving three steps:
\begin{description}
\item[Step 1:] Use the current state $\gamma$ to define a neighborhood
  $\neigh(\gamma)$ of proposal states.
\item[Step 2:] Choose a proposal state $\gamma'$ in $\neigh(\gamma)$
  according to some probability distribution $\SMAT(\gamma,\cdot)$
  over the neighborhood, e.g. the uniform distribution.
\item[Step 3:] Move to the new state $\gamma'$ with probability
  $\ACCEPT(\gamma, \gamma')$, and stay in the original state $\gamma$
  with probability \mbox{$1 - \ACCEPT(\gamma, \gamma')$}, where the
  acceptance ratio is given by
\begin{align}
\label{EqnAcceptance}
\ACCEPT(\gamma, \gamma') & \defn \min\big\{1,\frac{\posterior(\gamma'
  \mid \yobs)\, \SMAT(\gamma', \gamma)}{\posterior( \gamma \mid
  \yobs)\,\SMAT(\gamma, \gamma')} \big\}.
\end{align}
\end{description}
In this way, for any fixed choice of the neighborhood structure
$\neigh(\gamma)$, we obtain a Markov chain with transition probability
given by
\begin{align*}
\ProbMH(\gamma,\gamma') & =
\begin{cases} 
\SMAT(\gamma, \gamma')\, \ACCEPT(\gamma, \gamma') & \mbox{if $\gamma'
  \in \neigh(\gamma)$},\\
0 & \mbox{if $\gamma' \notin \neigh(\gamma) \cup \{\gamma\}$, and} \\
1-\sum_{\gamma'\neq \gamma} \ProbMH(\gamma, \gamma') & \mbox{if
  $\gamma'=\gamma$.}
\end{cases}
\end{align*}

The specific form of Metropolis-Hastings update analyzed in this paper
is obtained by randomly selecting one of the following two schemes to
update $\gamma$, each with probability $0.5$.
\begin{description}
\item[Single flip update:] Choose an index $j \in [\covdim]$ uniformly
  at random, and form the new state $\gamma'$ by setting $\gamma'_j =
  1 - \gamma_j$.
\item[Double flip update:] Define the subsets $S(\gamma) = \{ j \in
  [\usedim] \, \mid \, \gamma_j = 1 \}$ and let $\Sset^c(\gamma) = \{j
  \in [\usedim] \, \mid \gamma_j = 0 \}$.  Choose an index pair
  $(k,\ell) \in \Sset(\gamma) \times \Sset^c(\gamma)$ uniformly at
  random, and form the new state $\gamma'$ by flipping $\gamma_k$ from
  $1$ to $0$ and $\gamma_\ell$ from $0$ to $1$. (If the set
  $\Sset(\gamma)$ is empty, then we do nothing.)
\end{description}
This scheme can be understood as a particular of the general
Metropolis-Hastings scheme in terms of a neighborhood $\neigh(\gamma)$
to be all models $\gamma'$ that can be obtained from $\gamma$ by
either changing one component to its opposite (i.e., from $0$ to $1$,
or from $1$ to $0$) or switching the values of two components with
different values.  

Letting $\hamm(\gamma,\gamma')=\sum_{j=1}^\covdim \Ind(\gamma_j \neq
\gamma_j')$ denote the Hamming distance between $\gamma$ and
$\gamma'$.  the overall neighborhood is given by the union
$\neigh(\gamma) \defn \neigh_1(\gamma) \cup \neigh_2(\gamma)$, where
\begin{align*}
\neigh_1(\gamma) & \defn \big \{ \gamma' \, \mid \hamm(\gamma',\gamma)
=1 \big \}, \qquad \mbox{and} \\
\neigh_2(\gamma) & \defn \big \{\gamma' \, \mid \, \hamm(\gamma',
\gamma) = 2, \mbox{ and $\exists (k, \ell) \in S(\gamma) \times
  \Sset^c(\gamma)$ s.t. $\gamma'_k = 1 -\gamma_k$ and $\gamma'_\ell =
  1- \gamma_\ell$} \big \}.
\end{align*}

With these definitions, the transition matrix of the previously
described Metropolis-Hastings scheme takes the form
\begin{align}
\label{EqnMetropolisHastings}
\ProbMH(\gamma,\gamma') = \begin{cases}
\frac{1}{2 \, \covdim}\, \min \big \{ 1, \frac{\posterior(\gamma' \mid
  \yobs)}{\posterior(\gamma \mid \yobs)}\big\}, & \mbox{if $\gamma'
  \in \neigh_1(\gamma)$} \\
\frac{1}{2\, |\Sset(\gamma)| \, |\Sset^c(\gamma)|} \,\min \big \{1,
\frac{\posterior(\gamma' \mid \yobs)}{\posterior(\gamma \mid
  \yobs)}\big\},
& \mbox{if $\gamma' \in \neigh_2(\gamma)$} \\
0 & \mbox{if $\hamm(\gamma', \gamma) > 2$, and} \\
1-\sum_{\gamma'\neq \gamma} \ProbMH(\gamma,\gamma'), & \mbox{if
  $\gamma'=\gamma$. }
\end{cases}
\end{align}

%%%%%%%%%%%%%%%%%%%%%%%%%%%%%%%%%%%%%%%%%%%%%%%%%%%%%%%%%%%%%%%%%%%%%%%%%
\subsection{Background on mixing times}

Let $\MarkovChain$ be an irreducible, aperiodic Markov chain on the
discrete state space $\Mspace$, and described by the transition
probability matrix \mbox{$\PMAT \in \real^{|\Mspace|\times
    |\Mspace|}$} with stationary distribution $\pi$.  We assume
throughout that $\MarkovChain$ is reversible; i.e., it satisfies the
detailed balance condition \mbox{$\pi(\gamma) \PMAT(\gamma,\gamma') =
  \pi(\gamma') \PMAT(\gamma',\gamma)$} for all
$\gamma,\gamma'\in\Mspace$. It is easy to see that the previously
described Metropolis-Hastings matrix $\ProbMH$ satisfies this
reversibility condition. It is convenient to identify a reversible
chain with a weighted undirected graph $G$ on the vertex set
$\Mspace$, where two vertices $\gamma$ and $\gamma'$ are connected if
and only if the edge weight $\QMAT(\gamma,\gamma') \defn \pi(\gamma)
\PMAT(\gamma,\gamma')$ is strictly positive.

For $\gamma \in \Mspace$ and any subset $S \subseteq \Mspace$, we
write $\PMAT(\gamma,S)=\sum_{\gamma'\in S}
\PMAT(\gamma,\gamma')$. If $\gamma$ is the initial state of the
chain, then the total variation distance to the stationary
distribution after $t$ iterations is
\begin{align*}
\Delta_\gamma(t) = \|\PMAT^n(\gamma,\cdot)-\pi(\cdot)\|_{TV}\defn
\max_{S\subset\Mspace} \big|\PMAT^n(\gamma,S)-\pi(S)\big|.
\end{align*}
The $\epsilon$-mixing time is given by
\begin{align}
\label{EqnDefMixingTime}
\tau_\epsilon \defn \max_{\gamma\in\Mspace} \min \big \{t \in
\mathbb{N} \, \mid \Delta_\gamma(t')\leq\epsilon\mbox{ for all }t'\geq
t \big \},
\end{align}
which measures the number of iterations required for the chain to be
within distance $\epsilon\in(0,1)$ of stationarity.  The efficiency of
the Markov chain can be measured by the dependence of
$\tau_{\epsilon}$ on the difficulty of the problem, for example, the
dimension of the parameter space and the sample size. In our case, we
are interesed in the dependence of $\tau_{\epsilon}$ on the covariate
dimension $\covdim$ and the sample size $\numobs$.  Of particular
interest is whether the chain is \emph{rapidly mixing}, meaning that
the mixing time grows at most polynomially in the pair $(\covdim,
\numobs)$, or \emph{slowly mixing,} meaning that the mixing time grows
exponentially.

%%%%%%%%%%%%%%%%%%%%%%%%%%%%%%%%%%%%%%%%%%%%%%%%%%%%%%%%%%

\section{Main results and their consequences}
\label{SecMain}

The analysis of this paper applies to a particular family of
hierarchical Bayes models for variable selection.  Accordingly, we
begin by giving a precise description of this family of models, before
turning to statements of our main results and a discussion of their
consequences.  Our first result (Theorem~\ref{ThmBVSconsistency})
provides sufficient conditions for posterior concentration, whereas
our second result (Theorem~\ref{ThmMain}) provides sufficient
conditions for rapid mixing of the Metropolis-Hastings updates.

%%%%%%%%%%%%%%%%%%%%%%%%%%%%%%%%%%%%%%%%%%%%%%%%%%%%%%%%%%%%%%%%%%%%%%

\subsection{Bayesian hierarchical model for variable selection}
\label{SecHierarchical}

In addition to the standard linear model~\eqref{EqnStandardLinear},
the Bayesian hierarchical model analyzed in this paper involves three
other ingredients: a prior over the precision parameter $\phi$ (or
inverse noise variance) in the linear observation model, a prior on
the regression coefficients, and a prior over the binary indicator
vectors.  More precisely, it is given by
\begin{subequations}
\label{EqnHierarchicalBayes}
\begin{align}
\MODEL_{\gamma}:\qquad\qquad\qquad & \mbox{Linear model:} \qquad
\yobs=\design_{\gamma}\beta_{\gamma} + \wnoise,\quad \wnoise\sim
\mathcal{N} (0, \phi^{-1}I_\numobs) \\
& \mbox{Precision prior} \qquad \plprior(\phi) \propto
\frac{1}{\phi} \\
\label{EqnGPrior}
& \mbox{Regression prior} \qquad \big(\beta_{\gamma} \, \mid \, \gamma
\big) \sim \mathcal{N} (0, \hyperpara \, \phi^{-1}
(\design_{\gamma}^T\design_{\gamma})^{-1}) \\
\label{EqnPriorGamma}
& \mbox{Sparsity prior} \qquad \plprior(\gamma) \propto
\Big(\frac{1}{\covdim} \Big)^{\kappa |\gamma|} \Ind [|\gamma|\leq \maxsize].
\end{align}
\end{subequations}
For each model $\MODEL_\gamma$, there are three parameters to be
specified: the integer $\maxsize<\numobs$ is a prespecified upper
bound on the maximum number of important covariates, the
hyperparameter $\hyperpara > 0$ controls the degree of dispersion in
the regression prior, and the hyperparameter $\kappa>0$ penalizes
models with large size.  For a given integer $\maxsize \in \{1,
\ldots, \covdim \}$, we let $\MspaceMax = \{ \MODEL_\gamma \, \mid \,
|\gamma| \leq \maxsize \}$ the class of all models involving at most
$\maxsize$ covariates.

Let us make a few remarks on our choice of Bayesian model.  First, the
choice of covariance matrix in the regression prior---namely,
involving $X_\gamma^T X_\gamma$---is made for analytical convenience,
in particular in simplifying the posterior.  A more realistic choice
would be the independent prior
\begin{align*}
\beta_{\gamma}\, \mid \, \gamma \sim \mathcal{N} (0, \hyperpara \,
\phi^{-1} I_{|\gamma|}).
\end{align*} 
However, the difference between these choices will be negligible when
$g \gg \numobs$, which, as shown by our theoretical analysis, is the
regime under which the posterior is well-behaved.  Another popular
choice for the prior of $\beta_{\gamma}$ is the spike-and-slab
prior~\cite{Ishwaran2005}, where for each covariate $\design_j$, one
specifies the marginal prior for $\beta_j$ as a mixture of two normal
distributions, one with a substantially larger variance than the
other, and $\gamma_j$ can be viewed as the latent class indicator for
this mixture prior.  Our primary motivation in imposing Zellner's
$\hyperpara$-prior is in order to streamline the theoretical analysis:
it leads to an especially simple form of the marginal likelihood
function.  However, we note that our conclusions remain valid under
essentially the same conditions when the independent prior or the
spike-and-slab prior is used, but with much longer proofs.  The
sparsity prior on $\gamma$ is similar to the prior considered by
Narisetty and He~\cite{Narisetty2014} and Castillo et
al. \cite{Castillo2015}.  The $\covdim^{-\kappa}$ decay rate for the
marginal probability of including each covariate imposes a vanishing
prior probability on the models of diverging sizes. The only
difference is that we put a constraint $|\gamma| \leq \maxsize$ to
rule out models with too many covariates. As will be clarified in the
sequel, while this additional constraint is not needed for Bayesian
variable-selection consistency, it is necessary for rapid mixing of
the MCMC algorithm that we analyze.

Recall from our earlier set-up that the response vector $\yobs \in
\real^\numobs$ is generated from the standard linear model
\mbox{$\yobs = \design \betastar + \wnoise$,} where $\wnoise \sim
\neigh(0,\sigmazero^2 I_\numobs)$, $\betastar \in \real^\covdim$ is
the unknown regression vector, and $\sigmazero$ the unknown noise
standard deviation.  In rough terms, the goal of variable selection is
to determine the subset $S$ of ``influential'' covariates.  In order
to formalize this notion, let us fix a constant $\CB > 0$ depending on
$(\sigmazero, \numobs, \covdim)$ that quantifies the minimal signal
size requirement for a covariate to be ``influential''.  We then
define $\Sset = \Sset(\CB)$ to consist of the indices with relatively
large signal---namely
\begin{align}
\label{EqnDefnCB}
\Sset \defn \big \{ j \in [\covdim] \mid \, |\betastar_j| \geq \CB
\big\},
\end{align}
and our goal is to recover this subset.  Thus, the ``non-influential''
coefficients $\betastar_{\Sset^c}$ are allowed to be non-zero, but their
magnitudes are constrained.

We let $\gammastar$ be the indicator vector that selects the
influential covariates, and let $\sstar\defn |\gammastar|$ be the size
of the corresponding ``true'' model $\MODEL_{\gammastar}$.  Without
loss of generality, we may assume that the first $\sstar$ components
of $\gammastar$ are ones, and the rest are zeros.  We assume
throughout this section that we are in the high-dimensional regime
where $\covdim\geq \numobs$, since the low-dimensional regime where
$\numobs < \covdim$ is easier to analyze. For any symmetric matrix
$\QMAT$, let $\lammin(\QMAT)$ and $\lammax(\QMAT)$ denote its smallest
and largest eigenvalues.  Our analysis involves the following
assumptions:

\paragraph{Assumption A (Conditions on $\betastar$):}
The true regression vector has components \mbox{$\betastar =
  (\betastar_\Sset, \betastar_{\Sset^c})$} that satisfy the bounds
\begin{subequations}
\begin{align}
\begin{aligned}
{\mathbf{\mbox{Full $\betastar$ condition:}}} \qquad &
\big\|\frac{1}{\sqrt{\numobs}} \Xmat \betastar \big\|_2^2 \leq
\hyperpara \, \sigmazero^2 \, \frac{\log \covdim}{\numobs} \\
{\mathbf{\mbox{Off-support $\Sset^c$ condition:}}} \qquad &
\big\|\frac{1}{\sqrt{\numobs}} \Xmat_{\Sset^c} \betastar_{\Sset^c}
\big\|_2^2 \leq \Ltil \, \sigmazero^2 \, \frac{\log \covdim}{\numobs},
\end{aligned}
\end{align}
for some universal constant $\Ltil$.

In the simplest case, the true regression vector $\betastar$ is
$S$-sparse (meaning that $\beta^*_{S^c} = 0$), so that the off-support
condition holds trivially.  As for the full $\betastar$ condition, it
is known~\cite{Shang2011} that some form of upper bound on the norm
$\|\betastar\|_2$ in terms of the $g$-hyperparameter is required in
order to prove Bayesian model selection
consistency~\cite{Shang2011}. The necessity of such a condition is a
manifestation of the so-called information paradox of
$g$-priors~\cite{Liang2008}. \\

Our next assumption involves an integer parameter $\spindex$, which is
set either to a multiple of the true sparsity $\sstar$ (in order to
prove posterior concentration) or the truncated sparsity $\maxsize$
(in order to prove rapid mixing).

\paragraph{Assumption B (Conditions on the design matrix):} 
The design matrix has been normalized so that $\|X_j\|_2^2 = \numobs$
for all $j = 1, \ldots, \covdim$; moreover, letting $Z \sim N(0,
I_\numobs)$, there exist constants $\Cm > 0$ and $L < \infty$ such
that
\begin{align}
\begin{aligned}
{\mathbf{\mbox{Lower restricted eigenvalue ($\RE(\spindex)$):}}}
\qquad & \min_{|\gamma| \leq \spindex}
\lammin\Big(\frac{1}{\numobs}\design_{\gamma}^T\design_{\gamma}\Big)
\, \geq \, \Cm, \quad \mbox{and} \\
{\mathbf{\mbox{Sparse projection condition (\IC($\spindex$)):}}}
\qquad & \Exs_Z \Big[\max_{|\gamma|\leq \spindex} \max_{k \in
    [\covdim] \backslash \gamma} \frac{1}{\sqrt{\numobs}}
  \big|\inprod{ \big(I - \Proj_{\gamma}\big)\design_{k}}{Z}\big|\Big]
\leq \frac{1}{2} \sqrt{L \Cm \log \covdim},
\end{aligned}
\end{align}
where $\Proj_\gamma$ denotes projection onto the span of $\{X_j, j \in
\gamma \}$.  The lower restricted eigenvalue condition is a mild
requirement, and one that plays a role in the information-theoretic
limitations of variable selection~\cite{Martin2009}.  On the other
hand, the \HACKER condition can always be satisfied by choosing
$L=\order(\maxsize)$.  To see this, notice that
$\frac{1}{\sqrt{n}}\|(I-\Proj_{\gamma})X_k\| \leq 1$ and there are at
most $\covdim^{\maxsize}$ different choice of distinct pair
$(\gamma,k)$. Therefore, by the Gaussianity of $g_G$, the \HACKER
condition always holds with $L = 4\Cm^{-1} \maxsize$. On the other
extreme, if the design matrix $X$ has orthogonal columns, then $\big(I
- \Proj_{\gamma}\big)\design_{k}=\design_{k}$. As a consequence, due
to the same argument, the \HACKER condition holds with $L=4\Cm^{-1}$,
which depends neither on $\sstar$ nor on $\maxsize$.

%%%%%%%%%%%%%%%%%%%%%

\paragraph{Assumption C (Choices of prior hyperparameters):}
The noise hyperparameter $\hyperpara$ and sparsity penalty
hyperparameter $\kappa > 0$ are chosen such that
\begin{align}
\begin{aligned}
\hyperpara \asymp \covdim^{2\alpha} & \qquad \mbox{for some $\alpha >
  0$, and} \\
\kappa + \alpha \geq \calpha (L + \Ltil) + 2& \qquad \mbox{for some
  universal constant $\calpha > 0$.}
\end{aligned}
\end{align}
\noindent In the low-dimensional regime $\covdim = o(\numobs)$, the
$\hyperpara$-prior with either the unit information prior $\hyperpara
= \numobs$, or the choice $\hyperpara = \max\{\numobs, \covdim^2\}$
have been recommended~\cite{Kass1995,Fernandez2001,Spark2015}. In the
intermediate regime where $\covdim = \mathcal{O}(\numobs)$, Sparks et
al.~\cite{Spark2015} show that $\hyperpara$ must grow faster than
$\covdim\, \numobs^{-1} \log \covdim$ for the Bayesian linear model
without variable selection to achieve posterior consistency. These
considerations motivate us to choose the hyperparameter for the
high-dimensional setting as $\hyperpara \asymp \covdim^{2\alpha}$ for
some $\alpha > 0$, and our theory establishes the utility of this choice.

%%%%%%%%%%%%%%%%%%%%%%%

\paragraph{Assumption D (Sparsity control):}  
For a constant $\rcon > 4$, one of the two following conditions holds:
\begin{description}
\item[$\quad$ Version D$(\sstar)$:] We set $\maxsize \defn \usedim$
  in the sparsity prior~\eqref{EqnPriorGamma}, and the true sparsity
  $\sstar$ is bounded as $\sstar \; \leq \; \frac{1}{8 \rcon K} \Big
  \{ \frac{\numobs}{\log \usedim} - 16 \Ltil \sigmazero^2 \Big \}$ for
  some constant $K \geq 4+ \alpha + \UNICON \Ltil$.

\item[$\quad$ Version D$(\maxsize)$:] The sparsity parameter
  $\maxsize$ in the prior~\eqref{EqnPriorGamma} satisfies the sandwich
  relation
\begin{align}
\label{EqnSandwich}
\big(2\Cm^{-2} \,\omega(\Xmat) + 1\big) \sstar \; \leq \; \maxsize \;
\leq \; \frac{1}{8 \rcon \, K} \Big \{ \frac{\numobs}{\log \usedim} -
16 \Ltil \sigmazero^2 \Big \},
\end{align}
where $\omega(X) \defn \max \limits_{\gamma\in \Mspace}
\opnorm{(X_{\gamma}^TX_{\gamma})^{-1} X_{\gamma}^T X_{\gammastar
    \setminus \gamma}}^2$.
\end{description}
\end{subequations}

\vspace*{.1in}

Assumptions A, B, C and D are a common set of conditions assumed in the
existing literature (e.g., \cite{Shang2011,Narisetty2014}) for
establishing Bayesian variable-selection consistency; i.e., 
that the posterior probability of the true model
$\posterior(\gammastar|\yobs)\to 1$ as $\numobs\to\infty$.

%%%%%%%%%%%%%%%%%%%%%%%%%%%%%%%%%%%%%%%%%%%%%%%%%%%%%%%%%%%%%%%%%%%%%%%%%

\subsection{Sufficient conditions for posterior consistency}

Our first result characterizes the behavior of the (random) posterior
$\posterior(\cdot \mid Y)$.  As we mentioned in
Section~\ref{SectionBVS}, Bayesian variable-selection consistency does
not require that the sparsity prior~\eqref{EqnPriorGamma} be truncated
at some sparsity level much less than $\usedim$, so that we analyze
the hierarchical model with $\maxsize = \usedim$, and use the milder
Assumption D$(\sstar)$.  The reader should recall from
equation~\eqref{EqnDefnCB} the threshold parameter $\CB$ that defines
the subset $S = S(\CB)$ of influential covariates.

\btheos[Posterior concentration]
\label{ThmBVSconsistency}
Suppose that Assumption A, Assumption B with \mbox{$s = K \sstar$,}
Assumption C, and Assumption D$(\sstar)$ hold. If the threshold $\CB$ satisfies
\begin{align}
\label{EqnBetamin}
\CB^2 \geq \UNICON_0( L + \Ltil + \alpha
+ \kappa)\, \sigmazero^2 \, \frac{\log \covdim}{\numobs},
\end{align}
then we have $\posterior(\gammastar \mid \yobs) \geq 1 - \UNICON_1\,
  \covdim^{-1}$ with probability
at least $1 - \UNICON_2 \, \covdim^{-\UNICON_3}$.
\etheos

The threshold condition~\eqref{EqnBetamin} requires the set of
influential covariates to have reasonably large magnitudes; this type
of signal-to-noise condition is needed for establishing variable
selection consistency of any procedure~\cite{Martin2009}.  We refer to
it as the $\betamin$-condition in the rest of the paper.  Due to the
mildness of Assumption A (conditions on $\betastar$), the claim in the
theorem holds even when the true model is not exactly sparse:
Assumption A allows the residual $\betastar_{\Sset^c}$ to be nonzero
as long as it has small magnitude.

It is worth noting that the result of Theorem~\ref{ThmBVSconsistency}
covers two regimes, corresponding to different levels of
signal-to-noise ratio.  More precisely, it is useful to isolate the
following two mutually exclusive possibilities:
\begin{subequations}
\begin{align}
\label{EqnHighSNR}
\mbox{High SNR:} & \qquad \Sset = \big\{j \in [\covdim] \, \mid \,
\betastar_j \neq 0 \big \} \quad \mbox{and}\quad \min_{j \in \Sset}
|\betastar_j|^2 \geq \UNICON_0(\alpha + \kappa + L)\, \sigmazero^2 \,
\frac{\log \covdim}{\numobs}, \\
\label{EqnLowSNR}
\mbox{Low SNR:} & \qquad \Sset = \emptyset \quad \mbox{and} \quad  \big\|\frac{1}{\sqrt{\numobs}} \Xmat \betastar
\big\|_2^2  \leq
\Big(\frac{\alpha + \kappa - 2}{\calpha} - L \Big) \,\sigmazero^2\frac{\log\covdim}{\numobs}.
\end{align}
\end{subequations}
In terms of the parameter $\Ltil$ in Assumption A, The high SNR regime
corresponds to $\Ltil = 0$, whereas the low SNR regime corresponds to
$\Ltil = \frac{\alpha + \kappa - 2}{\calpha} - L$.  The intuition for
the low SNR setting is that the signal in every component is so weak
that the ``penalty'' induced by hyperparameters $(\hyperpara,\kappa)$
completely overwhelms it.  Theorem~\ref{ThmBVSconsistency} guarantees
that the posterior concentrates around the model $\MODEL_\gammastar$
under the high SNR condition, and under the null model $\MODEL_{\gamma_0}$
under the low SNR condition.  More precisely, we have:
\bcors
\label{CoroBVSconsistency}
Under the conditions of Theorem~\ref{ThmBVSconsistency}, with
probability at least $1 - \UNICON_2 \, \covdim^{-\UNICON_3}$:
\begin{enumerate}
\item[(a)] Under the high SNR condition~\eqref{EqnHighSNR}, we
  have $\posterior(\gammastar \mid \yobs) \geq 1 - \UNICON_1\,
  \covdim^{-1}$.
\item[(b)] Conversely, under the low SNR condition~\eqref{EqnLowSNR},
  we have $\posterior(\gamma_{0} \, \mid \yobs) \geq 1 - \UNICON_1\,
  \covdim^{-1}$. 
\end{enumerate}
\ecors

Corollary~\ref{CoroBVSconsistency} provides a complete
characterization of the high or low SNR regimes, but it does not cover
the intermediate regime, in which some component $\betastar_j$ of
$\betastar$ is sandwiched as
\begin{align}
\label{EqnIntermediateSNR}
\Big(\frac{\alpha + \kappa - 2}{\calpha} - L
\Big)\,\sigmazero^2\frac{\log\covdim}{\numobs} \; \leq|\betastar_j|^2
\; \leq \UNICON_0( \alpha + \kappa + L)\, \sigmazero^2 \frac{\log
  \covdim}{\numobs}.
\end{align} 
On one hand, Theorem~\ref{ThmBVSconsistency} still guarantees a form
of Bayesian variable selection consistency in this regime.  However,
the MCMC algorithm for sampling from the posterior can exhibit slow
mixing due to multimodality in the posterior.  In
Appendix~\ref{AppSlowMixing}, we provide a simple example that
satisfies the conditions of Theorem~\ref{ThmBVSconsistency}, so that
posterior consistency holds, but the Metropolis-Hastings updates have
mixing time growing exponentially in $\usedim$.  This example reveals
a phenomenon that might seem counter-intuitive at first sight: sharp
concentration of the posterior distribution need not lead to rapid
mixing of the MCMC algorithm.

%%%%%%%%%%%%%%%%%%%%%%%%%%%%%%%%%%%%%%%%%%%%%%%%%%%%%%%%%%%%%%%%%%%%%%%%%

\subsection{Sufficient conditions for rapid mixing}

With this distinction in mind, we now turn to developing sufficient
conditions for Metropolis-Hastings
scheme~\eqref{EqnMetropolisHastings} to be rapidly mixing.  As
discussed in Section~\ref{SecBackground}, this rapid mixing ensures
that the number of iterations required to converge to an
$\epsilon$-ball of the stationary distribution grows only polynomially
in the problem parameters.  The main difference in the conditions is
that we now require Assumption B---the RE and \HACKER conditions---to
hold with parameter $s = \maxsize$, as opposed to with the smaller
parameter $s = K \sstar \ll \maxsize$ involved in
Theorem~\ref{ThmBVSconsistency}.
\btheos[Rapid mixing guarantee]
\label{ThmMain}
Suppose that Assumption A, Assumption B with $s = \maxsize$,
Assumption C, and Assumption D($\maxsize$) all hold.  Then under
either the high SNR condition~\eqref{EqnHighSNR} or the low SNR
condition~\eqref{EqnLowSNR}, there are universal constants $c_1, c_2$
such that, for any $\epsilon \in (0,1)$, the $\epsilon$-mixing time of
the Metropolis-Hastings chain~\eqref{EqnMetropolisHastings} is upper
bounded as
\begin{align}
\label{EqnMixingBound}
\tau_{\epsilon} \leq c_1 \, \covdim \maxsize^2 \, \big (c_2 \alpha \,
( \numobs + \maxsize) \log \covdim + \log(1/\epsilon) + 2 \big)
\end{align}
with probability at least $1-4\covdim^{-c_1}$.
\etheos

According to our previous definition~\eqref{EqnDefMixingTime} of the
mixing time, Theorem~\ref{ThmMain} characterizes the worst case mixing
time, meaning the number of iterations when starting from the worst
possible initialization. If we start with a good intial state---for
example, the true model $\gammastar$ would be a nice though
impractical choice---then we can remove the $\numobs$ term in the
upper bound~\eqref{EqnMixingBound}.  In this way, the term
$\UNICON_1\UNICON_2 \alpha\, \numobs \covdim \maxsize^2 \log \covdim$
can be understood as the worst-case number of iterations required in
the burn-in period of the MCMC algorithm.

Theorem~\ref{ThmBVSconsistency} and Theorem~\ref{ThmMain} lead to the
following corollary, stating that after $\mathcal{O}(\alpha\, \numobs
\covdim \maxsize^2 \log \covdim)$ iterations, the MCMC algorithm will
output $\gammastar$ with high probability.

\bcors
\label{CoroMCMC}
Under the conditions of Theorem~\ref{ThmMain}, for any fixed iterate
$t$ such that
\begin{align*}
t \geq c_1 \, \covdim \maxsize^2 \, \big (c_2 \alpha \, ( \numobs +
\maxsize) \log \covdim + \log \covdim + 2 \big),
\end{align*}
the iterate $\gamma_t$ from the MCMC algorithm matches $\gammastar$
with probability at least $1 - \UNICON_2 \, \covdim^{-\UNICON_3}$.
\ecors

As with Corollary~\ref{CoroBVSconsistency}, Theorem~\ref{ThmMain} does
not characterize the intermediate regime in which some component
$\betastar_j$ of $\betastar$ satisfies the sandwich
inequality~\eqref{EqnIntermediateSNR}.  Based on our simulations, we
suspect that the Markov chain might be slowly mixing in this regime,
but we do not have a proof of this statement.

%%%%%%%%%%%%%%%%%%%%%%%%%%%%%%%%%%%%%%%%%%%%%%%%%%%%%%%%%%%%%%%%%%%%%%%%%%

\subsection{Illustrative simulations}

In order to illustrate the predictions of Theorem~\ref{ThmMain}, we
conducted some simulations.  We also provide an example for which a
frequentist method such as the Lasso fails to perform correct variable
selection while our Bayesian method succeeds.

%%%%%%%%%%%%%%%%%%%%%%%%%%%%%%%%%%%%%%%%%%%%%%%%%%%%%%%%%%%%%%%%%%%%%%%%

\subsubsection{Comparison of mixing times}

In order to study mixing times and their dependence on the model
structure, we performed simulations for linear models with random
design matrices, formed by choosing row $\rowdesign_i \in
\real^\covdim$ i.i.d. from a multivariate Gaussian distribution.  In
detail, setting the noise variance $\sigma^2 = 1$, we considered two
classes of linear models with random design matrices $X \in
\real^{\numobs \times \covdim}$, in each case formed with i.i.d. rows
$\rowdesign_i \in \real^\covdim$:
\begin{align*}
& \mbox{Independent design:} & \yobs \sim \neigh(\design\betastar,
  \sigma^2I_\numobs) & \mbox{ with $\rowdesign_i \sim \neigh(0,
    I_\covdim)$ i.i.d.}; \\
& \mbox{Correlated design:} & \yobs \sim \neigh(\design\betastar,
  \sigma^2 I_\numobs) & \mbox{ with $\rowdesign_i \sim \neigh(0,
    \Sigma)$ i.i.d.  and $\Sigma_{jk}=e^{-|j-k|}$}.
\end{align*}
In all cases, we choose a design vector $\betastar \in \real^\covdim$
with true sparsity $\sstar = 10$, taking the form
\begin{align*}
\betastar = \SNR \,\sqrt{\frac{\sigma^2\log \covdim}{\numobs}} \big(
2,\,-3,\,2,\, 2,\, -3 ,\,3 ,\, -2 ,\, 3,\,-2 ,\,3,\, 0 ,\, \cdots,\, 0
\big)^T\in\real^{\covdim},
\end{align*}
where $\SNR > 0$ is a signal-to-noise parameter.  Varying the
parameter $\SNR$ allows us to explore the behavior of the chains when
the model lies on the boundary of the $\betamin$-condition.  We
performed simulations for the SNR parameter $\SNR \in\{0.5, 1, 2,
3\}$, sample sizes \mbox{$\numobs \in \{300,900\}$,} and number of
covariates $\covdim \in \{500,5000\}$.  In all cases, we specify our
prior model by setting the dispersion hyperparameter $\hyperpara =
\covdim^3$ and the expected maximum model size $\maxsize=100$.

\begin{figure}[h]
\begin{center}
\begin{tabular}{ccc}
\widgraph{.45\textwidth}{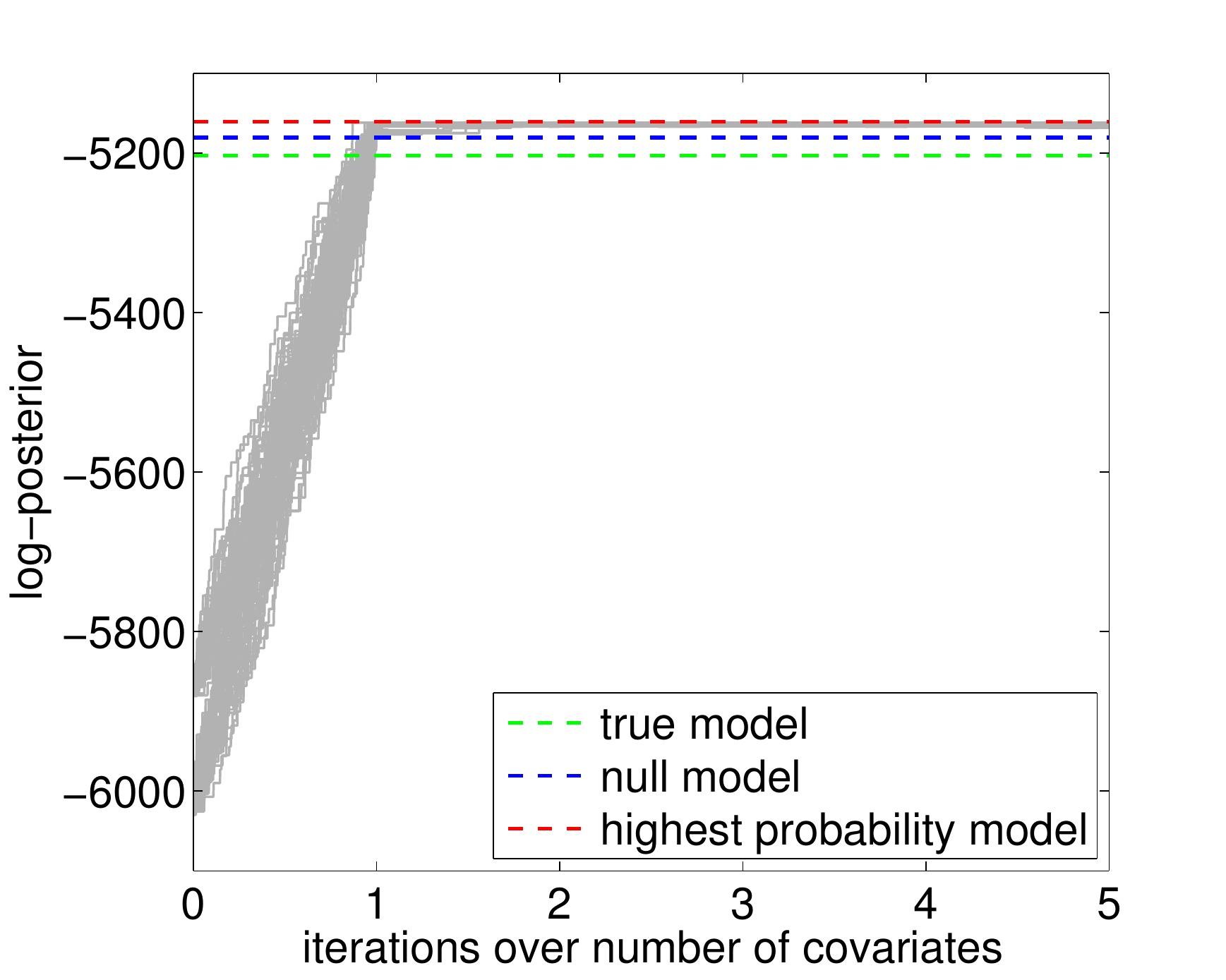} & &
\widgraph{.45\textwidth}{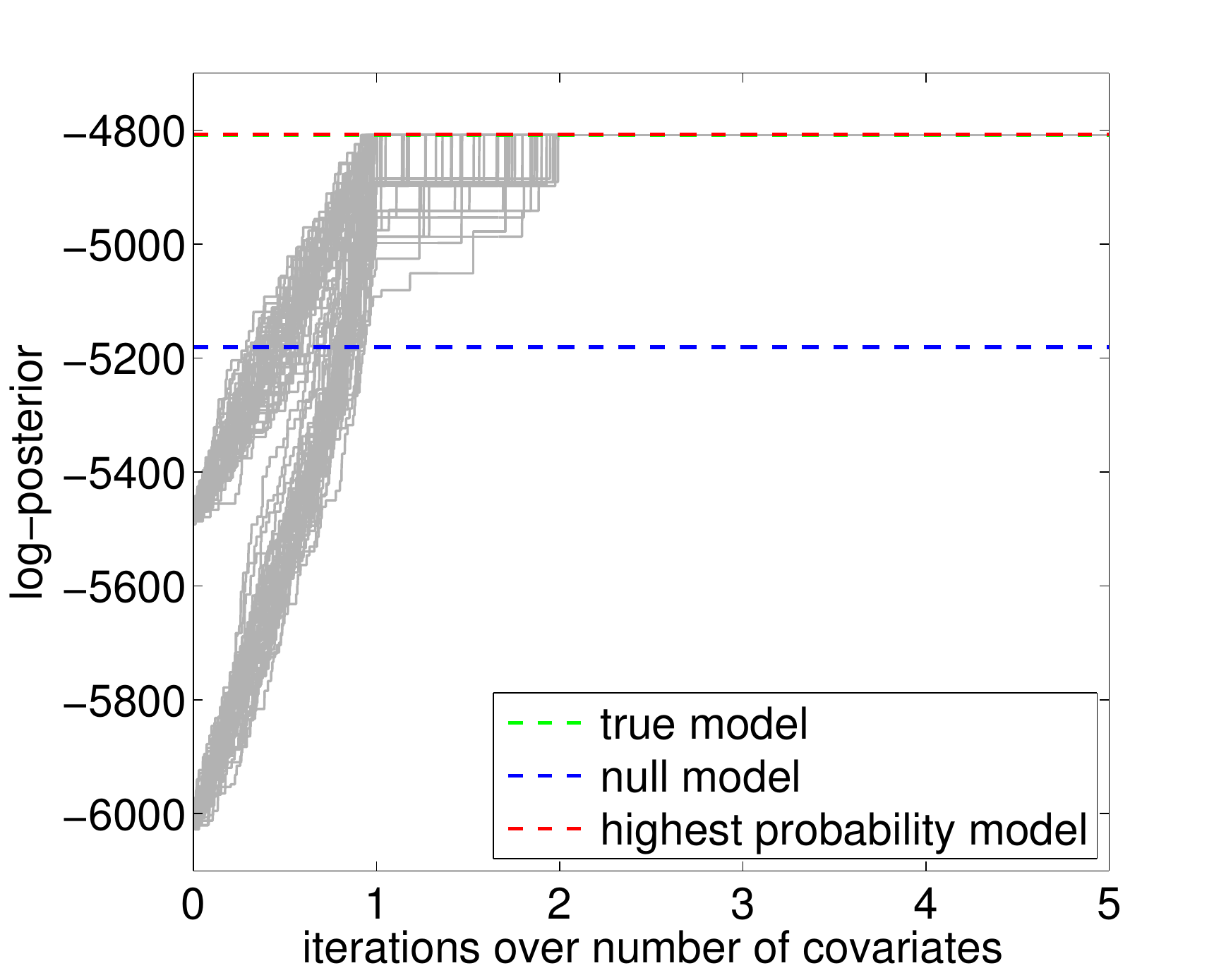} \\
(a) & & (b)
\end{tabular}
\end{center}
 \caption{Log-posterior probability versus the number of iterations
   (divided by the number of covariates $\covdim$) of $100$ randomly
   initialized Markov chains with \mbox{$\numobs = 500$,}
   \mbox{$\covdim=1000$} and \mbox{$\SNR \in \{1, 3\}$} in the
   independent design. In all cases, each grey curve corresponds to
   one trajectory of the chain ($100$ chains in total). Half of the
   chains are initialized at perturbations of the null model and half
   the true model. (a) Weak signal case: $\SNR=1$.  (b) Strong signal
   case: $\SNR=3$ (the posterior probability of the true model
   coincides with that of the highest probability model). }
\label{FigLogPosterior}
\end{figure}

Figure~\ref{FigLogPosterior} plots the typical trajectories of
log-posterior probability versus the number of iterations of the
Markov chain under the independent design. In the strong signal regime
($\SNR =3$), the true model receives the highest posterior
probability, and moreover the Metropolis-Hastings chain converges
rapidly to stationarity, typically within $3\covdim$ iterations.  This
observation is confirmation of our theoretical prediction of the
behavior when all nonzero components in $\betastar$ have relative high
signal-to-noise ratio ($\Sset = \{j: \beta_j\neq 0\}$).  In the
intermediate signal regime ($\SNR=1$), Bayesian variable-selection
consistency typically fails to hold, and here, we find that the chain
converges even more quickly to stationarity, typically within
$1.5\covdim$ iterations. This observation cannot be fully explained by
our theory.  A simulation to follow using a correlated design shows
that it is not a robust phenomenon: the chain can have poor mixing
performance in this intermediate signal regime when the design is
sufficiently correlated.

In order to gain further insight into the algorithm's performance, for
each pair $\{\design,\yobs\}$ we ran the Metropolis-Hastings random
walk based on six initializations: the first three of them are random
perturbations of the null model, whereas the remaining three are the
true model.  We made these choices of initialization because our
empirical observations suggest that the null model and the true model
tend to be near local modes of the posterior distribution. We run the
Markov chain for $20\covdim$ iterations and use the Gelman-Rubin (GR)
scale factor~\cite{Gelman1992} to detect whether the chains have
reached stationarity.  More precisely, we calculate the GR scale
factor for the coefficient of determination summary statistics
\begin{align*}
R_{\gamma}^2 = \frac{\yobs^T \Proj_{\gamma} \yobs}{\|\yobs\|_2^2},
\quad \mbox{for $\gamma\in \{0,1\}^\covdim$,}
\end{align*}
where $\Proj_\gamma$ denotes the projection matrix onto the span of
$\{X_j, j \in \gamma \}$.  Since the typical failing of convergence to
stationarity is due to the multimodality of the posterior
distribution, the GR scale factor can effectively detect the
problem. If the chains fail to converge, then the GR scale factor will
be much larger than $2$; otherwise, the scale factor should be close
to $1$. Convergence of the chain within at most $20\covdim$ iterations
provides empirical confirmation of our theoretical prediction that the
mixing time grows at most linearly in the covariate dimension
$\covdim$. (As will be seen in our empirical studies, the sample size
$\numobs$ and $\maxsize$ have little impact on the mixing time, as
long as $\maxsize$ remains small compared to $\numobs$.)

\begin{table}[t]
\centering
\begin{tabular}{c|c|C{1.8cm}C{1.8cm}C{1.8cm}C{1.8cm}}
  \hhline{======} $(\numobs,\covdim)$ & & $\SNR=0.5$ & $\SNR=1$ &
  $\SNR=2$ & $\SNR=3$ \\ \hline \multirow{3}{*}{$(500,1000)$} & SP
  & \bf{100} & \bf{100} & \bf{100} & \bf{100} \\ & H-T & 113.4 & 24.6
  & 0 & 0 \\ & N-T & 113.4 & 11.4 & -210.9 & -383.6 \\ \hline
  \multirow{3}{*}{$(500,5000)$} & SP & \bf{100} & \bf{100} & \bf{100}
  & \bf{100} \\ & H-T & 148.7 & 33.2 & 0 & 0 \\ & N-T & 148.7 & 17.4 &
  -216.6 & -395.9 \\ \hline \multirow{3}{*}{$(1000,1000)$} & SP &
  \bf{100} & \bf{100} & \bf{100} & \bf{100} \\ & H-T & 117.1 & 34.8 &
  0 & 0 \\ & N-T & 117.1 & -6.9 & -342.4 & -649.5 \\ \hline
  \multirow{3}{*}{$(1000,5000)$} & SP & \bf{100} & \bf{100} & \bf{100}
  & \bf{100} \\ & H-T & 160.4 & 32.8 & 0 & 0 \\ & N-T & 160.4 & -4.2 &
  -377.6 & -743.4 \\ \hhline{======}
\end{tabular}
\caption{Convergence behaviors of the Markov chain methods with sample
  sizes \mbox{$\numobs \in \{500,1000\}$,} ambient dimensions
  \mbox{$\covdim \in \{1000,5000\}$}, and $\SNR \in\{0.5,1,2,3\}$ in
  the independent design. SP: proportion of successful trials (in
  which GR$\leq 1.5$); H-T: log posterior probability difference
  between the highest probability model and the true model; N-T: log
  posterior probability difference between the null model and the true
  model. Each quantity is computed based on 20 simulated datasets.}
\label{table:IndptDesgin}
\end{table}

\begin{table}[h]
\centering
\begin{tabular}{c|c|C{1.8cm}C{1.8cm}C{1.8cm}C{1.8cm}}
  \hhline{======} $(\numobs,\covdim)$ & & $\SNR=0.5$ & $\SNR=1$ &
  $\SNR=2$ & $\SNR=3$ \\ \hline \multirow{3}{*}{$(500,1000)$} & SP
  & \bf{100} & \bf{95} & \bf{80} & \bf{100} \\ & H-T & 123.4 & 75.2 &
  0 & 0 \\ & N-T & 123.4 & 71.2 & -107.3 & -275.8 \\ \hline
  \multirow{3}{*}{$(500,5000)$} & SP & \bf{100} & \bf{15} & \bf{100} &
  \bf{100} \\ & H-T & 170.0 & 81.0 & 0 & 0 \\ & N-T & 170.0 & 78.7 &
  -102.1 & -288.9 \\ \hline \multirow{3}{*}{$(1000,1000)$} & SP &
  \bf{100} & \bf{100} & \bf{100} & \bf{100} \\ & H-T & 138.7 & 75.1 &
  0 & 0 \\ & N-T & 138.7 & -67.0 & -180.8 & -431.7 \\ \hline
  \multirow{3}{*}{$(1000,5000)$} & SP & \bf{100} & \bf{100} & \bf{100}
  & \bf{100} \\ & H-T & 161.8 & 61.9 & 0 & 0 \\ & N-T & 161.8 & -58.8
  & -204.2 & -445.4 \\ \hhline{======}
\end{tabular}
\caption{Convergence behavior of the Markov chain methods with sample
  size \mbox{$\numobs \in \{500, 1000\}$}, ambient dimension
  \mbox{$\covdim \in \{1000, 5000\}$}, and parameter \mbox{$\SNR \in
    \{0.5,1,2,3\}$} for the case of correlated design. SP: proportion
  of successful trials (in which GR $\leq 1.5$); H-T: log posterior
  probability difference between the highest probability model and the
  true model; N-T: log posterior probability difference between the
  null model and the true model. Each quantity is computed based on
  $20$ simulated datasets.}
\label{table:CorrDesgin}
\end{table}

We report the percentage of simulated datasets for which the GR scale
factor from six Markov chains is less than $1.5$ (success).  Moreover,
to see whether the variable-selection procedure based on the posterior
is consistent, we also compute the difference between the highest
posterior probability found during the Markov chain iterations and the
posterior probability of the true model (H-T) and the difference in
posterior probabilities between the null model and the true model
(N-T). If the true model receives the highest posterior probability, then H-T
would be $0$; if the null model receives the highest posterior probability,
then N-T would be the same as H-T.

Table~\ref{table:IndptDesgin} shows the results for design matrices
drawn from the independent ensemble. In this case, the Markov chain
method has fast convergence in all settings (it converges within
$20\covdim$ iterations). From the table, the setting $\SNR=0.5$
(respectively $\SNR\geq 2$) corresponds to the weak (respectively
strong) signal regime, while $\SNR=1$ is in the intermediate regime
where neither the null model nor the true model receives the highest
posterior probability. Table~\ref{table:CorrDesgin} shows the results
for design matrices drawn from the correlated ensemble. Now the Markov
chain method exhibits poor convergence behavior in the intermediate
regime $\SNR=1$ with $\numobs=500$, but still has fast convergence in
the weak and strong signal regimes. However, with larger sample size
$\numobs=1000$, the Markov chain has fast convergence in all settings
on $\covdim$ and $\SNR$. Comparing the results under the two different
designs, we find that correlations among the covariates increases the
difficulty of variable-selection tasks when Markov chain methods are
used.  Moreover, the results under the correlated design suggest that
there exists a regime, characterized by $\numobs$, $\covdim$ and
$\SNR$, in which the Markov chain is slowly mixing. It would be
interesting to see whether or not this regime characterizes some type
of fundamental limit on computationally efficient procedures for
variable selection.  We leave this question open as a possible future
direction.

%%%%%%%%%%%%%%%%%%%%%%%%%%%%%%%%%%%%%%%%%%%%%%%%%%%%%%%%%%%%%%%%%%%%%%%%%

\subsubsection{Bayesian methods versus the Lasso}

Our analysis reveals one possible benefit of a Bayesian approach as
opposed to $\ell_1$-based approaches such as the Lasso.  It is
well known that the performance of the Lasso and related
$\ell_1$-relaxations depends critically on fairly restrictive
incoherence conditions on the design matrix.  Here we provide an
example of an ensemble of linear regression problems for which the
Lasso fails to perform correct variable selection whereas the Bayesian
approach succeeds with high probability.

For Lasso-based methods, the irrepresentable condition
\begin{align}
\max_{|\gamma|=\sstar} \max_{ k \notin \gamma}
\|X_k^TX_{\gamma}(X_{\gamma}^T X_{\gamma})^{-1}\|_{1}<1
\end{align}
is both sufficient and necessary for variable-selection
consistency~\cite{Meinshausen2006,Zhao06,Wainwright09}.  In our theory
for the Bayesian approach, the analogous conditions are the upper
bound in Assumption D$(\maxsize)$ on the maximum model size, namely
\begin{align}
\maxsize\geq \big( 2\Cm^{-2} \, \omega(X) + 1 \big) \sstar,
\end{align}
as well as the \HACKER condition in Assumption B. Roughly speaking,
the first condition is needed to ensure that saturated models, i.e.,
models with size $\maxsize$, receive negligible posterior probability,
such that if too many unimportant covariates are included the removal
of some of them does not hurt the goodness of fit (see
Lemma~\ref{LemForwardSelection} in the Appendix).  This condition is
weaker than the irrepresentable condition since we can always choose
$\maxsize$ large enough so that $\maxsize\geq
\big(2\Cm^{-2}\,\omega(X)+1\big)\sstar$ holds, as long as Assumption B
is not violated.

As an example, consider a design matrix $\Xmat \in \real^{\numobs
  \times \covdim}$ that satisfies
\begin{align*}
\frac{1}{\numobs} \Xmat^T \Xmat = \Sigmabad \defn \begin{bmatrix} 1 &
  \mu &\mu &\cdots & \cdots & \mu\\ \mu & 1 & 0 & \cdots& \cdots
  &0\\ \mu & 0 & 1 & \cdots& \cdots &0\\ \vdots & \vdots & \vdots &
  \vdots & \vdots & \vdots \\ \mu & 0 & 0 & \cdots& \cdots & 1\\
\end{bmatrix} \in \real^{\covdim\times\covdim},
\end{align*}
with $\mu = (2\sqrt{\covdim})^{-1}$.  (When $\covdim>\numobs$, we may
consider instead a random design $\Xmat$ where the rows of $\Xmat$ are
generated i.i.d.\!  from the $\covdim$-variate normal distribution
$\mathcal{N} (0,\, \Sigmabad)$.)  This example was previously analyzed
by Wainwright~\cite{Martin2009}, who shows that it is an interesting
case in which there is a gap between the performance of $\ell_1$-based
variable-selection recovery and that of an optimal (but
computationally intractable) method based on searching over all
subsets.  For a design matrix of this form, we have
$\max_{|\gamma|=\sstar, \, k\notin\gamma}$ $\|\Xmat_k^T \Xmat_{\gamma}
(\Xmat_{\gamma}^T \Xmat_{\gamma})^{-1}\|_{1} \geq \sstar \mu$, so that
the irrepresentable condition fails if $\sstar > 2\sqrt{\covdim}$.
Consequently, by known results on the necessity of the irrepresentable
condition for Lasso~\cite{Zhao06,Wainwright09}, it will fail in
performing variable selection for this ensemble.

On the other hand, for this example, it can be verified that
Assumption D$(\maxsize)$ is satisfied with $\maxsize \geq 13 \sstar$,
and moreover, that the the RE$(s)$ condition in Assumption B holds
with $\Cm = 1/2$, whereas the \HACKER condition is satisfied with $L =
16 (1 + \maxsize^2\, \mu^2) = 16 + \frac{4\maxsize^2}{\covdim}$.  The
only consequence for taking larger values of $L$ is in the
$\betamin$-condition: in particular, the threshold $\CB$ is always
lower bounded by $L \frac{\log \covdim}{\numobs}$.  Consequently, our
theory shows that the Bayesian procedure will perform correct variable
selection with high probability for this ensemble.

\begin{figure}[h]
\begin{center}
\widgraph{.60\textwidth}{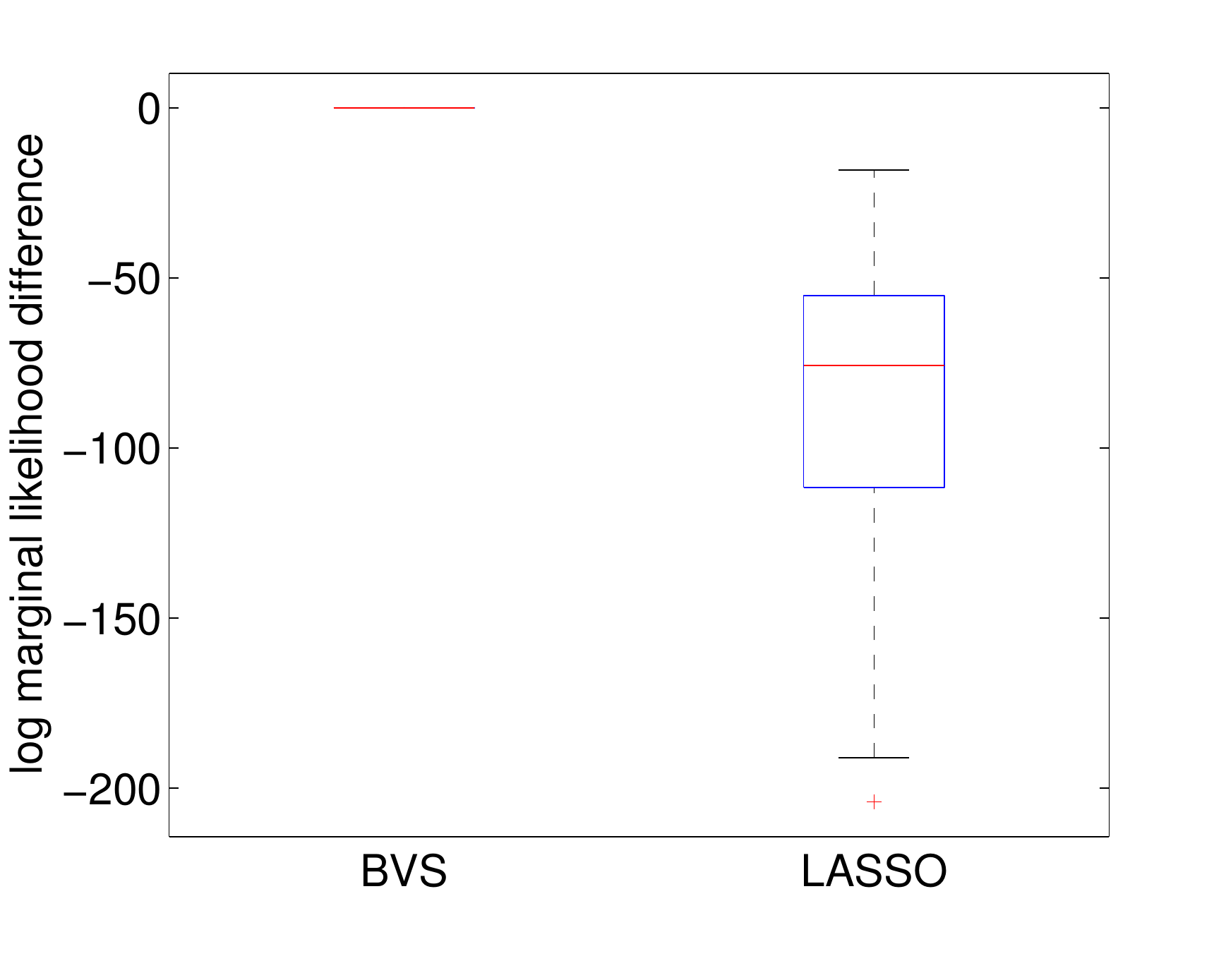}
\end{center}
 \caption{Boxplots indicating variable-selection performance of the
   Bayesian approach (BVS) and the Lasso. The boxplots are based on the
   logarithms of the ratio between the posterior probability of the
   selected model and the true model over $100$ replicates. The model
   selected by the Bayesian approach is the median probability
   model~\cite{Barbieri2004} and the regularization parameter of the
   Lasso is chosen by cross-validation.}
\label{FigComparison}
\end{figure}

To compare the performance of the Bayesian approach and the Lasso
under this setup, we generate our design matrix from a Gaussian
version of this ensemble; i.e., the rows of $\Xmat$ are generated 
i.i.d.\! from the $\covdim$-variate normal distribution $\mathcal{N}
(0,\, \Sigmabad)$. We choose $(\numobs,\covdim, \sstar) = (300, 80,
20)$ so that $\sstar\mu = 10/\sqrt{80} \approx 1.1 > 1$, i.e. the
irrepresentable condition fails.  Figure~\ref{FigComparison} shows the
variable-selection performance for the Bayesian approach and the Lasso
over $100$ replicates. We report the logarithm of the ratio between
the posterior probability (see equation~\eqref{EqnPostFormula}) of the
selected model and the true model, where we use the median probability
model \cite{Barbieri2004} as the selected model of the Bayesian
approach. If a variable-selection approach has good performance, then
we will expect this logarithm to be close to
zero. Figure~\ref{FigComparison} shows that the Bayesian approach
almost always selects the true model while the Lasso fails most of the
time, which is consistent with the theory.

%%%%%%%%%%%%%%%%%%%%%%%%%%%%%%%%%%%%%%%%%%%%%%%%%%%%%%%%%%%%%%%%%%%%%%%%%%%%

\section{Proofs}
\label{SecProofMain}

We now turn to the proofs of our main results, beginning with the
rapid mixing guarantee in Theorem~\ref{ThmMain}, which is the most
involved technically.  We then use some of the machinery developed in
Theorem~\ref{ThmMain} to prove the posterior consistency guarantee in
Theorem~\ref{ThmBVSconsistency}. Finally, by combining these two theorems we prove Corollary~\ref{CoroMCMC}. In order to promote readability, we
defer the proofs of certain more technical results to the appendices.

%%%%%%%%%%%%%%%%%%%%%%%%%%%%%%%%%%%%%%%%%%%%%%%%%%%%%%%%%%%%%%%%%%%%%%%%%

\subsection{Proof of Theorem~\ref{ThmMain}}

For the purposes of this proof, let $\widetilde{\PMAT}$ denote the
transition matrix of the original \mbox{Metropolis-Hastings}
sampler~\eqref{EqnMetropolisHastings}.  Now consider instead the
transition matrix $\PMAT \defn \widetilde{\PMAT}/2 + \mathbf{I}/2$,
corresponding to a lazy random walk that has a probability of at least
$1/2$ in staying in its current state.  By construction, the smallest
eigenvalue of $\PMAT$ will always be nonnegative, and as a
consequence, the mixing time of the Markov chain $\MarkovChain$ is
completely determined by the second largest eigenvalue $\lambda_2$ of
$\PMAT$.  The difference $\GAP(\PMAT) \defn 1 - \lambda_2$ is known as
the spectral gap, and for any lazy Markov chain, we have the sandwich
relation
\begin{align}
\label{EqnMixingTime}
\frac{1}{2} \frac{(1-\mbox{Gap}(\PMAT))}{ \mbox{Gap}(\PMAT)} \; \log
\big[1/(2\epsilon) \big] \; \leq \; \tau_{\epsilon} \leq \frac{
  \big(\log \big[1/\min \limits_{\gamma\in\Mspace}\pi(\gamma)\big]
  +\log (1/\epsilon) \big)}{\mbox{Gap}(\PMAT)} \: .
\end{align}
See the papers~\cite{Sinclair1992,Woodard2013} for bounds of this
form.

Using this sandwich relation, we claim that it suffices to show that
there are universal constants $(\UNICON_1,\UNICON_2)$ such that with
probability at least $1-4\covdim^{-c_1}$, the spectral gap of the lazy
transition matrix $\PMAT$ is lower bounded as
\begin{align}
\label{EqnSpectralBound}
\mbox{Gap}(\PMAT) \geq \frac{c_2}{\covdim\,\maxsize^2}.
\end{align}
To establish the sufficiency of this intermediate claim, we apply
Theorem~\ref{ThmBVSconsistency} and make use of the
expression~\eqref{EqnPostFormula} for the posterior distribution,
thereby obtaining that for $\gamma \in \Mspace$, the posterior
probability is lower bounded as
\begin{align*}
\posterior(\gamma \mid \yobs) & = \posterior(\gammastar \mid \yobs)
\cdot \frac{\posterior(\gamma \mid \yobs)} {\posterior(\gammastar \mid
  \yobs)}\\
& \geq e^{-2/\covdim} \cdot (\covdim \sqrt{1 +
  \hyperpara})^{-(|\gamma| - |\gammastar|)} \cdot \frac{\big(1 +
  \hyperpara(1 - R_{\gammastar}^2)\big)^{\numobs/2}}{\big(1 +
  \hyperpara (1 - R_{\gamma}^2) \big)^{\numobs/2}} \\
 & \geq e^{-2/\covdim} \cdot \covdim^{-(1 + \alpha/2)\maxsize} \cdot
\covdim^{-\alpha \numobs/2}
\end{align*}
with probability at least $1-4\covdim^{-c_1}$.  Combining the above
two displays with the sandwich relation~\eqref{EqnMixingTime}, we
obtain that there exist constants $(\CONTWO_1, \CONTWO_2)$ such that
for $\epsilon\in(0,1)$,
\begin{align*}
\tau_{\epsilon} \leq \CONTWO_1 \, \covdim \maxsize^2 \, \big( \CONTWO_2
\alpha \,(\numobs + \maxsize) \log \covdim + \log (1/\epsilon) + 2 \big)
\end{align*}
with probability at least $1-4\covdim^{-c_1}$.

\vspace*{.1in}

Accordingly, the remainder of our proof is devoted to establishing the
spectral gap bound~\eqref{EqnSpectralBound}, and we do so via a
version of the canonical path argument~\cite{Sinclair1992}.  Let us
begin by describing the idea of a canonical path ensemble associated
with a Markov chain. Given a Markov chain $\MarkovChain$ with state
space $\Mspace$, consider the weighted directed graph
$G(\MarkovChain)=(V,\, E)$ with vertex set $V = \Mspace$ and edge set
$E$ in which an ordered pair $e = (\gamma,\gamma')$ is included as an
edge with weight $\QMAT(e)=\QMAT(\gamma,\gamma')= \pi(\gamma) \PMAT(\gamma,\gamma')$ if and only if
$\PMAT(\gamma, \gamma') > 0$.  A \emph{canonical path ensemble}
$\PathEns$ for $\MarkovChain$ is a collection of paths that contains,
for each ordered pair $(\gamma,\gamma')$ of distinct vertices, a
unique simple path $\Path_{\gamma,\gamma'}$ in the graph that connects
$\gamma$ and $\gamma'$.  We refer to any path in the ensemble
$\PathEns$ as a canonical path.  

In terms of this notation, Sinclair~\cite{Sinclair1992} shows that for
any reversible Markov chain and any choice of canonical path $\Path$,
the spectral gap of $\PMAT$ is lower bounded as
\begin{align}
\label{EqnCanonicalPathBound}
\underbrace{\mbox{Gap}(\PMAT)}_{1-\lambda_2} & \geq
\frac{1}{\rho(\Path) \ell(\PathEns)},
\end{align}
where $\ell(\PathEns)$ corresponds to the length of a longest path in
the ensemble $\PathEns$, and the quantity \mbox{$\rho(\Path) \defn
  \max \limits_{e \in E} \frac{1}{\QMAT(e)} \sum
  \limits_{\Path_{\gamma,\gamma'}\ni e} \pi(\gamma)\pi(\gamma')$} is
known as the \emph{path congestion parameter}.

\vspace*{.1in}

%%%%%%%%%%%%%%%%%%%%%%%%%%%%%%%%%%%%%%%%%%%%%%%%%%%%%%%%%%%%%%%%%%%%%%%%%%%%

In order to apply this approach to our problem, we need to construct a
suitable canonical path ensemble $\PathEns$.  To begin with, let us
introduce some notation for operations on simple paths. For two 
given paths $\Path_1$ and $\Path_2$:
\begin{itemize}
\item Their intersection $\Path_1 \cap  \Path_2$ corresponds to the
  subset of overlapping edges.  (For instance, if $\Path_1=(1,1,1) \to
  (0,1,1) \to (0,0,1) \to (0,0,0)$ and $\Path_2=(0,0,1) \to (0,0,0)$,
  then $\Path_1 \cap \Path_2 = (0,0,1) \to (0,0,0)$.)
\item If $\Path_2\subset \Path_1$, then $\Path_1 \setminus \Path_2$
  denotes the path obtained by removing all edges in $\Path_2$ from
  $\Path_1$.  (With the same specific choices of $(\Path_1, \Path_2)$
  as above, we have $\Path_1\setminus \Path_2=(1,1,1)\to (0,1,1)\to
  (0,0,1)$.)
\item We use $\bar{\Path_1}$ to denote the reverse of $\Path_1$.
  (With the choice of $\Path_1$ as above, we have $\bar{\Path_1} =
  (0,0,0) \to (0,0,1) \to (0,1,1) \to (1,1,1)$.)
\item If the endpoint of $\Path_1$ and the starting point of $\Path_2$
  are the same, then we define the union $\Path_1\cup \Path_2$ as the
  path that connects $\Path_1$ and $\Path_2$ together.  (If
  $\Path_1=(0,0,0)\to (0,0,1)$ and $\Path_2 = (0,0,1) \to (0,1,1)$,
  then their union is given by $\Path_1\cup \Path_2 = (0,0,0) \to
  (0,0,1)\to (0,1,1)$.)

\end{itemize}

We now turn to the construction of our canonical path ensemble.  At a
high level, our construction is inspired by the variable-selection
paths carved out by greedy stepwise variable-selection procedures
(e.g., ~\cite{Zhang11,An08}).

%%%%%%%%%%%%%%%%%%%%%%%%%%%%%%%%%%%%%%%%%%%%%%%%%%%%%%%%%%%%%%%%%%%%
\paragraph{Canonical path ensemble construction for $\Mspace$:}

First, we construct the canonical path
$\Path_{\gamma,\gammastar}$ from any $\gamma\in\Mspace$ to the true
model $\gammastar$. The following construction will prove helpful. 
We call a set $\mathcal{R}$ of
canonical paths \emph{memoryless} with respect to the central state
$\gammastar$ if: (1) for any state $\gamma \in \Mspace$ satisfying
$\gamma\neq\gammastar$, there exists a unique simple path
$\Path_{\gamma,\gammastar}$ in $\mathcal{R}$ that connects $\gamma$
and $\gammastar$; (2) for any intermediate state $\gammatil\in\Mspace$
on any path $\Path_{\gamma,\gammastar}$ in $\mathcal{R}$, the unique
path $\Path_{\gammatil,\gammastar}$ in $\mathcal{R}$ that connects
$\gammatil$ and $\gammastar$ is the sub-path of
$\Path_{\gamma,\gammastar}$ starting from $\gammatil$ and ending at
$\gammastar$.  Intuitively, this memoryless property means that for
any intermediate state on any canonical path towards the central
state, the next move from this intermediate state towards the central
state does not depend on the history. A memoryless canonical path
ensemble has the property that in order to specify the canonical
path connecting any state $\gamma\in\Mspace$ and the central state
$\gammastar$, we only need to specify which state to move to from any
$\gamma\neq\gammastar$ in $\Mspace$; i.e., we need a transition function
$\Gfun:\Mspace\setminus\{\gammastar\}\to\Mspace$ that maps the current
state $\gamma\in\Mspace$ to a next state
$\Gfun(\gamma)\in\Mspace$. For simplicity, we define
$\Gfun(\gammastar)=\gammastar$ to make $\Mspace$ as the domain of
$\Gfun$. Clearly, each memoryless canonical path ensemble with respect
to a central state $\gammastar$ corresponds to a transition function
$\Gfun$ with $\Gfun(\gammastar)=\gammastar$, but the converse is not
true. For example, if there exist two states $\gamma$ and $\gamma'$ so
that $\Gfun(\gamma)=\gamma'$ and $\Gfun(\gamma')=\gamma$, then $\Gfun$
is not the transition function corresponding to any memoryless
canonical path ensemble. However, every valid transition function
$\Gfun$ gives rise to a unique memoryless canonical path set
consisting of paths connecting any $\gamma\in\Mspace$ to $\gammastar$,
with $\gammastar$ corresponding to the fixed point of $\Gfun$. 
We call function $\Gfun$ a valid transition function if there exists a
memoryless canonical path set for which $\Gfun$ is the corresponding
transition function. The next lemma provides a suffcient condition for
a function $\Gfun:\Mspace\setminus\{\gammastar\}\to\Mspace$ to be
valid, which motivates our construction to follow.  Recall that
$\hamm$ denotes the Hamming metric between a pair of binary strings.

\blems
\label{LemValidTranFun}
%%%
If a function $\Gfun:\Mspace\setminus\{\gammastar\}\to\Mspace$ satisfies that for any state $\gamma\in\Mspace\setminus\gammastar$, the Hamming distance between $\Gfun(\gamma)$ and $\gammastar$ is strictly less than the Hamming distance between $\gamma$ and $\gammastar$, then $\Gfun$ is a valid transition function.
\elems
\begin{proof}
Based on this function $\Gfun$, we can construct the canonical path
$\Path_{\gamma,\gammastar}$ from any state $\gamma\in\Mspace$ to
$\gammastar$ by defining $\Path_{\gamma,\gammastar}$ as $\gamma\to
\Gfun(\gamma)\to \Gfun^2(\gamma)\to\ldots\to
\Gfun^{k_{\gamma}}(\gamma)$, where $\Gfun^k\defn \Gfun\circ\ldots\circ
\Gfun$ denotes the $k$-fold self-composition of $\Gfun$ for any
$k\in\mathbb{N}$ and $k_{\gamma} \defn \min_{k} \{\Gfun^k(\gamma) =
\gammastar\}$.  In order to show that the set
$\{\Path_{\gamma,\gammastar}:\,\gamma\in\Mspace,\gamma\neq\gammastar\}$ is
a memoryless canonical path set, we only need to verify two things:
\begin{enumerate}
\item[(a)] for any $\gamma\neq\gammastar$, $\Path_{\gamma,\gammastar}$ is
  a well-defined path; i.e., it has finite length and ends at
  $\gammastar$, and
\item[(b)] for any $\gamma\neq\gammastar$, $\Path_{\gamma,\gammastar}$ is
  a simple path. 
\end{enumerate}
By our assumption, the function $F:\Mspace\to \real$ defined by
  $F(\gamma)=\hamm(\gamma,\gammastar)$ is strictly decreasing along the
  path $\Path_{\gamma,\gammastar}$ for $\gamma \neq \gammastar$.
Because $F$ only attains a finite number of values, there exists a
smallest $k_\gamma$ such that $\Gfun^{k+1}(\gamma) =\Gfun^k(\gamma)$
for each $k\geq k_\gamma$, implying that $\Gfun^{k_\gamma}(\gamma)$ is
a fixed point of $\Gfun$. Since $\gammastar$ is the
unique fixed point of $\Gfun$, we must have
$\Gfun^{k_\gamma}(\gamma)=\gammastar$, which proves the first claim. The
second claim is obvious since the function $F$ defined above is
strictly decreasing along the path $\Path_{\gamma,\gammastar}$, which
means that the states on the path $\Path_{\gamma,\gammastar}$ are all
distinct.
\end{proof}

\begin{figure}[t]
\centering \includegraphics[width=6in]{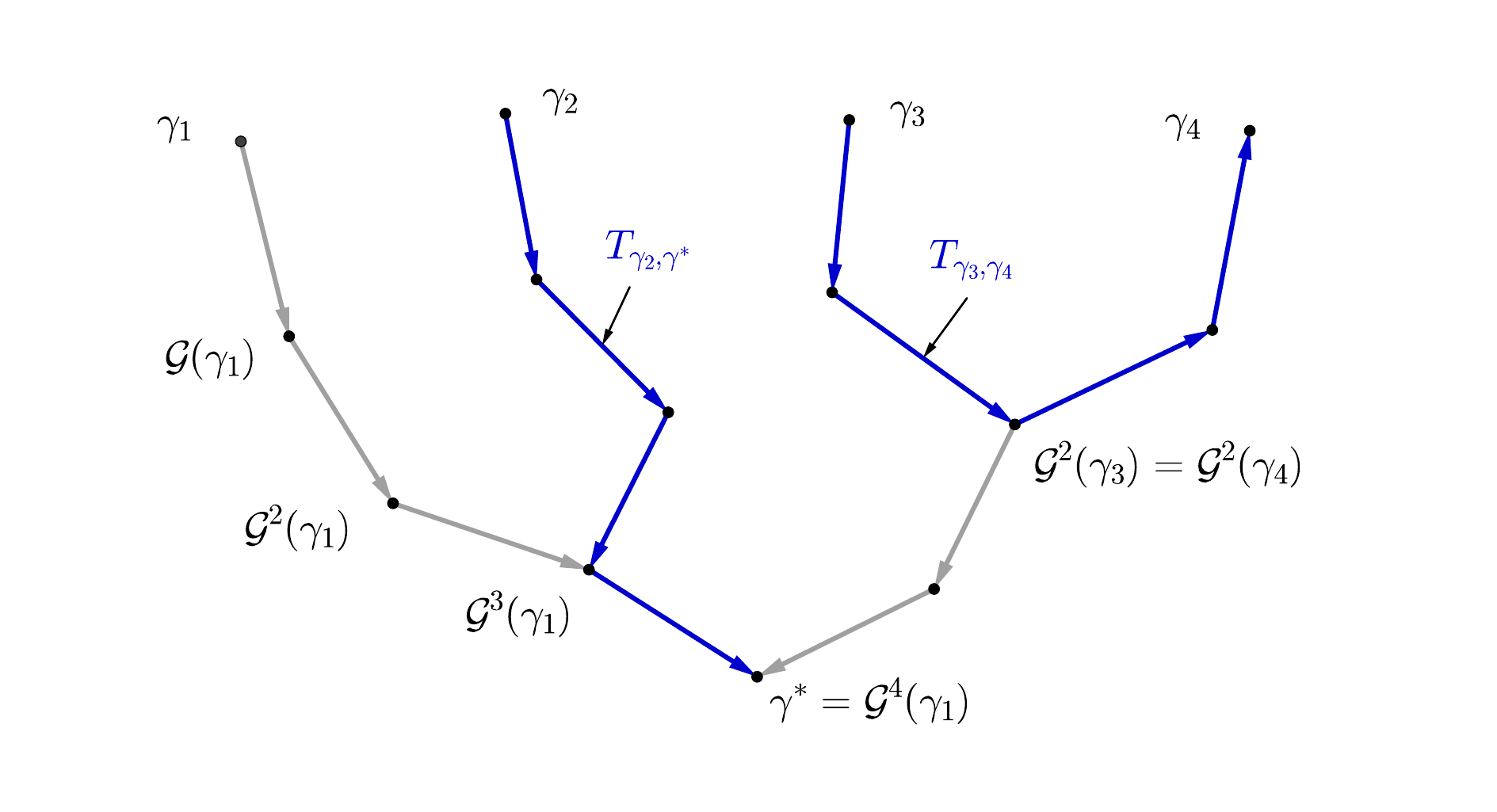}
\caption{Illustration of the construction of the canonical path ensemble. 
  In the plot, $\gamma^\ast$ is the central state, $\Gfun$ is the
  transition function and solid blue arrows indicate canonical
  paths $T_{\gamma_2,\gamma^\ast}$ and $T_{\gamma_3,\gamma_4}$.}
\label{FigIllustration}
\end{figure}

Equiped with this lemma, we start constructing a memoryless set of
canonical paths from any state $\gamma\in\Mspace$ to $\gammastar$ by
specifying a valid $\Gfun$ function. First, we introduce some
definitions on the states. A state $\gamma \neq \gammastar$ is called
\emph{saturated} if $|\gamma|=\maxsize$ and \emph{unsaturated} if
$|\gamma|<\maxsize$. We call a state $\gamma\neq\gammastar$
\emph{overfitted} if it contains all influential covariates,
i.e. $\gammastar\subset\gamma$, and \emph{underfitted} if it does not
contain at least one influential covariate. Recall the two updating
schemes in our Metropolis-Hastings (MH) sampler: single flip and
double flips. We accordingly construct the transition function $\Gfun$
as follows.
\begin{enumerate}
\item[(i)] If $\gamma \neq \gammastar$ is overfitted, then we define
  $\Gfun(\gamma)$ to be $\gamma'$, which is formed by deleting the
  least influential covariate from $\gamma$, i.e.  $\gamma'_j =
  \gamma_j$ for any $j\neq \ell_{\gamma}$ and $\gamma'_{\ell_{\gamma}}
  = 0$, where $\ell_{\gamma}$ is the index from the set $\gamma
  \setminus \gammastar$ of uninfluential covariates that minimizes
the difference
\begin{align*}
\|\Proj_{\gamma} X_{\gammastar} \betastar_{\gammastar}\|_2^2 -
\|\Proj_{\gamma\setminus \{\ell\}} X_{\gammastar}
\betastar_{\gammastar}\|_2^2,
\end{align*}
where $\Proj_\gamma$ denotes the projection onto the span of $\{X_j, j
\in \gamma \}$.  This transition remsembles the backward deletion step
in the stepwise variable-selection procedure and involves the single
flip updating scheme of the MH algorithm. By construction, if $\gamma
\neq \gammastar$ is overfitted, then
$\hamm(\Gfun(\gamma),\gammastar)=\hamm(\Gfun(\gamma),\gammastar) - 1$.
\item[(ii)] If $\gamma \neq \gammastar$ is underfitted and
  unsaturated, then we define $\Gfun(\gamma)$ to be $\gamma'$, which
  is formed by adding the influential covariate from $\gammastar
  \setminus \gamma$ that explains the most signal variation, i.e.
  $\gamma'_j=\gamma_j$ for any $j\neq j_{\gamma}$ and
  $\gamma'_{j_{\gamma}}=1$, where $j_{\gamma}$ is defined as the $j
  \in \gammastar \setminus \gamma$ that maximizes the quantity
  $\|\Proj_{\gamma\cup\{j\}}
  \design_{\gammastar}\beta_{\gammastar}^\ast\|_2^2$. This transition
  remsembles the forward selection step in the stepwise variable
  selection procedure and involves the single flip updating scheme of
  the MH algorithm. By construction, if $\gamma \neq \gammastar$ is
  underfitted and unsaturated, then
  $\hamm(\Gfun(\gamma),\gammastar)=\hamm(\Gfun(\gamma),\gammastar) -
  1$.
\item[(iii)] If $\gamma\neq\gammastar$ is underfitted and saturated,
  then we define $\Gfun(\gamma)$ to be $\gamma'$, which is formed by
  replacing the least influential unimportant covariate in $\gamma$
  with the most influential covariate from $\gammastar \setminus
  \gamma$, i.e. $\gamma'_j = \gamma_j$ for any $j\not\in
  \{j_{\gamma},k_{\gamma}\}$, $\gamma'_{j_{\gamma}} =1$ and
  $\gamma'_{k_{\gamma}}=0$, where $j_{\gamma}$ is defined in case 2
  and $k_{\gamma}\in\gamma\setminus\gammastar$ minimizes
  $\|\Proj_{\gamma\cup\{j\}} X_{\gammastar} \betastar_{\gammastar}\|_2^2 -
  \|\Proj_{\gamma \cup\{j\} \setminus\{k\}} X_{\gammastar}
  \betastar_{\gammastar}\|_2^2$. This transition step involves the 
  double-flip updating scheme of the MH algorithm. By construction, if
  $\gamma \neq \gammastar$ is underfitted and saturated, then
  $\hamm(\Gfun(\gamma),\gammastar)=\hamm(\Gfun(\gamma),\gammastar) -
  2$.
\end{enumerate}
By Lemma~\ref{LemValidTranFun}, this transition function $\Gfun$ is
valid and gives rise to a unique memoryless set of canonical paths
from any state $\gamma\in\Mspace$ to $\gammastar$. For example,
Fig~\ref{FigIllustration} shows such a memoryless set of canonical
paths for $\Mspace$ consisting of $14$ states, where
$T_{\gamma_2,\gamma^\ast}$ corresponds to the canonical path from
state $\gamma_2$ to the central state $\gamma^\ast$.

Based on this memoryless canonical path set, we can finish
constructing the canonical path ensemble $\PathEns$ by specifying the
path $\Path_{\gamma,\gamma'}$ connecting any distinct pair
$(\gamma,\gamma') \in \Mspace \times \Mspace$. More specifically, by
the memoryless property, the two simple paths
$\Path_{\gamma,\gammastar}$ and $\Path_{\gamma',\gammastar}$ share an
identical subpath towards $\gammastar$ from their first common
intermediate state.  Let $\Path_{\gamma \cap \gamma'}$ denote this
common subpath $\Path_{\gamma,\gammastar} \cap
\Path_{\gamma',\gammastar}$, and $\Path_{\gamma\setminus\gamma'} \defn
\Path_{\gamma,\gammastar}\setminus \Path_{\gamma \cap \gamma'}$ denote
the remaining path of $\Path_{\gamma,\gammastar}$ after removing the
segment $\Path_{\gamma \cap \gamma'}$. We define $\Path_{\gamma'
  \setminus \gamma}$ in a similar way as
$\Path_{\gamma',\gammastar}\setminus \Path_{\gamma \cap
  \gamma'}$. Then it is easy to see that the two remaining paths
$\Path_{\gamma \setminus \gamma'}$ and $\Path_{\gamma' \setminus
  \gamma}$ share the same endpoint. Therefore, it is valid to define
the path $\Path_{\gamma,\gamma'}$ as $\Path_{\gamma \setminus \gamma'}
\cup \bar{\Path}_{\gamma'\setminus \gamma}$. To understand this
construction, let us consider an example where $\Path_{\gamma,
  \gammastar}=(0,1,1,1) \to (1,1,0,1) \to (1,1,0,0)$ and
$\Path_{\gamma', \gammastar} = (1,0,0,1) \to (1,1,0,1) \to
(1,1,0,0)$. Their intersection is $\Path_{\gamma \cap
  \gamma'}=(1,1,0,1) \to (1,1,0,0)$ and the two remaining paths are
$\Path_{\gamma \setminus \gamma'} = (0,1,1,1) \to (1,1,0,1)$ and
$\Path_{\gamma'\setminus\gamma}=(1,0,0,1)\to (1,1,0,1)$. Consequently,
the path $\Path_{\gamma,\gamma'}$ from $\gamma$ to $\gamma'$ is
$(0,1,1,1) \to (1,1,0,1) \to (1,0,0,1)$ by our construction. For
example, path $T_{\gamma_3,\gamma_4}$ in Fig~\ref{FigIllustration}
illustrates the construction of the path connecting
$(\gamma_3,\gamma_4)$ when $\Mspace$ is composed of $14$ states.

We call $\gamma$ a \emph{precedent} of $\gamma'$ if $\gamma'$ is on the
canonical path $\Path_{\gamma,\gammastar} \in \PathEns$, and a pair of
states $\gamma, \gamma'$ \emph{adjacent} if the canonical path
$\Path_{\gamma, \gamma'}$ is $e_{\gamma, \gamma'}$, the edge in $E$
connecting $\gamma$ and $\gamma'$. For $\gamma\in\Mspace$, let
\begin{align}
\label{EqnPrecedent}
\PREC(\gamma) \defn \{\gammabar \mid \gamma\in
\Path_{\gammabar,\gammastar}\}
\end{align}
denote the set of all its precedents. Use the notation $|\Path|$ to denote the length of a path $\Path$.
The following lemma provides some important properties of the contructed canonical path ensemble that will be used later.

%%%%%
\blems
\label{LemPathProperties}
For any distinct pair $(\gamma, \gamma') \in \Mspace \times \Mspace$:
\begin{enumerate}
\item[(a)]  We have
\begin{subequations}
\begin{align}
|\Path_{\gamma, \gammastar}| & \leq \hamm(\gamma, \gammastar) \; \leq
\; \maxsize, \quad \mbox{and} \\
|\Path_{\gamma,\gamma'}| & \leq \hamm(\gamma,\gammastar) +
\hamm(\gamma',\gammastar) \leq 2 \maxsize.
\end{align}
\end{subequations}
\item[(b)] If $\gamma$ and $\gamma'$ are adjacent (joined by edge
  $e_{\gamma, \gamma'}$) and $\gamma$ is a precedent of $\gamma'$, then
\begin{align*}
\{(\gammabar, \gammabar') \mid \, \Path_{\gammabar, \gammabar'} \ni
e_{\gamma, \gamma'}\} \subset \PREC(\gamma) \times \Mspace,
\end{align*}
\end{enumerate}
\elems
%%%%
%%%%%%%%%%%%%%%%%%%%%%%%%%%%%%%%%%%%%%%%%%%%%%%%%%%%%%%%%%%%%%%%%%%%%%

\begin{proof}
The first claim follows since the function $F:\Mspace\to \real$ defined by
  $F(\gamma)=\hamm(\gamma,\gammastar)$ is strictly decreasing along the
  path $\Path_{\gamma,\gammastar}$ for $\gamma \neq \gammastar$. Now we prove
  the second claim. For any pair $(\gammabar,\gammabar')$ such that
  $\Path_{\gammabar,\gammabar'} \ni e_{\gamma,\gamma'}$, either
  $e_{\gamma,\gamma'}\in \Path_{\gammabar\setminus \gammabar'}$ or
  $e_{\gamma',\gamma}\in \Path_{\gammabar'\setminus \gammabar}$ should
  be satisfied since $\Path_{\gammabar,\gammabar'} = \Path_{\gammabar
    \setminus \gammabar'} \cup \bar{\Path}_{\gammabar',\gammabar}$ by
  our construction. Because $\gamma$ is a precedent of $\gamma'$, we
  can only have $e_{\gamma,\gamma'}\in \Path_{\gammabar\setminus
    \gammabar'}$. This shows that $\gamma$ is on the path
  $\Path_{\gammabar,\gammastar}$ and $\gammabar \in \PREC(\gamma)$.
\end{proof}
%%%

According to Lemma~\ref{LemPathProperties}\,(b), the path congestion
parameter $\rho(\Path)$ of the canonical path $\Path$ satisfies
\begin{align}
\label{EqnPCP}
\rho(\Path) \leq \max_{(\gamma,\gamma')\in \Gammastar}
\frac{1}{\QMAT(\gamma,\gamma')}\, \sum_{\gammabar\in \PREC(\gamma),\,
  \gammabar'\in\Mspace} \Stat (\gammabar) \Stat(\gammabar') =
\max_{(\gamma, \gamma')\in \Gammastar}
\frac{\Stat[\PREC(\gamma)]}{\QMAT(\gamma, \gamma')},
\end{align}
where the maximum is taken over the set
\begin{align*}
\Gammastar & \defn \Big \{ (\gamma,\gamma') \in \Mspace \times \Mspace
\, \mid \Path_{\gamma,\gamma'} = e_{\gamma,\gamma'} \mbox{ and $\gamma
  \in \PREC(\gamma')$} \Big \}.
\end{align*}
Here we used the fact that the weight function $\QMAT$ of a reversible
chain satisfies \mbox{$\QMAT(\gamma, \gamma') =
  \QMAT(\gamma',\gamma)$} so as to be able to restrict the range of
the maximum to pairs $(\gamma,\gamma')$ where $\gamma \in
\PREC(\gamma')$.

For the lazy form of the Metropolis-Hastings
walk~\eqref{EqnMetropolisHastings}, given any pair $(\gamma,\gamma')$
such that $\PMAT(\gamma, \gamma') > 0$, we have
\begin{align*}
\QMAT(\gamma, \gamma') & = \frac{1}{2} \posterior (\gamma \mid \yobs)
\PMAT(\gamma,\gamma') \\
& \geq \frac{1}{2\, \covdim \, \maxsize} \, \posterior(\gamma \mid
\yobs) \min \big\{1, \frac{ \posterior(\gamma'\mid
  \yobs)}{\posterior(\gamma \mid \yobs)} \big\} = \frac{1}{2\, \covdim
  \, \maxsize} \min \big \{\posterior(\gamma' \mid \yobs), \,
\posterior ( \gamma \mid \yobs )\big\}.
\end{align*}
Substituting this lower bound into our upper bound~\eqref{EqnPCP} on
the path congestion parameter yields
\begin{align}
\rho(\Path) & \leq 2 \, \covdim \, \maxsize
\,\max_{(\gamma,\gamma')\in \Gammastar} \frac{\posterior[\PREC(\gamma)
    \mid \yobs]} {\min\big\{\posterior(\gamma \mid
  \yobs),\posterior(\gamma' \mid \yobs)\big\}} \notag \\
\label{EqnEspresso}
& = 2 \, \covdim \, \maxsize \,\max_{(\gamma,\gamma') \in \Gammastar}
\Big \{\max \Big\{1, \frac{\posterior(\gamma \mid
  \yobs)}{\posterior(\gamma' \mid \yobs)} \Big\} \cdot
\frac{\posterior[\PREC(\gamma) \mid \yobs]}{\posterior(\gamma \mid
  \yobs)} \Big\}.
\end{align}
In order to prove that $\rho(\Path) = \order(\covdim \maxsize)$ with
high probability, it suffices to show that the two terms inside the
maximum are $\order(1)$ with high probability.  In order to do so, we
make use of two auxiliary lemmas.

Given the constant $\rcon \geq 4$ and the noise vector $\wnoise \sim
\neigh(0,\sigmazero^2I_\numobs)$, consider the following events
\begin{subequations}
\label{EqnKeyEvents}
\begin{align}
\Aevent_\numobs & \defn \Big\{\max_{\substack{{(\gamma_1, \gamma_2)
      \in \Mspace \times \Mspace} \\ {\gamma_2 \subset
      \gamma_1}}}\frac{\wnoise^T (\Proj_{\gamma_1} - \Proj_{\gamma_2})
  \wnoise}{|\gamma_1| - |\gamma_2|} \leq L \sigmazero^2 \log \covdim
\Big\}, \\
\Bevent_\numobs & \defn \Big\{
\max_{\gamma\in\Mspace}\frac{\wnoise^T\Proj_{\gamma}
  \wnoise}{|\gamma|} \leq r \sigmazero^2\log \covdim\Big\}, \quad
\mbox{and} \\
\Cevent_\numobs & \defn \Big\{ \Big|
\frac{\|\wnoise\|_2^2}{\numobs\sigmazero^2} - 1 \Big| \leq \frac{1}{2}
\Big\}, \quad \mbox{and} \quad \Devent_\numobs\defn \Big\{
\frac{\|\yobs\|_2^2}{\hyperpara} \leq \frac{5 \numobs
  \sigmazero^2}{\sstar} \Big\}.
\end{align}
\end{subequations}
%%%%%%%%%%%%%%%%%%%%%%%%%%%%%%%%%%%%%%%%%%%%%%%%%%%%%%%%%%%%%%%%%%%%%%

Our first auxiliary lemma guarantees that, under the stated
assumptions of our theorem, the intersection of these events holds
with high probability:
\blems
\label{LemConcentrationIneq}
%%%
Under the conditions of Theorem~\ref{ThmMain}, we have
\begin{align}
\Prob(\Aevent_\numobs \cap \Bevent_\numobs \cap \Cevent_\numobs \cap
\Devent_\numobs) & \geq 1 - 6 \covdim^{-\UNICON}.
\end{align}
%%%
\elems
\noindent We prove this lemma in Section~\ref{SecConcentrationInq} to
follow. \\

Our second auxiliary lemma ensures that when these four events hold,
then the two terms on the right-hand side of the upper
bound~\eqref{EqnEspresso} are controlled.
\blems
\label{LemPosteriorConcentration}
Suppose that, in addition to the conditions of Theorem~\ref{ThmMain},
the compound event $\Aevent_\numobs \cap \Bevent_\numobs \cap
\Cevent_\numobs \cap \Devent_\numobs$ holds.  Then for all $\gamma
\neq \gammastar$, we have
\begin{subequations}
\begin{align}
\label{EqnPosteriorA}
\frac{\posterior(\gamma \mid \yobs)}{\posterior(\Gfun(\gamma) \mid
  \yobs)} & \leq
\begin{cases} 
\covdim^{-2}, & \mbox{if $\gamma$ is overfitted},\\
\covdim^{-3}, & \mbox{if $\gamma$ is underfitted},
\end{cases}
\end{align}
and moreover, for all $\gamma$,
\begin{align}
\label{EqnPosteriorB}
\frac{\posterior[\PREC(\gamma) \mid \yobs]}{\posterior(\gamma \mid
  \yobs)} \leq c \qquad \mbox{for some universal constant $c$.}
\end{align}
\end{subequations}

\elems
\noindent We prove this lemma in
Section~\ref{SecPosteriorConcentration} to follow. \\

\noindent Combining Lemmas~\ref{LemConcentrationIneq}
and~\ref{LemPosteriorConcentration} with our earlier
bound~\eqref{EqnEspresso}, we conclude that \mbox{$\rho(\Path) \leq 2
  \, c \, \covdim \, \maxsize$.}  By
Lemma~\ref{LemPathProperties}\,(a), our path ensemble $\PathEns$ has
maximal length $\ell(\PathEns) \leq 2 \maxsize$, and hence the
canonical path lower bound~\eqref{EqnCanonicalPathBound} implies that
$\GAP(\PMAT) \geq \frac{1}{4 c \covdim \, \maxsize^2}$, as claimed in
inequality~\eqref{EqnSpectralBound}.  This completes the proof of the
theorem. \\

\vspace*{.1in}

The only remaining detail is to prove
Lemmas~\ref{LemConcentrationIneq} and~\ref{LemPosteriorConcentration},
and we do so in the following two subsections.

%%%%%%%%%%%%%%%%%%%%%%%%%%%%%%%%%%%%%%%%%%%%%%%%%%%%%%%%%%%%%%%%%%%%%%%%%%

\subsection{Proof of Lemma~\ref{LemConcentrationIneq}}
\label{SecConcentrationInq}

We split the proof up into separate parts, one for each of the events
$\Aevent_\numobs, \Bevent_\numobs, \Cevent_\numobs$ and
$\Devent_\numobs$.

\paragraph{Bound on $\mprob[\Cevent_\numobs]$:}
Since $\|\wnoise\|_2^2/\sigmazero^2\sim \chi^2_\numobs$, a standard
tail bound for the $\chi^2_\numobs$ distribution
(e.g.,~\cite{Laurent2000}, Lemma 1) yields
\begin{align}
\label{EqnConcenBound3}
\Prob \big[ \Cevent_\numobs \big] & \geq 1- 2 e^{
  -\frac{\numobs}{25}}.
\end{align}

\paragraph{Bound on $\mprob[\Bevent_\numobs]$:}
For each state $\gamma \in \Mspace$, the random variable $w^T
\Proj_{\gamma} w / \sigmazero^2$ follows a chi-squared distribution
with $|\gamma|$ degrees of freedom.  For each integer $\ell \in \{1,
\ldots, \maxsize \}$, the model space $\Mspace$ contains
${\covdim\choose \ell}$ models of size $\ell$.  Therefore, by a union
bound, we find that
\begin{align}
\Prob[\Bevent_\numobs] \geq 1- \sum_{\ell=1}^{\maxsize} {\covdim
  \choose \ell} \, \Prob(\chi^2_\ell \geq \rcon \ell \log \covdim) &
\geq 1- \sum_{l=1}^{\maxsize} e^{-(\rcon/4-1)\, \ell \log
  \covdim}\notag \\
& \geq 1- 2 e^{-(\rcon/4-1)\, \log \covdim} \notag \\
& = 1- 2 \covdim^{-(\rcon/4-1)}.
\label{EqnConcenBound1}
\end{align}

\paragraph{Bound on $\mprob[\Devent_\numobs]$:}
Given the linear observation model, we have
\begin{align*}
\|\yobs\|_2^2 = \| \design \betastar + \wnoise\|_2^2 \leq 2\|
\design\betastar \|^2 + 2 \|\wnoise\|_2^2.
\end{align*}
Combining this with inequality~\eqref{EqnConcenBound3}, we obtain
\begin{align*}
\Prob \big[ \|\yobs\|_2^2 \geq 2\|\design\betastar \|_2^2 + 3 \numobs
  \sigmazero^2 \big] & \leq 2 e^{- \frac{\numobs}{25}} \leq
\covdim^{-\maxsize(r/4-1)}
\end{align*}
for large $\numobs$ and some constant $C>0$, where we have used
Assumption D.  By Assumptions A and D, we have $\|\design\betastar
\|_2^2 \leq 2\numobs \sigmazero^2 \hyperpara/\sstar $, implying that
\begin{align}
\label{EqnConcenBound4}
\Prob \big[ \Devent^c_\numobs \big] \leq \Prob \big[ \|\yobs\|_2^2
  \geq 2\|\design\betastar \|_2^2 + 3 \numobs \sigmazero^2 \big] \leq
\covdim^{-\maxsize(r/4-1)}.
\end{align}

\paragraph{Bound on $\mprob[\Aevent_\numobs]$:}  To control this probability,
we require two auxiliary lemmas.

\blems
\label{LemUnderfit}
Under Assumption B, for any distinct pair $(\gamma,
\gammabar)\in \Mspace \times \Mspace$ satisfying $\gamma\subset
\gammabar$, we have
\begin{align*}
\lammin\Big(\frac{1}{\numobs} \design_{\gammabar\setminus \gamma}^T
(I_\numobs-\Proj_{\gamma}) \design_{\gammabar\setminus\gamma}\Big)
\geq \Cm.
\end{align*}
\elems
\begin{proof}
By partitioning the matrix $\design_{\gammabar}$ into a block form
$(\design_{\gamma},\design_{\gammabar\setminus\gamma})$ and using the
formula for the inverse of block matrices, one can show that the lower
right corner of $\big( \numobs^{-1} \design_{\gammabar}^T
\design_{\gammabar} \big)^{-1}$ is $\big(\numobs^{-1}
\design_{\gammabar\setminus \gamma}^T (I_\numobs - \Proj_{\gamma})
\design_{\gammabar \setminus \gamma} \big)^{-1}$, which implies the
claimed bound.
\end{proof}

\blems
\label{LemProjection}
For $\gamma \in \Mspace$ and $k \notin \gamma$, we have
\begin{align*}
\Proj_{\gamma\cup\{k\}} - \Proj_{\gamma} = \frac{(I - \Proj_{\gamma})
  \Xmat_k \Xmat_k^T (I - \Proj_{\gamma})} {\Xmat_k^T (I -
  \Proj_{\gamma}) \Xmat_k}.
\end{align*}
\elems
\begin{proof}

By the block matrix inversion formula~\cite{Horn85}, we have
\begin{align*}
\begin{bmatrix}
\Xmat_{\gamma}^T \Xmat_{\gamma} & \Xmat_{\gamma}^T \Xmat_k\\ \Xmat_k^T
\Xmat_{\gamma} & \Xmat_k^T \Xmat_k\\
\end{bmatrix}^{-1} = 
\begin{bmatrix}
B + a B \Xmat_{\gamma}^T \Xmat_{k} \Xmat_k^T \Xmat B & -a B
\Xmat_{\gamma}^T \Xmat_k\\ -a \Xmat_k^T \Xmat_{\gamma} B & a \\
\end{bmatrix},
\end{align*}
where $B = ( \Xmat_{\gamma}^T \Xmat_{\gamma})^{-1} \in \real^{|\gamma|
  \times| \gamma|}$ and $a=( \Xmat_k^T (I - \Proj_{\gamma}) \Xmat_k
)^{-1} \in \real$.  Then simple linear algebra yields
\begin{align*}
\Proj_{\gamma\cup\{k\}} - \Proj_{\gamma} & =
\begin{bmatrix}
\Xmat_{\gamma} & \Xmat_k\\
\end{bmatrix}
\begin{bmatrix}
\Xmat_{\gamma}^T \Xmat_{\gamma} & \Xmat_{\gamma}^T \Xmat_k
\\ \Xmat_k^T \Xmat_{\gamma} & \Xmat_k^T \Xmat_k \\
\end{bmatrix}^{-1}
\begin{bmatrix}
\Xmat_{\gamma}^T \\ \Xmat_k^T \\
\end{bmatrix}
-\Proj_{\gamma}\\ & = a (I - \Proj_{\gamma}) \Xmat_k \Xmat_k^T (I -
\Proj_{\gamma}),
\end{align*}
which is the claimed decomposition.
\end{proof}

\vspace*{.05in}

Returning to our main task, let us define the event
\begin{align*}
\Aevent'_\numobs \defn \Big \{\max_{\substack{{\gamma \in \Mspace, \,
      k \in \{1, \ldots, \covdim\}} \\ {\mbox{s.t.} \, k \notin
      \gamma}}}
w^T(\Proj_{\gamma \cup \{k\}} - \Proj_{\gamma}) \wnoise \leq L
\sigmazero^2 \log \covdim \Big\}.
\end{align*}
By construction, we have $\Aevent_\numobs' \subseteq \Aevent_\numobs$
so that it suffices to lower bound $\Prob(\Aevent_\numobs')$.
Lemma~\ref{LemProjection} implies that
\begin{align}
\label{EqnNoiseBound}
\wnoise^T(\Proj_{\gamma\cup\{k\}} - \Proj_{\gamma}) \wnoise & =
\frac{\big|\inprod{ \big(I - \Proj_{\gamma}\big)\design_{k}}{
    \wnoise}\big|^2/\numobs} {\design_{k}^T \big(I -
  \Proj_{\gamma}\big) \design_{k}/\numobs}.
\end{align}
Now we show that with probability at least $1-\covdim^{-\UNICON}$, the
above quantity is uniformly bounded by $L\sigmazero^2\log \covdim$
over all $(\gamma,k) \in \Mspace \times \{1,\ldots,\covdim\}$
satisfying $|\gamma| \leq \maxsize$ and $k \notin \gamma$, which
yields the intermediate result
\begin{align}
\label{EqnConcenBound2}
\Prob(\Aevent_n) \geq \Prob(\Aevent_\numobs')\geq 1 -
\covdim^{-\UNICON}.
\end{align}
Now Lemma~\ref{LemUnderfit} implies that $\frac{1}{\numobs}
\design_{k}^T \big(I - \Proj_{\gamma}\big) \design_{k}\geq \Cm$, and
therefore, if we define the random variable
\begin{align*}
V(Z) & \defn \max_{\substack{{\gamma\in \Mspace, \, k \in \{1, \ldots,
      \covdim\}} \\ {\mbox{s.t.}\, k\notin\gamma}}}
\frac{1}{\sqrt{n}}\big|\inprod{ \big(I -
  \Proj_{\gamma}\big)\design_{k}}{Z}\big|, \qquad \mbox{where $Z \sim
  N(0, I_\numobs)$},
\end{align*}
then it suffices to show that $V(Z) \leq \sqrt{L\Cm\log p}$ with
probability at least $1-\covdim^{-\UNICON}$.  For any two vectors $Z,
Z' \in\real^\numobs$, we have
\begin{align*}
|V(Z) - V(Z')| & \leq \max_{\substack{{\gamma\in \Mspace, \, k \in
      \{1, \ldots, \covdim\}} \\ {\mbox{s.t.}\, k\notin\gamma}}}
\frac{1}{\sqrt{n}}\big|\inprod{ \big(I -
  \Proj_{\gamma}\big)\design_{k}}{ Z - Z'}\big| \\ &\leq
\frac{1}{\sqrt{\numobs}} \| \big(I -
\Proj_{\gamma}\big)\design_{k}\|_2 \, \|Z - Z'\|_2 \; \leq \; \|Z -
Z'\|_2,
\end{align*}
where we have used the normalization condition of Assumption B in the
last inequality. Consequently, by concentration of measure for
Lipschitz functions of Gaussian random variables \cite{Ledoux01}, we
have
\begin{align}
\label{EqnLedoux}
\mprob \big[ V(Z) \geq \Exs[V(Z)] + t \big] \leq e^{- \frac{t^2}{2}}.
\end{align}
By the \HACKER condition in Assumption B, the expectation satisfies
$\Exs[V(Z)] \leq \sqrt{L\Cm\log\covdim}\,/2$, which combined with
\eqref{EqnLedoux} yields the claimed bound~\eqref{EqnConcenBound2}
with $\UNICON\leq L\Cm/8$.

%%%%%%%%%%%%%%%%%%%%%%%%%%%%%%%%%%%%%%%%%%%%%%%%%%%%%%%%%%%%%%%%%%%%%

\subsection{Proof of Lemma~\ref{LemPosteriorConcentration}}
\label{SecPosteriorConcentration}

We defer the proof of the claim~\eqref{EqnPosteriorA} to
Appendix~\ref{AppLemForwardSelectionRatio}, as it is somewhat
technically involved. It is worth mentioning that its proof uses some
auxiliary results in Lemma~\ref{LemForwardSelection} in
Appendix~\ref{AppLemForwardSelectionRatio}, which characterizes some
key properties of the state $\Gfun(\gamma)$ selected by the transition
function $\Gfun$ via the greedy criterion.

It remains to prove the second bound~\eqref{EqnPosteriorB} in
Lemma~\ref{LemPosteriorConcentration}, and we split our analysis into
two cases, depending on whether $\gamma$ is underfitted or overfitted.

\subsubsection{Case $\gamma$ is underfitted} 

In this case, the bound~\eqref{EqnPosteriorA} implies that
$\frac{\posterior(\gamma \mid \yobs)}{\posterior(\Gfun(\gamma) \mid
  \yobs)}\leq \covdim^{-3}$.  For each $\gammabar \in \PREC(\gamma)$,
where $\PREC(\gamma)$ is defined in Lemma~\ref{LemPathProperties}(e),
we know $\gamma\in \Path_{\gammabar,\gamma} \subset
\Path_{\gammabar,\gammastar}$. Let the path $\Path_{\gammabar,\gamma}$
be $\gamma_0 \to \gamma_1 \to \cdots \to \gamma_s$, where $s =
|\Path_{\gammabar,\gamma}|$ is the length of
$\Path_{\gammabar,\gamma}$, and $\gamma_0 = \gammabar$ and $\gamma_s =
\gamma$ are the two endpoints. Since any intermediate state
$\gammatil$ on path $\Path_{\gammabar, \gamma}$ is also underfitted,
inequality~\eqref{EqnPosteriorA} ensures that
\begin{align*}
\frac{\posterior(\gammabar \mid \yobs)}{\posterior(\gamma \mid\yobs)}
= \prod_{\ell=1}^s \frac{\posterior(\gamma_{\ell-1} \mid
  \yobs)}{\posterior(\gamma_l \mid \yobs)} \leq \covdim^{-3s} =
\covdim ^{ - 3 \, |\Path_{\gammabar, \gamma}|}.
\end{align*}
Now for each $s \in \{0,\ldots,\sstar\}$, we count the total number of
states $\gammabar$ in $\PREC(\gamma)$ that satisfies
$|\Path_{\gammabar,\gamma}|=s$. By construction, at each intermediate
state in a canonical path, we either add a new influential covariate
by the single flip updating scheme of the MH algorithm, or add a new
influential covariate and delete an unimportant covariate by the
double-flip updating scheme. As a consequence, any state in $\Mspace$
has at most $(\sstar+1)\,\covdim$ adjacent precedents, imlying that
the total number of states $\gammabar$ in $\PREC(\gamma)$ with path
length $|\Path_{\gammabar,\gamma}|=s$ is upper bounded by $(\sstar +
1)^s \,\covdim^s$. Consequently, we have by the preceding display that
under the event $\Aevent_\numobs \cap \Bevent_\numobs \cap
\Cevent_\numobs \cap \Devent_\numobs$
\begin{align}
\label{EqnUnderfit}
\frac{\posterior[\PREC(\gamma)|\yobs]} {\posterior(\gamma \mid \yobs)}
= \sum_{\gammabar \in \mathcal{s}(\gamma)} \frac{\posterior(\gammabar
  \mid \yobs)} {\posterior(\gamma|\yobs)} \leq \sum_{s=0}^{\sstar}
\covdim^s\, (\sstar + 1)^s \, \covdim^{-3s} \leq \sum_{s = 0}^\infty
\covdim^{-s} \leq \frac{1}{1 - 1/\covdim}.
\end{align}
The above argument is also valid for $\gamma=\gammastar$.

%%%%%%%%%%%%%%%%%%%%%%%%%%%%%%%%%%%%%%%%%%%%%%%%%%%%%%%%%%%%%%%

\subsubsection{Case $\gamma$ is overfitted} 

In this case, we bound the ratio $\frac{\posterior[\PREC(\gamma) |
    \yobs]} {\posterior(\gamma | \yobs)}$ by dividing the set
$\PREC(\gamma)$ into two subsets:
\begin{enumerate}
\item[(a)] Overfitted models: $\MODEL_1 = \{\gamma'\in \PREC(\gamma):
  \gamma' \supset \gammastar\}$, all models in $\PREC(\gamma)$ that
  include all influential covariates.
\item[(b)] Underfitted models: $\MODEL_2=\{\gamma'\in \PREC(\gamma):
  \gamma' \not\supset \gammastar\}$, all models in $\PREC(\gamma)$
  that miss at least one influential covariate.
\end{enumerate}
First, we consider the ratio $\posterior(\MODEL_1 \mid
\yobs)/\posterior(\gamma \mid \yobs)$.  For each model $\gammabar\in
\MODEL_1$, according to our construction of the canonical path, all
intermediate states on path $\Path_{\gammabar,\gamma} = \gamma_0 \to
\gamma_1 \to \cdots \to \gamma_k$ correspond to overfitted models
(only involve the first flipping updating scheme of the MH algorithm),
where endpoints $\gamma_0 = \gammabar$ and $\gamma_k = \gamma$, and
$k$ denotes the length of path $\Path_{\gammabar,\gamma}$. As a
consequence, inequality~\eqref{EqnPosteriorA} implies that
\begin{align*}
\frac{\posterior(\gammabar\mid \yobs)} {\posterior(\gamma \mid \yobs)}
& = \prod_{s = 1}^k \frac{\posterior(\gamma_{s-1} \mid \yobs)}
     {\posterior(\gamma_s \mid \yobs)} \leq \covdim^{-2k}.
\end{align*}
Since there are at most $\covdim^k$ states $\gammabar$ in $\MODEL_1$
satisfying $|\gammabar| - |\gamma|=k$, we obtain that under the event
$\Aevent_\numobs \cap \Bevent_\numobs \cap \Cevent_\numobs \cap
\Devent_\numobs$
\begin{equation}
\label{EqnSumM1}
\begin{aligned}
\frac{\posterior(\MODEL_1 \mid \yobs)} {\posterior(\gamma \mid \yobs)}
&\leq \sum_{k=0}^{p- |\gamma|} \covdim^k \,\covdim^{-2k} \leq
\sum_{k=0}^\infty \covdim^{-k} \leq \frac{1}{1 - 1/\covdim} \; \leq 2.
\end{aligned}
\end{equation}

Second, we consider the ratio $\posterior(\MODEL_2 \mid \yobs) /
\posterior(\gamma \mid \yobs)$.  For fixed $\gammabar\in \MODEL_2$,
let $f(\gammabar)$ be the first state along the path
$\Path_{\gammabar,\gamma}$ that contains all influential
covariates. Since the overfitted state $\gamma$ contains all
influential covariates, $f(\gammabar)$ exists and is
well-defined. Moreover,
this construction ensure that $f(\gammabar) \in \MODEL_1$ and
$\gammabar \subset \PREC(f(\gammabar)) \setminus \{f(\gammabar)\}$.
Applying inequality~\eqref{EqnPosteriorA} then yields
\begin{align*}
\frac{\posterior(\MODEL_2 \mid \yobs)} {\posterior(\gamma \mid \yobs)}
= \sum_{\gammabar \in \MODEL_2} \frac{\posterior(\gammabar \mid
  \yobs)} {\posterior(\gamma \mid \yobs)} & = \sum_{\gammabar \in
  \MODEL_2} \frac{\posterior(f(\gammabar) \mid \yobs)}
{\posterior(\gamma \mid \yobs)} \cdot \frac{\posterior(\gammabar \mid
  \yobs)} {\posterior(f(\gammabar) \mid \yobs)} \\
& \leq \sum_{\substack{{\exists \,\gammabar\in \MODEL_2} \\{\mbox{such
        that}\, \gammatil = f(\gammabar)}}} \frac{\posterior(\gammatil
  \mid \yobs)} {\posterior(\gamma \mid \yobs)}\, \sum_{\gammabar \in
  \PREC(\gammatil) \setminus \{\gammatil\} }
\frac{\posterior(\gammabar \mid \yobs)} {\posterior(\gammatil \mid
  \yobs)} \\
& = \sum_{\substack{{\exists \,\gammabar\in \MODEL_2} \\
{\mbox{such that } \, \gammatil = f(\gammabar)}}}
\frac{\posterior(\gammatil \mid \yobs)} {\posterior(\gamma \mid
  \yobs)} \cdot \Big( \frac{\posterior[\PREC(\gammatil) \mid \yobs]}
     {\posterior(\gammatil \mid \yobs)} - 1\Big).
\end{align*}
Then by treating $\gammatil=f(\gammabar)\in \MODEL_1$ as the $\gamma$
in inequality~\eqref{EqnUnderfit} and inequality~\eqref{EqnSumM1}, we
obtain that under the event $\Aevent_\numobs \cap \Bevent_\numobs \cap
\Cevent_\numobs \cap \Devent_\numobs$
\begin{equation}
\label{EqnSumM2}
\begin{aligned}
\frac{\posterior(\MODEL_2 \mid \yobs)} {\posterior(\gamma \mid
  \yobs)}& \leq \sum_{\substack{{\exists \,\gammabar\in
      \MODEL_2}\\{\mbox{s.t.}\, \gammatil = f(\gammabar)}}}
\frac{\posterior(\gammatil \mid \yobs)} {\posterior(\gamma \mid
  \yobs)} \cdot \Big \{ \frac{1}{1 - 1/\covdim} - 1 \Big \} \\
& \leq \frac{2}{\covdim} \sum_{\gammatil \in \MODEL_1}
\frac{\posterior(\gammatil \mid \yobs)}{\posterior(\gamma \mid \yobs)}
\\
& = \frac{2}{\covdim} \, \frac{\posterior(\MODEL_1 \mid \yobs)}
  {\posterior(\gamma \mid \yobs)}\\
& \leq \frac{4}{\covdim}
\end{aligned}
\end{equation}
Combining inequality~\eqref{EqnSumM1} and inequality~\eqref{EqnSumM2},
we obtain that that under the event $\Aevent_\numobs \cap
\Bevent_\numobs \cap \Cevent_\numobs \cap \Devent_\numobs$, the
posterior ratio is upper bounded as
\begin{align}
\label{EqnOverfit}
\frac{\posterior[\PREC(\gamma) \mid \yobs]} {\posterior(\gamma \mid
  \yobs)} & = \frac{\posterior(\MODEL_1 \mid \yobs)}{\posterior(\gamma
  \mid \yobs)} + \frac{\posterior(\MODEL_2 \mid
  \yobs)}{\posterior(\gamma \mid \yobs)} \leq 6.
\end{align}
The above argument is also valid for $\gamma=\gammastar$, and this
completes the proof of inequality~\eqref{EqnPosteriorB}.

%%%%%%%%%%%%%%%%%%%%%%%%%%%%%%%%%%%%%%%%%%%%%%%%%%%%%%%%%%%%%%%%%%%%%%%%%%

\subsection{Proof of Theorem~\ref{ThmBVSconsistency}}

We divide the analysis into two steps. In the first step, we show that
the total posterior probability assigned to models with size
$\mathcal{O}(\sstar)$ other than $\gammastar$ is small. In the second
step, we use the fact that all large models receive small prior
probabilities to show that the remaining models should also receive
small posterior probability.

\paragraph{Step 1:}
Let $\MODEL_S\defn \{\gamma\in\{0,1\}^\covdim:\, |\gamma|\leq
K\sstar,\,\gamma\neq\gammastar\}$ denote the set of all models with
moderate sizes, where $K\geq 1$ some constant to be determined in step
2.  Consider the quantity
\begin{align}
\frac{\posterior(\MODEL_S \mid \yobs)} {\posterior(\gammastar \mid
  \yobs)} = \sum_{\gamma\in \MODEL_S} \frac{\posterior(\gamma \mid
  \yobs)}{\posterior(\gammastar \mid
  \yobs)}. \label{EqnPostDecomposition}
\end{align}

Similar to Lemma~\ref{LemConcentrationIneq}, we modify the definition
of the four events $\Aevent_\numobs,\, \Bevent_\numobs,\,
\Cevent_\numobs$ and $\Devent_\numobs$ by replacing $\Mspace$ with
$\MODEL_S$. Following the proof of Lemma~\ref{LemConcentrationIneq},
it is straightforward to show that these four events satisfy
\begin{align}
\label{EqnEventProb}
\Prob \Big[ \Aevent_\numobs \cap \Bevent_\numobs \cap \Cevent_\numobs
  \cap \Devent_\numobs \Big] & \geq 1 - 6 \covdim^{-\UNICON}.
\end{align}
The following auxiliary lemma ensures that when these four events
hold, then the posterior ratios on the right hand side of
equation~\eqref{EqnPostDecomposition} are well controlled.
%%%%%%%%%%%%%%%%
\blems
\label{LemConsistencyRatio}
Under Assumptions A--D and under the event $\Aevent_\numobs \cap
\Bevent_\numobs \cap \Cevent_\numobs \cap \Devent_\numobs$, the
posterior ratio of any $\gamma$\,($\neq\gammastar$) in $\MODEL_S$ is
bounded as
\begin{align*}
\frac{\posterior(\gamma \mid \yobs)}{\posterior(\gammastar \mid
  \yobs)} & \leq
\begin{cases}
\covdim^{-2 |\gamma\setminus \gammastar|}, & \mbox{if $\gamma$ is
  overfitted},\\
\covdim^{-2|\gamma| - 2}, & \mbox{if $\gamma$ is underfitted}.
\end{cases}
\end{align*}
\elems
%%%%%%%%%%%%%%
\noindent We prove this lemma in
Appendix~\ref{SectionLemConsistencyRatio}.

Equipped with this lemma, a simple counting argument yields that under
the event $\Aevent_\numobs \cap \Bevent_\numobs \cap \Cevent_\numobs
\cap \Devent_\numobs$,
\begin{align*}
\frac{\posterior(\MODEL_S \mid \yobs)} {\posterior(\gammastar \mid
  \yobs)} \overset{(i)}{\leq} &\, \sum_{k = 1}^\infty \covdim^{k}
\covdim^{-2k} + \sum_{\ell=0}^\infty \covdim^{l} \covdim^{-2l - 2}
\leq 3\covdim^{-1},
\end{align*}
where in step (i), we used the fact that there are at most
$\covdim^{k}$ overfitted models $\gamma$ with $|\gamma\setminus
\gammastar| = k$ and at most $\covdim^{\ell}$ underfitted models
$\gamma$ with $|\gamma| = \ell$.  Combining this with
inequality~\eqref{EqnEventProb}, we obtain that with probability at
least $1 - 6 \covdim^{-\UNICON}$,
\begin{align}
\label{EqnConsistencyA}
\posterior(\MODEL_S \mid \yobs) \leq 3\covdim^{-1} \,
\posterior(\gammastar \mid \yobs) \leq 3\covdim^{-1}.
\end{align}
%%%%%%%

\paragraph{Step 2:} 
Let $\MODEL_L\defn \{\gamma\in\{0,1\}^\covdim:\, |\gamma| \geq K\sstar
+ 1\}$ denote the set of large models. By Bayes' theorem, we can
express the posterior probability of $\MODEL_L$ as
\begin{align}
\posterior (\MODEL_L \mid Y) = \frac{\sum_{\gamma \in \MODEL_L}
  \int_{\theta,\phi} \frac{ d\Prob_{\beta, \phi, \gamma} }{ d\Prob_{0}
  } (Y) \, \posterior(d\theta, d\phi, \gamma)}{\sum_{\gamma \in
    \{0,1\}^\covdim} \int_{\theta,\phi} \frac{ d\Prob_{\beta, \phi,
      \gamma} }{ d\Prob_{0} } (Y) \, \posterior(d\theta, d\phi,
  \gamma)}, \label{EqnPosteriorL}
\end{align}
where $\Prob_{\beta, \phi, \gamma}$ and $\Prob_{0}$ stand for
probability distribution of $Y$ under parameters
$(\beta,\phi,\gamma)$ and the true data generating model,
respectively. We bound the numerator and denominator separately.

First consider the numerator.  According to our specification of the
sparsity prior~\eqref{EqnPriorGamma} for the binary indicator vector
$\gamma$, the prior probability of $\MODEL_L$ satisfies
\begin{align*}
\posterior (\MODEL_L) = \sum_{\gamma:\, |\gamma| > K\sstar + 1}
\posterior (\gamma) \leq \covdim^{-K \sstar - 1}.
\end{align*}
By Fubini's theorem we have the following bound for the expectation
of the numerator:
\begin{align*}
\Exs_0 \Big[ \sum_{\gamma \in \MODEL_L} \int_{\theta,\phi} \frac{
    d\Prob_{\beta, \phi, \gamma} }{ d\Prob_{0} } (Y)
  \posterior(d\theta, d\phi, \gamma) \Big] = &\, \sum_{\gamma \in
  \MODEL_L} \int_{\theta,\phi} \Exs_0 \Big[ \frac{d\Prob_{\beta, \phi,
      \gamma} }{ d\Prob_{0} } (Y) \Big] \posterior(d\theta, d\phi,
\gamma) \\
= & \, \sum_{\gamma \in \MODEL_L} \int_{\theta,\phi}
\posterior(d\theta, d\phi, \gamma) = \posterior(\MODEL_L) \leq
\covdim^{-K \sstar - 1},
\end{align*}
where we have used the fact that $\Exs_0 \big[
\frac{d\Prob_{\beta, \phi, \gamma} }{ d\Prob_{0} } (Y)\big] = 1$.
Therefore, by applying Markov's inequality we have
\begin{align}\label{EqnPosteriorNum}
\Prob_0 \Big[ \sum_{\gamma \in \MODEL_L} \int_{\theta,\phi} \frac{
    d\Prob_{\beta, \phi, \gamma} }{ d\Prob_{0} } (Y)
  \posterior(d\theta, d\phi, \gamma) \leq \covdim^{-K \sstar/2 - 1}
  \Big] \geq 1 - \covdim^{-K \sstar/2 }.
\end{align}

By the expression~\eqref{EqnMarginalLIK} of the marginal likelihood
function, we can bound the denominator from below by
\begin{align*}
& \int_{\theta,\phi} \frac{ d\Prob_{\beta, \phi, \gammastar} }{
    d\Prob_{0} } (Y) \, \posterior(d\theta, d\phi, \gammastar) =
  \frac{\mathcal{L}_\numobs(\yobs | \, \gammastar) \, \posterior
    (\gammastar)}{ d\Prob_{0} (Y)} \\
= & \, \frac{\Gamma\big( \frac{n}{2} \big)\, (1 +
  \hyperpara)^{\numobs/2}} {\pi^{\numobs/2} } \frac{(1 + \hyperpara)^{
    -\sstar/2}}{( \|\yobs\|_2^2 + \hyperpara\, \| ( I -
  \Proj_{\gammastar} ) \, \wtil \|_2^2 )^{n/2}} \cdot \frac{ \UNICON
  \, \covdim^{ - 2 \sstar}} { d\Prob_{0} (Y)},
\end{align*}
where $\wtil = \wnoise + \design_{\Sset^c}\betastar_{\Sset^c} \sim
\mathcal{N} (\design_{\Sset^c}\betastar_{\Sset^c}, \sigmazero^2)$.
Under the true data-generating model $\Prob_0$, the density for
$\yobs$ is $\sigmazero^{-\numobs} (2\pi)^{-\numobs/2}
\exp\{-(2\sigmazero^2)^{-1} \|\wnoise\|_2^2\}$. By applying the the
lower bound \mbox{$\Gamma(\numobs/2) \geq (2\pi)^{1/2}\, (\numobs/2 -
  1)^{\numobs /2 - 1/2} e^{-\numobs/2 + 1}$} and using the fact that
the projection operator $I - \Proj_{\gammastar} $ is non-expansive,
we obtain
\begin{align*}
& \int_{\theta,\phi} \frac{ d\Prob_{\beta, \phi, \gammastar} }{
    d\Prob_{0} } (Y) \, \posterior(d\theta, d\phi, \gammastar) \geq
  \UNICON \, \covdim^{ - 2 \sstar} (1 + \hyperpara)^{ -\sstar/2} (1 +
  \hyperpara^{-1})^{\numobs/2} \\ &\qquad\qquad\qquad\qquad \cdot
  \exp\big\{(2\sigmazero^2)^{-1} \big(\|\wnoise\|_2^2 - \|\wtil\|_2^2
  - \|\yobs\|_2^2/\hyperpara \big) \big\} \underbrace{\big(
    u^{-\numobs/2} e^{u/2} \big)}_{f(u)} \; \underbrace{\big(
    \numobs^{\numobs/2} e^{-\numobs/2}\big)}_{1/f(\numobs)}
\end{align*}
where $u = \sigmazero^{-2} (\|\wtil\|_2^2 +
\|\yobs\|_2^2/\hyperpara)$. Since $\hyperpara^{-1} \lesssim
\numobs^{-1}$ and the function $f(u)=u^{-\numobs/2} e^{u/2}$ attains
its minimum at $u = \numobs$, we further obtain
\begin{align*}
& \int_{\theta,\phi} \frac{ d\Prob_{\beta, \phi, \gammastar} }{
    d\Prob_{0} } (Y) \, \posterior(d\theta, d\phi, \gammastar) \geq
  \UNICON \, \covdim^{ - 2 \sstar} (1 + \hyperpara)^{ -\sstar/2} \,
  \exp\big\{(2\sigmazero^2)^{-1} \big(\|\wnoise\|_2^2 - \|\wtil\|_2^2
  - \|\yobs\|_2^2/\hyperpara \big) \big\},
\end{align*}
with a different universal constant $\UNICON$.

The off-support $\Sset^c$ condition in Assumption A and the high probability bound for the event $\Cevent_\numobs \cap\Devent_\numobs$ in Lemma~\ref{LemConcentrationIneq} imply that the last exponential term is of order
$\covdim^{-\UNICON_1 \Ltil}$ for some universal constant $\UNICON'$
with probability at least $1-\covdim^{-\UNICON_2}$. Therefore, for $K
\geq 4 + \alpha + 2\UNICON_1\Ltil$, we have
\begin{align}\label{EqnPosteriorDen}
& \int_{\theta,\phi} \frac{ d\Prob_{\beta, \phi, \gammastar} }{ d\Prob_{0} } (Y)  \,
\posterior(d\theta, d\phi, \gammastar) 
\geq  \UNICON \, \covdim^{ - K\sstar/2}.
\end{align}
Combining equations~\eqref{EqnPosteriorL}, \eqref{EqnPosteriorNum}
and~\eqref{EqnPosteriorDen}, we obtain that
\begin{align}
\label{EqnConsistencyB}
\posterior (\MODEL_L \mid Y) \leq \UNICON\, \covdim^{- 1}
\end{align}
holds with probability at least $1 - 2\, \covdim^{-\UNICON'}$.

Finally, inequalities~\eqref{EqnConsistencyA} and
\eqref{EqnConsistencyB} in steps 1 and 2 together yield that
\begin{align*}
\posterior (\gammastar \mid Y) = 1 - \posterior (\MODEL_S \mid Y) -
\posterior (\MODEL_L \mid Y) \geq 1 - \UNICON_3\, \covdim^{- 1},
\end{align*}
holds with probability at least $1 - 8\, \covdim^{-\UNICON'}$, which
completes the proof.

%%%%%%%%%%%%%%%%%%%%%%%%%%%%%%%%%%%%%%%%%%%%%%%%%%%%%%%%%%%%%%%%%%%%%%%%%%%

\subsection{Proof of Corollary~\ref{CoroMCMC}}
Let $\Prob_t$ denote the probability distribution of iterate
$\gamma_t$ in the MCMC algorithm.  According to the definition of
$\epsilon$-mixing time, for any $t \geq \tau_{1/\covdim}$, we are
guaranteed that $\big| \Prob_t (\gammastar) - \posterior (\gammastar)
\big| \leq \frac{1}{\covdim}$.  By Theorem~\ref{ThmBVSconsistency},
the posterior probability of $\gammastar$ satisfies $\posterior
(\gammastar) \geq 1 - \UNICON_1 \, \covdim^{-1}$ with probability at
least $1 - \UNICON_2 \,\covdim^{ - \UNICON_3}$.  By
Theorem~\ref{ThmMain}, the $\covdim^{-1}$-mixing time $\tau_{1 /
  \covdim}$ satisfies
\begin{align*}
\tau_{1 / \covdim} \leq c_1 \, \covdim \maxsize^2 \, \big (c_2 \alpha \,
( \numobs + \maxsize) \log \covdim + \log \covdim  + 2 \big)
\end{align*}
with probability at least $1-4\covdim^{-c_1}$.  Combining the three
preceding displays, we find that $\Prob_t (\gammastar) \geq 1 -
(\UNICON_1 + 1) \, \covdim^{-1}$, as claimed.

%%%%%%%%%%%%%%%%%%%%%%%%%%%%%%%%%%%%%%%%%%%%%%%%%%%%%%%%%%%%%%%%%%%%%%%%%%%

\section{Discussion}
\label{SecDiscussion}

In this paper, we studied the computational complexity of MCMC methods
for high-dimensional Bayesian linear regression under a sparsity
constraint. We show that under a set of conditions that guarantees
Bayesian variable-selection consistency, the corresponding MCMC
algorithm achieves rapid mixing. Our result on the computational
complexity of Bayesian variable-selection example provides insight
into the dynamics of the Markov chain methods applied to statistical
models with good asymptotic properties. It suggests that contraction
properties of the posterior distribution are useful not only in
guaranteeing desirable statistical properties such as parameter
estimation or model selection consistency, but they also have
algorithmic benefits in certifying the rapid mixing of the Markov
chain methods designed to draw samples from the posterior.

As a future direction, it is interesting to investigate the mixing
behavior of the MCMC algorithm when Bayesian variable selection
fails. For example, slow mixing behavior is observed empirically in
the intermediate SNR regime in our simulated example and it would be
interesting to understand this result theoretically. Another
interesting direction is to consider the computational complexity of
MCMC methods for models more complex than linear regression, for
example, high-dimensional nonparametric additive regression. A third
direction is to investigate whether the upper bound on mixing time
provided in Theorem~\ref{ThmMain} is sharp up to constants.

\subsection*{Acknowledgements}

Authors YY, MJW and MIJ were partially supported by Office of Naval
Research MURI grant N00014-11-1-0688.  YY and MJW were additionally
supported by National Science Foundation Grants CIF-31712-23800 and 
DMS-1107000.

%%%%%%%%%%%%%%%%%%%%%%%%%%%%%%%%%%%%%%%%%%%%%%%%%%%%%%%%%%%%%%%%%%%%%%%%

\appendix

%%%%%%%%%%%%%%%%%%%%%%%%%%%%%%%%%%%%%%%%%%%%%%%%%%%%%%%%%%%%%%%%%%%%%%%%%%%

\section{Further details on Metropolis-Hastings}
\label{AppMAP}

In Appendix~\ref{AppMAP}, we show that under the specified model, the
maximum a posteriori solution (MAP) of the Bayesian variable-selection
problem is equivalent to the following optimizatio problem with
$\ell_0$-penalty
\begin{align*}
\widehat{\gamma} = \arg \min_{|\gamma|\leq \maxsize} \Big \{
\frac{\numobs}{2}\log \Big [ 1 + \hyperpara \big(1-
  \frac{\yobs^T\Proj_{\gamma} \yobs}{\|\yobs\|_2^2} \big)
  \Big] + \lambda |\gamma| \Big \},
\end{align*} 
where $\Proj_{\gamma} = \design_{\gamma} (\design_{\gamma}^T
\design_{\gamma})^{-1} \design_{\gamma}^T$ is the projection onto the
column space of $\design_{\gamma}$, and the regularization parameter
\mbox{$\regu \defn \frac{1}{2} \log (1+\hyperpara) + \kappa \log \covdim$.}
Here the penalty $\regu |\gamma|$ comes from two sources:
the penalty $\kappa \log \covdim\,|\gamma|$ on $\gamma$
and the Occam's razor penalty $\frac{1}{2}\log(1+\hyperpara)\,|\gamma|$ 
due to the integration over the model parameter
$\beta_{\gamma}$. Therefore, choosing an appropriate hyperparameter
$\kappa$ in the Bayesian approach is equivalent to choosing a
corresponding regularization parameter $\regu$ in the penalization method: a
small $\kappa$ could make the posterior include uninfluential
covariates due to noise; a large $\kappa$ requires the
signal-to-noise ratio $\betastar_j/\sigma$ of influential covariates
to be large enough so that they can be selected out by the posterior.

%%%%%%%%%%%%%%%%%%%%%%%%%%%%%%%%%%%%%%%%%%%%%%%%%%%%%%%%%%%%%%%%%%%%%%%%%%
\subsection{Connections to MAP estimates}

Under the Bayesian model specified by
equation~\eqref{EqnHierarchicalBayes}, we can obtain a closed-form
expression for the marginal likelihood of the indicator vector
$\gamma$ by integrating out $\beta_{\gamma}$ and $\phi$:
\begin{align}
\mathcal{L}_\numobs(\yobs | \, \gamma) \defn \posterior( Y | \,
\gamma)= &\,\int d\Prob_{\beta, \phi, \gamma} (Y)\, \posterior(
d\beta, d\phi \mid \gamma) \notag\\ = &\, \frac{\Gamma\big(
  \frac{n}{2} \big)\, (1 + \hyperpara)^{\numobs/2} } {\pi^{\numobs/2 }
  \,\|\yobs\|_2^\numobs} \frac{(1 + \hyperpara)^{ -|\gamma|/2}}{(1 +
  \hyperpara (1 - R_{\gamma}^2) )^{n/2}}, \label{EqnMarginalLIK}
\end{align}
where $\Gamma(\cdot)$ the Gamma function, $\Prob_{\beta, \phi, \gamma} $ is the 
distribution of $Y$ under parameters $(\beta,\phi,\gamma)$, and $R^2_{\gamma}$ is the coefficient of determination for the model
$\MODEL_{\gamma}$
\begin{align*}
R^2_{\gamma} = \frac{\yobs^T\Proj_{\gamma} \yobs}{\|\yobs\|_2^2},
\end{align*}
with $\Proj_{\gamma} = \design_{\gamma} (\design_{\gamma}^T
\design_{\gamma})^{-1} \design_{\gamma}^T$ the projection onto the
column space of $\design_{\gamma}$. When there is no confusion, we
identify the variable inclusion vector $\gamma$ and the linear model
$\MODEL_{\gamma}$ associated with it.  Let $\Mspace \defn \{\gamma:\,
|\gamma|\leq \maxsize\}$ denote the entire model space, which is a
subset of the $\covdim$-dimensional hypercube $\{0,1\}^\covdim$ under our
identification. Then, by Bayes' theorem the posterior probability of
$\gamma$ is given by
\begin{align}
\label{EqnPostFormula}
\posterior(\gamma \mid \yobs) = C \cdot\frac{1}{\covdim^{\kappa |\gamma|}}
\cdot \frac{(1+\hyperpara)^{-|\gamma|/2}}{ (1 + \hyperpara (1 -
  R_{\gamma}^2) )^{\numobs/2}} \, \Ind[\gamma\in\Mspace],
\end{align}
where $C$ is a normalization constant.

According to the preceeding display, the maximum a posteriori (MAP) 
solution of the Bayesian variable-selection problem is equivalent to 
the following penalized optimization problem
$\ell_0$-penalty
\begin{align*}
\widehat{\gamma} = \arg \min_{|\gamma|\leq \maxsize} \Big \{
\frac{\numobs}{2}\log \Big [ 1 + \hyperpara \big(1-
  \frac{\yobs^T\Proj_{\gamma} \yobs}{\|\yobs\|_2^2} \big)
  \Big] + \lambda |\gamma| \Big \},
\end{align*} 
where the regularization parameter
\mbox{$\regu \defn \frac{1}{2} \log (1+\hyperpara) + \kappa \log \covdim$.}
Conversely, if we have a variable-selection
procedure based on the penalization method with $\ell_0$-penalty
\begin{align*}
\widehat{\gamma} = \mbox{arg}\min_{\gamma:\, |\gamma|\leq \maxsize}
 \big\{f(\yobs,\gamma) + \regu |\gamma| \big\},
\end{align*}
where $f(\yobs,\gamma)$ is some function reflecting the goodness of
fit by using model $\MODEL_{\gamma}$, then we can construct a
pseudo-posterior distribution
\begin{align*}
\widetilde{\pi}_\numobs(\gamma|\yobs) = \widetilde{C} \cdot e^{-\regu
  |\gamma|} \cdot e^{- f(\yobs,\gamma)} \, \Ind \big[ |\gamma| \leq
  \maxsize \big]
\end{align*}
with $\widetilde{C}$ is a normalization constant, and conduct
Bayesian inference based on $\widetilde{\pi}_n$. For example, when
$f(\yobs,\gamma)$ is the negative profile log-likelihood
$\frac{\numobs}{2}
\log \big( \frac{1}{\numobs} \|( I - \Proj_{\gamma}) \yobs\|_2^2 \big)$, where
the regression coefficient $\beta_{\gamma}$ and the precision
parameter $\phi$ have been profiled out from the log-likelihood given
$\gamma$ by maximization, then the choice of $\regu = 1$
corresponds to the Akaike information criterion, or AIC for
short~\cite{Akaike1974}), $\regu=\log \numobs/2$ the Bayesian
information criterion, or BIC for short~\cite{Schwarz1978}), and
$\regu = \alpha \log \covdim$ ($\alpha\geq 1$) the high-dimensional
BIC~\cite{Wang11}. Our conclusion on the MCMC complexity of Bayesian
variable selection with Zellner's $\hyperpara$-prior applies
to BIC in the low-dimensional regime where
$\covdim=o(\numobs)$, and to high-dimensional BIC in the
high-dimensional regime where $\covdim\geq\numobs$. Because of the
Bayesian interpretation for $\posterior(\gamma|\yobs)$ in
equation~\eqref{EqnPostFormula}, we will focus on this posterior
distribution over the model space $\Mspace$.  

%%%%%%%%%%%%%%%%%%%%%%%%%%%%%%%%%%%%%%%%%%%%%%%%%%%%%%%%%%%%%%%%%%%%%

\subsection{Example of slow mixing}

\label{AppSlowMixing}

Suppose $\covdim=\numobs$ and $\hyperpara = \covdim^{2\alpha}$ with
$\alpha > 1$. Let $\yobs = \wnoise \sim \neigh(0,I_\numobs)$. We claim
that if we use the untruncated distribution $\posterior(\gamma) =
Cp^{-\kappa |\gamma|}$ as the prior for the variable-selection
indicator vector over the entire space $\{0,1\}^\covdim$, then the
mixing time of the Markov chain with transition probability specified
by formula~\eqref{EqnMetropolisHastings} grows exponentially in
$\numobs$ with probability at least $1/2$ with respect to the
randomness of $\wnoise$. Moreover, it is easy to check that this
example satisfies the conditions in Theorem~\ref{ThmBVSconsistency},
which imply Bayesian variable-selection consistency. As a consequence,
this example suggests that although a size constraint $|\gamma| \leq
\maxsize$ such as the one in the sparsity prior~\eqref{EqnPriorGamma}
is not needed for Bayesian model selection consistency, it is
necessary for MCMC to mix rapidly.

\paragraph{Proof of slow mixing:}
We use the following conductance argument: for any reversible Markov
chain $\MarkovChain$ over a finite state space, the spectral gap is
upper bounded as
\begin{align}
\label{EqnConductance}
1-\lambda_2\leq 2 \Phi_\MarkovChain,
\end{align}
where the quantity
\begin{align}
\label{EqnDefnConductance}
\Phi_\MarkovChain & \defn \min_{A \subset \Mspace:\,
  0<\pi(B)<1}\Phi_\MarkovChain(A), \quad \mbox{where
  $\Phi_\MarkovChain(A) \defn \frac{\sum_{\gamma\in A} \pi(\gamma)
    P(\gamma,A^c)}{\pi(A) \,\pi(A^c)}$}
\end{align}
is called the conductance of $\MarkovChain$~\cite{Sinclair1992}.

Now we analyze the mixing time of the Markov chain in the previous
example. Use the notation $\mathbf{1}$ to denote the full model. Under
the prior choice in the theorem, the posterior has an expression as
\begin{align}
\label{EqnPostP}
\posterior(\gamma|\yobs) \propto \frac{1} {\covdim^{|\gamma|}} \cdot
\frac{(1 + \hyperpara)^{-|\gamma|/2}}  {(1 +
  \hyperpara(1 - R_{\gamma}^2))^{\numobs/2}  }  \quad \mbox{for $\gamma \in
  \{0,1\}^\covdim$.}
\end{align}
Now we apply
inequality~\eqref{EqnConductance} and
equation~\eqref{EqnDefnConductance} with $B = \{\mathbf{1}\}$ to obtain
\begin{align*}
1-\lambda_2 \leq  2\Phi_{\MarkovChain} \leq
\frac{2\sum_{i=1}^\numobs  \posterior(\mathbf{1} \mid \yobs)  P(\mathbf{1},
  \mathbf{1}_{-j})} {\posterior(\mathbf{1} \mid \yobs)\,
  (1 - \posterior(\mathbf{1} \mid \yobs))} = \frac{2\sum_{i=1}^\numobs
  P(\mathbf{1}, \mathbf{1}_{-j})}  {1-\posterior(\mathbf{1} \mid \yobs)},
\end{align*}
where we have used the fact that under the transition probability specification~\eqref{EqnMetropolisHastings}, the only ``neighbor'' of $\mathbf{1}$ is
$\mathbf{1}_{-j}$ for $j=1,\ldots,\numobs$, i.e.\!
$\gamma $\,($ \neq \mathbf{1}$) satisfies $P(\mathbf{1}, \gamma)>0$ if
and only if $\gamma = \mathbf{1}_{-j}$ for some
$j \in \{1,\ldots,\numobs\}$.  Using~\eqref{EqnMetropolisHastings} and the last
display, we can further obtain
\begin{align}
\label{EqnGap}
1-\lambda_2\leq \frac{\sum_{i=1}^\numobs \min \big \{1,
  \frac{\posterior(\mathbf{1}_{-j} \mid \yobs)}  {\posterior(\mathbf{1}
    \mid \yobs)}\big\}}{\numobs (1 - \posterior(\mathbf{1} \mid
  \yobs))}.
\end{align}

We consider the numerator of the right-hand side in equation~\eqref{EqnGap}
first. Since the true model is the null model, we have
\begin{align}
\label{EqnFullRatio}
\frac{\posterior(\mathbf{1}_{-j} \mid \yobs)}{ \posterior ( \mathbf{1}
  \mid \yobs)} = \frac{\numobs^{\kappa} ( 1 + \hyperpara)^{1/2}}{\big ( 1 +
  \hyperpara \wnoise^T ( I - P_{-j}) \wnoise/ (\wnoise^T \wnoise)
  \big)^{\numobs/2}},
\end{align}
where $P_{-j}$ is the projection onto $\design_{-j}$. Since
$\{\wnoise^T(I - P_{-j})\wnoise\}_{j=1}^\numobs$ are $\chi^2$ random
variables, by the union bound and the tail probability of $\chi^2$
distribution, we have that for constant $\UNICON_1$ sufficiently
small,
\begin{align*}
\Prob\Big( \min_{j=1,\ldots,n}  \wnoise^T (I - P_{-j}) \wnoise \geq
\frac{\UNICON_1} {\numobs^2} \Big) \geq 1-
\numobs  \cdot  \frac{1}{4\numobs} =\frac{3}{4}.
\end{align*}
Moreover, by a standard tail bound for the $\chi_\numobs^2$ distribution
(\cite{Laurent2000}, Lemma 1) we have
\begin{align*}
\Prob \Big(\wnoise^T\wnoise \leq \frac{3}{2}\numobs\Big)\geq
\frac{3}{4}.
\end{align*}
Combining the last two displays, we obtain
\begin{align*}
\Prob\Big( \min_{j=1,\ldots,\numobs}  \frac{\wnoise^T ( I - P_{-j})
  \wnoise} {\|\wnoise\|_2^2} \geq \frac{\UNICON_2} {\numobs^3}\Big) \geq
\frac{1}{2},
\end{align*}
where $\UNICON_2 = 2\UNICON_1/3$. Combining the above with
equation~\eqref{EqnFullRatio}, we obtain
\begin{align}
\label{EqnGapP1}
\Prob\Big( \min_{j=1,\ldots,\numobs}\frac{\posterior(\mathbf{1}_{-j}
  \mid \yobs)} {\posterior(\mathbf{1} \mid \yobs)} \leq
e^{-\UNICON_3\hyperpara \numobs^{-2}} \Big) \geq \frac{1}{2},
\end{align}
where we have used the inequality $1 + x \leq  e^x$ for $x \in \real$ and
$\UNICON_3$ is some universal constant.

Now we consider the denominator of the right-hand side in \eqref{EqnGap}. Recall that
$\mathbf{0}$ denotes the indicator vector in $\{0,1\}^\covdim$
associated with the null model.  By equation~\eqref{EqnPostP}, we have
\begin{align*}
\frac{\posterior(\mathbf{1} \mid \yobs)}{\posterior(\mathbf{0} \mid
  \yobs)} = \frac{(1 + \hyperpara)^{-\numobs/2}}  {\numobs^{\kappa \numobs}} \Big/
\frac{1}{(1 + \hyperpara)^{\numobs/2}} = \numobs^{-\kappa \numobs} \leq
\frac{1}{2},
\end{align*}
for $\numobs\geq 2$. This implies
\begin{align}
\label{EqnGapP2}
\posterior(\mathbf{1}|\yobs) = \frac{\posterior(\mathbf{1} \mid
  \yobs)} {\sum_{\gamma\in \Mspace} \posterior(\gamma \mid \yobs)} \leq
\frac{\posterior(\mathbf{1} \mid \yobs)}{\posterior(\mathbf{0} \mid
  \yobs)} \leq \frac{1}{2}.
\end{align}

Combining equations~\eqref{EqnGap}, ~\eqref{EqnGapP1}
and~\eqref{EqnGapP2} yields
\begin{align*}
\Prob \Big(1- \lambda_2\leq 2 e^{-\UNICON_3 \hyperpara \numobs^{-2}}
\Big) \geq \frac{1}{2},
\end{align*}
which completes the proof of the claimed result.

%%%%%%%%%%%%%%%%%%%%%%%%%%%%%%%%%%%%%%%%%%%%%%%%%%%%%%%%%%%%%%%%%%%%%%%

\section{Proof of inequality~\eqref{EqnPosteriorA} in Lemma~\ref{LemPosteriorConcentration}}
\label{AppLemForwardSelectionRatio}

Since $\gamma\neq\gammastar$, we know that $\gamma'\defn
\Gfun(\gamma)\neq \gamma$.  We divide the proof into the following
three disjoint cases:
\begin{itemize}
\item model $\gamma$ is overfitted,
\item model $\gamma$ is underfitted and unsaturated,
\item model $\gamma$ is underfitted and saturated.
\end{itemize}

\subsection{Case $\gamma$ is overfitted}

Let $\ell_{\gamma}$ be the index selected from the set
$\gamma\setminus\gammastar$ of uninfluential covariates in our
construction of the transition function $\Gfun$,
i.e. $\gamma'=\gamma\setminus\{\ell_{\gamma}\}$. We can express the
posterior probability ratio as
\begin{align}
\frac{\posterior(\gamma \mid \yobs)}{\posterior(\gamma' \mid \yobs)} &
= \frac{1}{ \covdim^\kappa \sqrt{1 + \hyperpara}} \cdot
\Big(\frac{1 + \hyperpara(1-R_{\gamma'}^2)}{1 + \hyperpara ( 1 -
  R_{\gamma}^2)}\Big)^{\numobs/2} \notag \\
& = \frac{1}{\covdim^\kappa \sqrt{1 + \hyperpara}} \cdot \Big(1 +
\frac{R_{\gamma}^2 - R_{\gamma'}^2}{g^{-1} + (1 - R_{\gamma}^2)}\Big)^{n/2}
\notag.
\end{align}
Since all influential covariates are included in models
$\MODEL_{\gamma}$ and $\MODEL_{\gamma'}$, we have
\begin{align*}
1 - R_{\gamma}^2 = \frac{\|(I - \Proj_{\gamma})
   (X_{(\gammastar)^c}\betastar_{(\gammastar)^c} + \wnoise)  \|_2^2}   
  {\|\yobs\|_2^2} \quad \mbox{and} \quad
1 - R_{\gamma'}^2 = \frac{\|(I - \Proj_{\gamma'})  (X_{(\gammastar)^c}
  \betastar_{(\gammastar)^c} + \wnoise)\|_2^2}{\|\yobs\|_2^2}.
\end{align*}
Applying the Cauchy-Schwarz inequality yields
\begin{align*}
1 - R_{\gamma}^2 & \geq \frac{ \frac{1}{2} \|(I - \Proj_{\gamma})
  \wnoise\|_2^2 - \|(I - \Proj_{\gamma}) \Xmat_{(\gammastar)^c}
  \betastar_{(\gammastar)^c} \|_2^2}{\|\yobs\|_2^2} \\
& \overset{(i)}{\geq} \frac{\| (I - \Proj_{\gamma}) \wnoise\|_2^2 -
  2\|X_{(\gammastar)^c}
  \betastar_{(\gammastar)^c} \|_2^2}{2\|\yobs\|_2^2}
\overset{(ii)}{\geq} \frac{\|( I - \Proj_{\gamma}) \wnoise\|_2^2 -
  2\tilde{L} \sigmazero^2 \log \covdim}{2\|\yobs\|_2^2},
\end{align*}
where in step (i) we used the fact that the projection is a non-expansive
mapping and in step (ii) we used the assumption $\|X_{(\gammastar)^c}
\betastar_{(\gammastar)^c}\|_2^2\leq \tilde{L} \sigmazero^2 \log
\covdim$.  Similarly, since $\gamma' \subset \gamma$, we obtain
the following inequality for the quantity $R_{\gamma}^2 - R_{\gamma'}^2$
\begin{align*}
R_{\gamma}^2 - R_{\gamma'}^2 & = \frac{\|(\Proj_{\gamma} -
  \Proj_{\gamma'}) (X_{(\gammastar)^c} \betastar_{(\gammastar)^c} +
  \wnoise)\|_2^2}  {\|\yobs\|_2^2} \leq \frac{2\| (\Proj_{\gamma} -
  \Proj_{\gamma'}) \wnoise\|_2^2 + 2 \Ltil
  \sigmazero^2 \log \covdim} {\|\yobs\|_2^2}.
\end{align*}
On the event $\Aevent_\numobs \cap \Bevent_\numobs \cap
\Cevent_\numobs$, we have 
\begin{align*}
\|(I - \Proj_{\gamma}) \wnoise\|_2^2\geq \frac{1}{2}\numobs - r
\maxsize \log \covdim \geq \frac{3 \numobs}{8}, \quad \mbox{and} \quad
\|(\Proj_{\gamma} - \Proj_{\gamma'}) \wnoise\|_2^2 \leq (|\gamma | -
|\gamma'|) L \sigmazero^2 \log \covdim = L \sigmazero^2 \log \covdim,
\end{align*}
where we have used the assumption $4 r \, \maxsize \log \covdim \leq
\numobs $.  Combining these two inequalities with the preceding two
displays, we obtain that the posterior probability ratio on the event
$\Aevent_\numobs \cap \Bevent_\numobs \cap \Cevent_\numobs$ is bounded
as
\begin{align}
\frac{\posterior(\gamma \mid \yobs)}{\posterior(\gamma' \mid \yobs)} &
\leq\frac{1}{\covdim^\kappa \sqrt{g}} \cdot \Big( 1 + \frac{2 (L + \Ltil )
  \log \covdim}{\numobs/8})^{\numobs/2} \notag \\
\label{EqnPosRatio}
 & \overset{(i)}{\leq} \frac{1}{\covdim^\kappa \sqrt{g}} \cdot \exp \Big\{
\frac{16 ( L + \Ltil )\log \covdim}{\numobs} \, \frac{\numobs}{2}
\Big\} = \covdim^{8 (L + \Ltil) - 1 - \alpha - \kappa} \overset{(ii)}{\leq}
\covdim^{-2},
\end{align}
where in step (i) we used the inquality $1+x \leq e^x$ for $x\in\real$
and the last step follows since $\sqrt{\hyperpara}\asymp
\covdim^{\alpha}$ with $\alpha \geq 8(L + \tilde{L}) + 1 - \kappa$
according to our choice of the hyperparameter, which completes the
proof of the overfitted case.

%%%%%%%%%%%%%%%%%%%%%%%%%%%%%%%%%%%%%%%%%%%%%%%%%%%%%%%%%%%%%%%%%%%%

\subsection{Case $\gamma$ is underfitted and unsaturated}

This case happens only when $\sstar\geq 1$.  Let $j_{\gamma}$ be the
index in our construction of the transition function $\Gfun$, i.e. the
index from the set $\gammastar\setminus\gamma$ that maximizes
$\|\Proj_{\gamma\cup\{j\}} \design_{\gammastar} \betastar_{\gammastar}
\|_2^2$ over $j \in \gammastar \setminus \gamma$. Then we have
$\gamma' = \gamma \cup \{j_{\gamma}\}$, implying that $\Proj_{\gamma'}
- \Proj_{\gamma}$ is a projection operator. Therefore, we can write
\begin{align*}
1 - R_{\gamma}^2 - ( 1 - R_{\gamma'}^2) & = \frac{\yobs^T (
  \Proj_{\gamma'} - \Proj_{\gamma}) \yobs}{\|\yobs\|_2^2} \\
& = \frac{\big\| (\Proj_{\gamma'} - \Proj_{\gamma})
  \design_{\gammastar} \beta_{\gammastar}^\ast + (\Proj_{\gamma'} -
  \Proj_{\gamma}) X_{(\gammastar)^c} \betastar_{(\gammastar)^c} +
  (\Proj_{\gamma'} - \Proj_{\gamma}) \wnoise \big\|_2^2}{\|\yobs\|_2^2}
\\
& \geq \frac{\big(\|(\Proj_{\gamma'} - \Proj_{\gamma})
  \design_{\gammastar} \beta_{\gammastar}^\ast \|_2 - \| (
  \Proj_{\gamma'} - \Proj_{\gamma}) X_{(\gammastar)^c}
  \betastar_{(\gammastar)^c}\|_2 - \|(\Proj_{\gamma'} - \Proj_{\gamma})
  \wnoise\|_2 \big)^2}{\|\yobs\|_2^2}.
\end{align*}
By the $\betamin$-condition and Lemma~\ref{LemForwardSelection} at the end of this appendix, we have
\begin{align}
\label{EqnStrongSignal0}
\|(\Proj_{\gamma'} - \Proj_{\gamma}) \design_{\gammastar}
\beta_{\gammastar}^\ast\|_2^2 \geq  \Cm \frac{\|\big(I - \Proj_{\gamma}\big)
 \design_{\gammastar}  \beta_{\gammastar}^\ast\|^2}{|\gammastar\setminus\gamma|}\geq \numobs \Cm^2 \frac{
  \|\betastar_{\gammastar\setminus\gamma} \|_2^2}{\sstar} \geq \Cm^2 \CB
\sigmazero^2 \log \covdim,
\end{align}
where $\CB\defn \UNICON_0 (L + \Ltil + \alpha + \kappa)$ denotes the
coefficient in the $\betamin$-condition of the theorem. The last
display shows that at least an amount of $\Cm^2 \CB\sigmazero^2 \log
\covdim$ variation in the true signal $\design_{\gammastar}
\beta_{\gammastar}^\ast$ can be explained by adding the influential
covariate $\design_{j_{\gamma}}$ into the current model $\gamma$.
Combining the above two displays, we obtain that for $\numobs$
sufficiently large, so that $\Cm\sqrt{\CB}\geq 2 \sqrt{L + \Ltil}$,
the following holds under the event $\Aevent_\numobs$
\begin{align*}
1 - R_{\gamma}^2-( 1 - R_{\gamma'}^2) \geq \frac{\numobs \Cm^2
  \|\betastar_{\gammastar \setminus \gamma}\|_2^2}{4 \sstar \|\yobs\|_2^2}.
\end{align*}
Similarly, we have that under the event $\Cevent_\numobs$,
\begin{align}
1 - R_{\gamma}^2 & = \frac{Y^T(I - \Proj_{\gamma})Y}{\|\yobs\|_2^2}
\notag \\
& \leq
\frac{\big(\|( I - \Proj_{\gamma})  \design_{\gammastar}  \beta_{\gammastar}^\ast\|_2 +
  \|( \Proj_{\gamma'} - \Proj_{\gamma})  X_{(\gammastar)^c}  \betastar_{(\gammastar)^c}\|_2
  + \|(I - \Proj_{\gamma})  \wnoise \|_2\big)^2}  {\|\yobs\|_2^2} \notag \\
& \leq \frac{2\|(I - \Proj_{\gamma}) \design_{\gammastar}
  \beta_{\gammastar}^\ast \|_2^2 + 4 \|(\Proj_{\gamma'} -
  \Proj_{\gamma})
  \design_{(\gammastar)^c}\betastar_{(\gammastar)^c}\|_2^2 + 4 \|(I -
  \Proj_{\gamma}) \wnoise \|_2^2} {\|\yobs\|_2^2} \notag \\
\label{EqnR^2Bound}
& \overset{(i)}{\leq} \frac{2\|(I - \Proj_{\gamma}) \design_{\gammastar} \beta_{\gammastar}^\ast \|_2^2 + 4\tilde{L} \sigmazero^2 \log \covdim +
  3\numobs \sigmazero^2}{\|\yobs\|_2^2}
  \leq \frac{4\|(I - \Proj_{\gamma}) \design_{\gammastar} \beta_{\gammastar}^\ast \|_2^2}{\|\yobs\|_2^2}
\end{align}
where in step (i) we have used Assumption A, the fact that $\| ( I -
\Proj_{\gamma}) \wnoise\|_2^2 \leq \|\wnoise\|_2^2$, and
the last step uses $2\| ( I -\Proj_{\gamma}) \design_{\gammastar} \beta_{\gammastar}^\ast\|_2^2 \geq  2\numobs \Cm \|\betastar_{\gammastar \setminus
  \gamma}\|_2^2 \geq 4 \Ltil \sigmazero^2 \log \covdim + 3 \numobs
\sigmazero^2$ for $\CB \geq (4 \Ltil + 3)/ \Cm$.  Combining the above two
displayed inequalities, we obtain that for $\CB$ sufficiently large,
so that $\Cm^2 \CB \geq 64 \, (\alpha + \kappa + 3)$, the posterior
probability ratio $\frac{\posterior(\gamma
  \mid\yobs)}{\posterior(\gamma' \mid \yobs)}$ under the event
$\Aevent_\numobs \cap \Cevent_\numobs \cap \Devent_\numobs$ satisfies
\begin{align}
\frac{\posterior(\gamma \mid \yobs)}{\posterior(\gamma' \mid \yobs)}
&= \covdim^\kappa \sqrt{1 + \hyperpara} \cdot \Big( 1 -\frac{\|\yobs\|_2^2 ( 1
  - R_{\gamma}^2 - ( 1 - R_{\gamma'}^2) )}{\|\yobs\|_2^2/g + \yobs^T (
  I - \Proj_{\gamma} ) \yobs} \Big)^{\numobs/2} \notag \\
& \leq \covdim^\kappa \sqrt{1 + \hyperpara} \cdot \Big(1 - \frac{\Cm
  \|\big(I - \Proj_{\gamma}\big)\design_{\gammastar}  \beta_{\gammastar}^\ast\|^2
   / (4 \sstar)} {4\numobs \sigmazero^2/\sstar + 4
  \Cm \|\big(I - \Proj_{\gamma}\big)\design_{\gammastar}  \beta_{\gammastar}^\ast\|^2} \Big)^{\numobs/2} \notag \\
& \overset{(i)}{\leq} \covdim^\kappa \sqrt{1 + \hyperpara} \cdot \Big ( 1 -
\min \Big \{ \frac{\Cm}{32 \sstar},\, \frac{\Cm^2 \CB
  \log \covdim}{32 \numobs } \Big\}\Big)^{\numobs/2}
\notag\\
&\overset{(ii)}{\leq}   \covdim^\kappa \sqrt{1 + \hyperpara} \cdot
\Big( 1 - (\alpha + \kappa+ 3)\log \covdim \cdot \frac{2}{\numobs} \Big)^{\numobs/2}
\notag\\
& \leq \covdim^\kappa \cdot \covdim^{\alpha} \cdot \covdim^{-(\alpha +
  \kappa + 3)} \notag \\
& = \covdim^{-3},
\label{EqnForwardPosRatio}
\end{align}
where in step (i) we have used the inequality $a/(b+a)\geq \min\{1/2,
a/(2b)\}$ for any $a,b>0$ and inequality~\eqref{EqnStrongSignal0}, and
step (ii) follows by our assumption on $\CB$ and Assumption D on $\sstar$.

%%%%%%%%%%%%%%%%%%%%%%%%%%%%%%%%%%%%%%%%%%%%%%%%%%%%%%%%%%%%%%%%%%%%%%%%%%

\subsection{Case $\gamma$ is underfitted and saturated}

This case happens only when $\sstar\geq 1$.  Let $j_{\gamma}$ and
$k_{\gamma}$ be the indices defined in the construction of
$\Gfun(\gamma)$ in the underfitted and saturated case. Then we have
$\gamma' = \gamma \cup \{j_{\gamma}\} \setminus \{k_{\gamma}\}$. Let
$v_1 = (\Proj_{\gamma \cup \{j_{\gamma}\}} - \Proj_\gamma)
\design_{\gammastar} \beta_{\gammastar}^\ast$ and $v_2 =
(\Proj_{\gamma \cup \{j_{\gamma}\}} - \Proj_{\gamma'})
\design_{\gammastar} \beta_{\gammastar}^\ast$. Then
Lemma~\ref{LemForwardSelection}, stated and proved in
Appendix~\ref{AppLemForwardSelection}, guarantees that
\begin{equation}
\label{EqnBounda}
\begin{aligned}
\|v_1\|_2^2 & \geq \Cm \frac{\|\big(I - \Proj_{\gamma}\big)
 \design_{\gammastar}  \beta_{\gammastar}^\ast\|^2}{|\gammastar\setminus \gamma|} 
 \geq \numobs \Cm^2 \frac{\|\beta_{\gammastar \setminus
    \gamma}^\ast\|_2^2}{\sstar} \qquad \mbox{and}\\
 \|v_2\|_2^2 & \leq
\numobs\, \omega(\Xmat)
\frac{\|\beta_{\gammastar\setminus\gamma}^\ast\|_2^2}{\maxsize - \sstar}\leq
\frac{1}{2} \numobs
\Cm^2 \frac{\|\beta_{\gammastar \setminus \gamma}^\ast\|_2^2} {\sstar}
\leq \frac{1}{2}\Cm \frac{\|\big(I - \Proj_{\gamma}\big)
 \design_{\gammastar}  \beta_{\gammastar}^\ast\|^2}{|\gammastar\setminus \gamma|},
\end{aligned}
\end{equation}
under Assumption D on $\maxsize$. This inequality shows that a larger
proportion of the true signal $\design_{\gammastar}
\beta_{\gammastar}^\ast$ can be explained when the unimportant
covariate $X_{k_{\gamma}}$ is replaced with the influential covariate
$X_{j_{\gamma}}$ in the current model $\gamma$. By letting $\wtil =
\wnoise + X_{(\gammastar)^c} \betastar_{(\gammastar)^c}$ be the
effective noise, we have
\begin{align}
& 1 - R_{\gamma}^2 - ( 1 - R_{\gamma'}^2) = \frac{\yobs^T
    (\Proj_{\gamma'}-\Proj_{\gamma}) \yobs}
  {\|\yobs\|_2^2}=\frac{\yobs^T (\Proj_{\gamma \cup \{j_{\gamma}\}} -
    \Proj_{\gamma}) \yobs - \yobs^T(\Proj_{\gamma \cup \{j_{\gamma}\}}
    - \Proj_{\gamma'}) \yobs} {\|\yobs\|_2^2} \notag \\
  = &\, \frac{\|v_1\|_2^2 + 2 v_1^T \wtil + \wtil^T (\Proj_{\gamma
      \cup \{j_{\gamma}\}} - \Proj_\gamma)\wtil - \Big \{ \|v_2\|_2^2
    + 2v_2^T \wtil + \wtil^T (\Proj_{\gamma \cup \{j_{\gamma}\}} -
    \Proj_{\gamma'}) \wtil \Big \}} {\|\yobs\|_2^2} \notag \\
\geq & \, \frac{1}{\|\yobs\|_2^2}  \, \Big\{\|v_1\|_2 \, (\|v_1\|_2 - 2 \|(\Proj_{\gamma \cup
    \{j_{\gamma}\}} - \Proj_{\gamma})\wtil\|_2) - \|v_2\|_2\cdot(\|v_2\|_2 +
  2 \| ( \Proj_{\gamma \cup \{j_{\gamma}\}} - \Proj_{\gamma'})
  \wtil\|_2) \notag\\
  \label{EqnBoundb}
  &\qquad\qquad\qquad\qquad\qquad\qquad\qquad\qquad\qquad\qquad\qquad
  - \| (\Proj_{\gamma \cup \{j_{\gamma}\}}
  - \Proj_{\gamma'}) \wtil \|_2^2 \Big\},
\end{align}
where in the last step we applied the Cauchy-Schwarz inequality to the
two cross terms $v_1^T \wnoise = v_1^T
(\Proj_{\gamma\cup\{j_{\gamma}\}} - \Proj_{\gamma})\wnoise$ and $v_2^T
\wnoise = v_1^T (\Proj_{\gamma\cup\{j_{\gamma}\}} -
\Proj_{\gamma'})\wnoise$.  Note that under the event
$\Aevent_\numobs$, the Cauchy-Schwarz inequality guarantees that
\begin{align*}
\|(\Proj_{\gamma\cup\{j_{\gamma}\}} - \Proj_{\gamma}) \wtil\|_2^2 \leq
&\, 2 \|(\Proj_{\gamma \cup \{j_{\gamma}\}} - \Proj_{\gamma})
\wnoise\|_2^2 + 2 \|\design_{(\gammastar)^c}
\betastar_{(\gammastar)^c} \|_2^2 \leq 2 (L + \tilde{L}) \sigmazero^2
\log \covdim, \quad \mbox{and} \\ \|(\Proj_{\gamma \cup
  \{j_{\gamma}\}} - \Proj_{\gamma'}) \wtil\|_2^2 \leq&\,
2\|(\Proj_{\gamma \cup \{j_{\gamma}\}} - \Proj_{\gamma'})\wnoise\|_2^2
+ 2\| \design_{(\gammastar)^c} \betastar_{(\gammastar)^c} \|_2^2 \leq
2(L + \tilde{L}) \sigmazero^2 \log \covdim.
\end{align*}
Let $A^2 \defn \Cm \frac{\|\big(I - \Proj_{\gamma}\big)
 \design_{\gammastar}  \beta_{\gammastar}^\ast\|^2}{|\gammastar\setminus \gamma|} 
 \geq \numobs \Cm^2
     \min_{j\in\gammastar} |\betastar_j|^2 \geq \Cm^2 \CB \sigmazero^2
     \log \covdim$. Then, for $\CB$ large enough so that $\Cm^2 \CB \geq 32(L +
\tilde{L})$, we have, by the $\betamin$-condition and the preceding display, that under the event
$\Aevent_\numobs$
\begin{equation}
\label{EqnBoundc}
\begin{aligned}
\|(\Proj_{\gamma\cup\{j_{\gamma}\}} -
\Proj_{\gamma})\wtil\|_2 \leq \frac{A}{4} \quad\mbox{and}\quad
\|(\Proj_{\gamma\cup\{j_{\gamma}\}} -
\Proj_{\gamma'})\wtil\|_2\leq \frac{A}{4}.
\end{aligned}
\end{equation}
By the definition of $A$, we can also write inequality~\eqref{EqnBounda} as
\begin{align}
\label{EqnBoundd}
\|v_1\| \geq A \quad \mbox{and} \quad \|v_2\| \leq \frac{A}{\sqrt{2}}.
\end{align}
By plugging in the bounds~\eqref{EqnBoundc}
and~\eqref{EqnBoundd} into inequality~\eqref{EqnBoundb}, we obtain
\begin{align*}
1 - R_{\gamma}^2 - (1 - R_{\gamma'}^2) \geq &\,
\frac{A\cdot(A - A/4) - (A/\sqrt{2}) \cdot (A/\sqrt{2} + A/4) - A^2/16} {\|\yobs\|_2^2}
\geq \frac{A^2} {8\|\yobs\|_2^2},
\end{align*}
Combining this with inequality~\eqref{EqnR^2Bound}, we obtain that for $\CB$
sufficiently large so that $\Cm^2 \CB \geq 192 $, the following
holds under the event $\Aevent_\numobs \cap \Cevent_\numobs \cap
\Devent_\numobs$
\begin{align}
\frac{\posterior(\gamma \mid \yobs)} {\posterior(\gamma' \mid \yobs)} &
= \Big(1 - \frac{\|\yobs\|_2^2 (1 -
  R_{\gamma}^2 - (1 - R_{\gamma'}^2))} {\|\yobs\|_2^2/g +
  \yobs^T(I - \Proj_{\gamma}) \yobs }\Big)^{n/2} \notag\\ &\leq
\Big(1 - \frac{\Cm \|\big(I - \Proj_{\gamma}\big)  \design_{\gammastar}  \beta_{\gammastar}^\ast\|_2^2 / (8
  \sstar)}{4\numobs \sigmazero^2 / \sstar + 4\|\big(I - \Proj_{\gamma}\big)  \design_{\gammastar}  \beta_{\gammastar}^\ast\|_2^2} \Big)^{\numobs/2}
\notag \\
& \leq \Big(1 - \frac{3\log \covdim} {\numobs/2}\Big)^{\numobs/2} \leq
\covdim^{-3},\notag
\end{align}
where the last two steps follows by the same argument as for the steps
(i) and (ii) in inequality~\eqref{EqnForwardPosRatio}.

%%%%%%%%%%%%%%%%%%

\subsection{Lemma~\ref{LemForwardSelection} and its proof}
\label{AppLemForwardSelection}

Recall the definition of $j_{\gamma}$, $k_{\gamma}$ and
$\ell_{\gamma}$ in the construction of the transition function $\Gfun$
after Lemma~\ref{LemValidTranFun}. The first result in following lemma
shows that at least an amount of $\numobs \Cm^2
\|\beta_{\gammastar\setminus\gamma}^\ast\|_2^2/\sstar$ variation in
the true signal $\design_{\gammastar} \beta_{\gammastar}^\ast$ can be
explained by adding $\design_{j_{\gamma}}$ into the current model
$\gamma$. The second result shows that removing $\design_{k_{\gamma}}$
from the model $\gamma\cup \{j_{\gamma}\}$ incurs a loss in the
explained variation of at most $\numobs \, \omega(X)
\|\beta_{\gammastar\setminus\gamma}^\ast\|_2^2/(\maxsize -
\sstar)$. As a result, if $\maxsize$ satisfies the condition
$\maxsize\geq ( 2\Cm^{-2} \omega(X) + 1) \sstar$ in Assumption D, then
it is favorable to replace the unimportant covariate
$\design_{k_{\gamma}}$ with the influential covariate
$\design_{j_{\gamma}} $ in the current model $\gamma$.

\blems
\label{LemForwardSelection}
%%%
Under the conditions and notation of
Lemma~\ref{LemPosteriorConcentration}, we have:
\begin{enumerate}
\item[(a)] If $\gamma$ is underfitted, then
\begin{align*}
&\|\Proj_{\gamma \cup \{j_{\gamma}\}} \design_{\gammastar}
  \beta_{\gammastar}^\ast\|_2^2 - \|\Proj_{\gamma}
  \design_{\gammastar} \beta_{\gammastar}^\ast\|_2^2 \geq \numobs
  \Cm^2 \frac{\|\beta_{\gammastar\setminus\gamma}^\ast\|_2^2}{\sstar}.
\end{align*}
\item[(b)] If $\gamma$ is underfitted and saturated, then
\begin{align*}
\|\Proj_{\gamma \cup \{j_{\gamma}\}}
\design_{\gammastar} \beta_{\gammastar}^\ast\|_2^2 - \|\Proj_{\gamma \cup \{j_{\gamma}\} \setminus \{k_{\gamma}\}}
\design_{\gammastar} \beta_{\gammastar}^\ast\|_2^2 \leq
\numobs \, \omega(X) \frac{ \|\beta_{\gammastar\setminus\gamma}^\ast\|_2^2}
{\maxsize - \sstar}.
\end{align*}
\end{enumerate}
%%%
\elems

\begin{proof}
For each $\ell\in \gammastar \setminus \gamma$,
Lemma~\ref{LemProjection} yields
\begin{align*}
\|\Proj_{\gamma\cup\{\ell\}}
\design_{\gammastar} \beta_{\gammastar}^\ast\|_2^2 - \|\Proj_{\gamma}
\design_{\gammastar} \beta_{\gammastar}^\ast \|_2^2 &=
(\beta_{\gammastar}^\ast)^T  \design_{\gammastar}^T
\big( \Proj_{\gamma\cup\{\ell\}} - \Proj_{\gamma} \big)
\design_{\gammastar} \beta_{\gammastar}^\ast  \\ 
&=  \frac{(\beta_{\gammastar}^\ast)^T  \design_{\gammastar}^T
  \big(I - \Proj_{\gamma}\big)  \design_\ell \design_\ell^T
  \big(I - \Proj_{\gamma} \big)
  \design_{\gammastar} \beta_{\gammastar}^\ast} {\design_\ell^T  
  \big(I - \Proj_{\gamma}\big) \design_\ell}\\ &\geq
(\beta_{\gammastar}^\ast)^T \design_{\gammastar}^T
\big(I - \Proj_{\gamma}\big) \frac{\design_\ell \design_\ell^T} {\numobs}
\big(I - \Proj_{\gamma}\big)
\design_{\gammastar} \beta_{\gammastar}^\ast,
\end{align*}
where the first step follows by the idempotence of projection matrices
and the last step follows by the normalization condition in Assumption B.
By summing the preceding inequality over $\ell \in
\gammastar\setminus\gamma$, we obtain
\begin{align*}
&\sum_{\ell \in
    \gammastar\setminus\gamma} \big(\|\Proj_{\gamma \cup \{\ell\}}
  \design_{\gammastar} \beta_{\gammastar}^\ast\|_2^2 - \|\Proj_{\gamma}
  \design_{\gammastar} \beta_{\gammastar}^\ast\|_2^2 \big) \\ 
  \geq & \,
  (\beta_{\gammastar}^\ast)^T \design_{\gammastar}^T \big(I -
  \Proj_{\gamma}\big) \frac{\design_{\gammastar\setminus\gamma} 
    \design_{\gammastar\setminus\gamma}^T} {\numobs} \big( I -
  \Proj_{\gamma}\big) \design_{\gammastar}\beta_{\gammastar}^\ast\\
 =&\, (\beta_{\gammastar}^\ast)^T  \design_{\gammastar}^T
 \big(I - \Proj_{\gamma}\big)
 \frac{\design_{\gammastar\cup\gamma}  \design_{\gammastar\cup\gamma}^T}{\numobs}
 \big(I - \Proj_{\gamma}\big)
 \design_{\gammastar}  \beta_{\gammastar}^\ast\\
 \overset{(i)}{\geq} &\, \Cm
 (\beta_{\gammastar}^\ast)^T  \design_{\gammastar}^T  \big(I - \Proj_{\gamma}\big)
 \design_{\gammastar}  \beta_{\gammastar}^\ast \\
 = &\, \Cm (\beta_{\gammastar \setminus \gamma}^\ast)^T
 \design_{\gammastar \setminus \gamma}^T \big ( I - \Proj_{\gamma}
 \big) \design_{\gammastar \setminus \gamma}  \beta_{\gammastar
   \setminus \gamma}^\ast
    \overset{(ii)}{\geq}   \numobs \Cm^2
 \|\beta_{\gammastar \setminus \gamma}^\ast\|_2^2,
\end{align*}
where in step (i) we used the fact that the vector $\big( I
- \Proj_{\gamma} \big) \design_{\gammastar}\beta_{\gammastar}^\ast$
belongs to the column space of $\design_{\gammastar\cup\gamma}$ and
applied Lemma~\ref{LemUnderfit}, and step (ii) follows by applying
Lemma~\ref{LemUnderfit}.  Since $j_{\gamma}$ maximizes
$\|\Proj_{\gamma\cup\{\ell\}}
\design_{\gammastar} \beta_{\gamma_\ast}^\ast\|_2^2$ over $\ell \in
\gammastar\setminus\gamma$, the preceding inequality implies
\begin{align*}
\|\Proj_{\gamma \cup \{j_{\gamma}\}}
\design_{\gammastar} \beta_{\gammastar}^\ast\|_2^2 - \|\Proj_{\gamma}
\design_{\gammastar} \beta_{\gammastar}^\ast\|_2^2 \geq
\numobs \Cm^2 \frac{\|\beta_{\gammastar\setminus\gamma}^\ast \|_2^2} 
 {| \gammastar\setminus\gamma |} \geq
\numobs \Cm^2 \frac{\| \beta_{\gammastar\setminus\gamma}^\ast\|_2^2}
  {\sstar}.
\end{align*}
This proves the first claimed inequality.

Denote the subset $\gamma \cup \{j_{\gamma}\}$ by $\gammatil$. For any
$\gamma' \in\Mspace$, denote by $\betahat(\gamma')$ the least-squares
solution to the problem
\begin{align}
\label{EqnModifiedLS}
\min_{\beta \in \real^\covdim, \, \beta_j = 0, \, j \notin \gamma'}
\|X \beta - X_{\gammastar \setminus \gammatil } \betastar_{\gammastar
  \setminus \gammatil}\|_2^2.
\end{align}
Given this definition, some simple linear algebra leads to
\begin{align}
\label{EqnLSD}
\|\design \betahat(\gamma') - X_{\gammastar \setminus \gammatil}
\betastar_{\gammastar \setminus \gammatil}\|_2^2 = \|(I - \Proj_{\gamma'})
X_{\gammastar \setminus \gammatil} \betastar_{\gammastar \setminus
  \gammatil} \|_2^2.
\end{align}
Since $k_{\gamma} \notin \gammastar$, we have
\begin{align*}
&\|\Proj_{\gammatil} \design_{\gammastar}
  \beta_{\gammastar}^\ast\|_2^2 -
  \|\Proj_{\gammatil\setminus\{k_{\gamma}\}} \design_{\gammastar}
  \beta_{\gammastar}^\ast\|_2^2\\ =&\, \|(\Proj_{\gammatil} -
  \Proj_{\gammatil\setminus\{k_{\gamma}\}}) \, \design_{\gammastar}
  \beta_{\gammastar}^\ast\|_2^2\\ \overset{(i)}{=} &\,
  \|(\Proj_{\gammatil} - \Proj_{\gammatil\setminus\{k_{\gamma}\}}) \,
  \design_{\gammastar\setminus\gammatil}
  \beta_{\gammastar\setminus\gammatil}^\ast\|_2^2 \\ =&\, \|(I -
  \Proj_{\gammatil\setminus\{k_{\gamma}\}})
  \design_{\gammastar\setminus\gammatil}
  \beta_{\gammastar\setminus\gammatil}^\ast\|_2^2 - \|(I -
  \Proj_{\gammatil}) \design_{\gammastar\setminus\gammatil}
  \beta_{\gammastar\setminus\gammatil}^\ast\|_2^2 \\ \overset{(ii)}{=}
  & \, \|X\betahat(\gammatil\setminus\{k_{\gamma}\}) -
  X_{\gammastar\setminus\gammatil }
  \betastar_{\gammastar\setminus\gammatil}\|_2^2 - \|X\betahat(\gammatil)
  - X_{\gammastar\setminus\gammatil}
  \betastar_{\gammastar\setminus\gammatil}\|_2^2,
\end{align*}
where in step (i) we used the fact that for $k \notin \gammastar$,
$\Proj_{\gammatil}  X_{\gammastar \cap \gammatil} =
\Proj_{\gammatil \setminus \{k\} } X_{\gammastar \cap \gammatil}$, and step
(ii) follows by equation~\eqref{EqnLSD}.  This shows that the
second claimed inequality is equivalent to
\begin{align*}
\|X\betahat(\gammatil \setminus \{k_{\gamma}\}) - X_{\gammastar\setminus\gammatil}\betastar_{\gammastar\setminus\gammatil}\|_2^2 - \|X\betahat(\gammatil) - X_{\gammastar\setminus\gammatil}  \betastar_{\gammastar\setminus\gammatil}\|_2^2
\leq \numobs \omega(X)
\frac{\|\beta_{\gammastar\setminus\gammatil}^\ast\|_2^2} {\maxsize - \sstar}.
\end{align*}
We use $\betahat_j(\gammatil)$ to denote the $j$th component of
$\betahat(\gammatil)$.  By the optimality of $\betahat(\gammatil)$ for
the least-squares problem~\eqref{EqnModifiedLS}, we have
\begin{align*}
X_k^T (X_{\gammatil} \, \betahat(\gammatil) -
X_{\gammastar\setminus\gammatil} \betastar_{\gammastar \setminus
  \gammatil}) = 0, \quad \mbox{for all $k \in \gammatil$.}
\end{align*}
Therefore, for each $k\in \gammatil$, we have
\begin{align*}
\|X \betahat( \gammatil \setminus \{k\} ) -
X_{\gammastar\setminus\gammatil}\betastar_{\gammastar\setminus\gammatil}\|_2^2
& = \|X \betahat(\gammatil) - X_{\gammastar\setminus\gammatil}
\betastar_{\gammastar\setminus\gammatil} - X_k\betahat_k
(\gammatil)\|_2^2 \\ 
& = \|X \betahat(\gammatil) - X_{\gammastar\setminus\gammatil}
\betastar_{\gammastar\setminus\gammatil}\|_2^2 +
\|X_k\betahat_k(\gammatil)\|_2^2 \\ 
& \leq \|X \betahat(\gammatil) - X_{\gammastar\setminus\gammatil}
\betastar_{\gammastar\setminus\gammatil}\|_2^2 +
n|\betahat_k(\gammatil)|^2,
\end{align*}
where in the last step we used the optimality of
$\betahat(\gammatil\setminus\{k\})$ and the normalization assumption
$\|X_k\|_2^2 = \numobs$.  Then, by the definition of $k_{\gamma}$ as
the index $k$ in $\gamma \setminus \gammastar$ that minimizes
$\|X_{\gammastar} \betastar\|_2^2 - \|\Proj_{\gammatil \setminus \{k\}
} X_{\gammastar}\betastar\|_2^2 = \|X
\betahat(\gammatil\setminus\{k\}) - X_{\gammastar\setminus\gammatil}
\betastar_{\gammastar\setminus\gammatil}\|_2^2$, we have
\begin{align*}
\|X  \betahat(\gammatil\setminus\{k_{\gamma}\}) - X_{\gammastar \setminus \gammatil}\betastar_{\gammastar\setminus\gammatil}\|_2^2
&=\min_{k \in \gamma\setminus\gammastar}
\|X  \betahat(\gammatil\setminus\{k\}) - X_{\gammastar\setminus\gammatil}\betastar_{\gammastar\setminus\gammatil}\|_2^2  \\ 
&\leq
\|X \betahat(\gammatil) - X_{\gammastar\setminus\gammatil} \betastar_{\gammastar\setminus\gammatil}\|_2^2
+ \numobs \min_{k \in \gamma\setminus\gammastar}  |\betahat_k(\gammatil)|^2  \\ 
&  \leq
\|X \betahat(\gammatil) - X_{\gammastar\setminus\gammatil}  \betastar_{\gammastar\setminus\gammatil}\|_2^2
+
\numobs \frac{\|\betahat(\gammatil)\|_2^2}  {|\gamma\setminus\gammastar|}\\
 & =
\|X  \betahat(\gammatil) - X_{\gammastar\setminus\gammatil}  \betastar_{\gammastar\setminus\gammatil}\|_2^2
+ \numobs  \frac{\|\betahat(\gammatil)\|_2^2}  {\maxsize - \sstar},
\end{align*}
where last step follows since
$|\gamma\setminus\gammastar| = \maxsize - \sstar$ by the saturation of
$\gamma$. By our definition and Assumption D,
\begin{align*}
\|\betahat(\gammatil)\|_2^2 = \|(X_{\gammatil}^T X_{\gammatil})^{-1} X_{\gammatil}^T
X_{\gammastar\setminus\gammatil}  \betastar_{\gammastar\setminus\gammatil}\|_2^2
&\leq
\opnorm{(X_{\gammatil}^T  X_{\gammatil})^{-1}  X_{\gammatil}^T
  X_{\gammastar\setminus\gammatil}}^2 \,
\|\betastar_{\gammastar\setminus\gammatil}\|_2^2  \\ 
& \leq \omega(X) \,
\|\betastar_{\gammastar\setminus\gamma}\|_2^2.
\end{align*}
Combining the last two displays yields the second claimed inequality.
\end{proof}

%%%%%%%%%%%%%%%%%%%%%%%%%%%%%%%%%%%%%%%%%%%%%%%%%%%%%%%%%%%%%%%%%%%%%%%%%

\section{Proof of Lemma~\ref{LemConsistencyRatio}}
\label{SectionLemConsistencyRatio}

We divide the proof into two cases: $\gamma$ is overfitted and
underfitted.

\paragraph{Case $\gamma$ is overfitted:}

Let $k = |\gamma\setminus\gammastar|$ be the number of unimportant covariates selected by $\gamma$. Since $\gamma\neq \gammastar$, we have $k\geq 1$. Then, we can express the posterior probability ratio as
\begin{align}
\frac{\posterior(\gamma \mid \yobs)}{\posterior(\gammastar \mid \yobs)} &
= \frac{1}{ \covdim^{\kappa k} (1 + \hyperpara)^{k/2}} \cdot
\Big(\frac{1 + \hyperpara(1-R_{\gammastar}^2)}{1 + \hyperpara ( 1 -
  R_{\gamma}^2)}\Big)^{\numobs/2} \notag \\
& = \frac{1}{ \covdim^{\kappa k} (1 + \hyperpara)^{k/2}} \cdot \Big(1 +
\frac{R_{\gamma}^2 - R_{\gammastar}^2}{g^{-1} + (1 - R_{\gamma}^2)}\Big)^{n/2}
\notag.
\end{align}
Since all influential covariates are included in models
$\MODEL_{\gamma}$ and $\MODEL_{\gammastar}$, we have
\begin{align*}
1 - R_{\gamma}^2 = \frac{\|(I - \Proj_{\gamma})
   (X_{(\gammastar)^c}\betastar_{(\gammastar)^c} + \wnoise)  \|_2^2}   
  {\|\yobs\|_2^2} \quad \mbox{and} \quad
1 - R_{\gammastar}^2 = \frac{\|(I - \Proj_{\gammastar})  (X_{(\gammastar)^c}
  \betastar_{(\gammastar)^c} + \wnoise)\|_2^2}{\|\yobs\|_2^2}.
\end{align*}
Applying the Cauchy-Schwarz inequality yields
\begin{align*}
1 - R_{\gamma}^2 & \geq \frac{ \frac{1}{2} \|(I - \Proj_{\gamma})
  \wnoise\|_2^2 - \|(I - \Proj_{\gamma}) \Xmat_{(\gammastar)^c}
  \betastar_{(\gammastar)^c} \|_2^2}{\|\yobs\|_2^2} \\
& \overset{(i)}{\geq} \frac{\| (I - \Proj_{\gamma}) \wnoise\|_2^2 -
  2\|X_{(\gammastar)^c}
  \betastar_{(\gammastar)^c} \|_2^2}{2\|\yobs\|_2^2}
\overset{(ii)}{\geq} \frac{\|( I - \Proj_{\gamma}) \wnoise\|_2^2 -
  2\tilde{L} \sigmazero^2 \log \covdim}{2\|\yobs\|_2^2},
\end{align*}
where in step (i) we used the fact the projection is a non-expansive
mapping and in step (ii) we used the assumption $\|X_{(\gammastar)^c}
\betastar_{(\gammastar)^c}\|_2^2\leq \tilde{L} \sigmazero^2 \log
\covdim$.  Similarly, since $\gammastar \subset \gamma$, we can obtain
the following inequality for the quantity $R_{\gamma}^2 - R_{\gammastar}^2$
\begin{align*}
R_{\gamma}^2 - R_{\gammastar}^2 & = \frac{\|(\Proj_{\gamma} -
  \Proj_{\gammastar}) (X_{(\gammastar)^c} \betastar_{(\gammastar)^c} +
  \wnoise)\|_2^2}  {\|\yobs\|_2^2} \leq \frac{2\| (\Proj_{\gamma} -
  \Proj_{\gammastar}) \wnoise\|_2^2 + 2 \Ltil
  \sigmazero^2 \log \covdim} {\|\yobs\|_2^2}.
\end{align*}
Under the event $\Aevent_\numobs \cap \Bevent_\numobs \cap
\Cevent_\numobs$, we have $\|(I - \Proj_{\gamma}) \wnoise\|_2^2\geq
\frac{1}{2}\numobs - r K \sstar \log \covdim \geq \frac{3 \numobs}{8}$
and $\|(\Proj_{\gamma} - \Proj_{\gammastar}) \wnoise\|_2^2 \leq k L \sigmazero^2 \log
\covdim$, where we have used the assumption $4 r \, K\sstar \log
\covdim \leq \numobs $.  Combining these two inequalities with the
preceding two displays, we obtain that the posterior probability ratio
under the event $\Aevent_\numobs \cap \Bevent_\numobs \cap
\Cevent_\numobs$ is bounded as
\begin{align}
\frac{\posterior(\gamma \mid \yobs)}{\posterior(\gammastar \mid \yobs)} &
\leq \frac{1}{ \covdim^{\kappa k} (1 + \hyperpara)^{k/2}} \cdot \Big( 1 + \frac{2 (L + \Ltil )
  \log \covdim}{\numobs/8})^{\numobs/2} \notag \\
 & \overset{(i)}{\leq} \frac{1}{ \covdim^{\kappa k} (1 + \hyperpara)^{k/2}} \cdot \exp \Big\{
\frac{16 ( L + \Ltil )\log \covdim}{\numobs} \, \frac{\numobs}{2}
\Big\} = \covdim^{8 (L + \Ltil) - k (\alpha + \kappa)} \label{EqnOverfit}\\
&\overset{(ii)}{\leq}
 \covdim^{-2 k}, \notag
\end{align}
where in step (i) we used the inquality $1+x \leq e^x$ for $x\in\real$
and step (ii) follows since $\sqrt{\hyperpara}\asymp
\covdim^{\alpha}$ with $\alpha \geq 8(L + \tilde{L}) + 2 - \kappa$ according to our
choice of the hyperparameter. This proves the first part. Now we
consider the underfitted case.

%%%%%%%%%%%

\paragraph{Case $\gamma$ is underfitted:}

This case happens only when $\sstar\geq 1$.  Let $\gammatil = \gamma
\cup \gammastar$. Denote $k = |\gammastar\setminus\gamma|$ and $\ell =
|\gamma|$, then $|\gammatil \setminus \gamma| = k$,
$|\gammatil\setminus\gammastar| = k + \ell -\sstar$, and $|\gammatil|
= k + \ell \leq (K + 1) \sstar$. Since $\gamma \subset \gammatil$, we
can write
\begin{align*}
1 - R_{\gamma}^2 - ( 1 - R_{\gammatil}^2) & = \frac{\yobs^T (
  \Proj_{\gammatil} - \Proj_{\gamma}) \yobs}{\|\yobs\|_2^2} \\
& = \frac{\big\| (\Proj_{\gammatil} - \Proj_{\gamma})
  \design_{\gammastar} \beta_{\gammastar}^\ast + (\Proj_{\gammatil} -
  \Proj_{\gamma}) X_{(\gammastar)^c} \betastar_{(\gammastar)^c} +
  (\Proj_{\gammatil} - \Proj_{\gamma}) \wnoise \big\|_2^2}{\|\yobs\|_2^2}
\\
& \geq \frac{\big(\|(\Proj_{\gammatil} - \Proj_{\gamma})
  \design_{\gammastar} \beta_{\gammastar}^\ast \|_2 - \| (
  \Proj_{\gammatil} - \Proj_{\gamma}) X_{(\gammastar)^c}
  \betastar_{(\gammastar)^c}\|_2 - \|(\Proj_{\gammatil} - \Proj_{\gamma})
  \wnoise\|_2 \big)^2}{\|\yobs\|_2^2}.
\end{align*}
By the $\betamin$-condition and Lemma~\ref{LemUnderfit}, we have 
\begin{align}
\label{EqnStrongSignal}
\|(\Proj_{\gammatil} - \Proj_{\gamma}) \design_{\gammastar}
\beta_{\gammastar}^\ast\|_2^2 =& \, \|(I - \Proj_{\gamma})
\design_{\gammastar} \beta_{\gammastar}^\ast\|_2^2 \geq \, \numobs
\Cm^2 \, \|\betastar_{\gammastar\setminus\gamma} \|_2^2 \geq \Cm^2 \CB \,
k \sigmazero^2 \log \covdim,
\end{align}
where $\CB \defn \UNICON_0 (L + \Ltil + \alpha + \kappa)$ denotes the
coefficient in the $\betamin$-condition of the theorem.  Consequently,
as long as $\UNICON_0$ is sufficiently large, we can ensure that
$\Cm\sqrt{\CB}\geq 4 \sqrt{L + \Ltil}$, and hence, under the event
$\Aevent_\numobs$, we find that
\begin{align*}
1 - R_{\gamma}^2-( 1 - R_{\gamma'}^2) \geq \frac{\|\big(I - \Proj_{\gamma}\big)  \design_{\gammastar}  \beta_{\gammastar}^\ast\|^2}{4 \|\yobs\|_2^2}.
\end{align*}
Similarly, under the event $\Cevent_\numobs$, we have
\begin{align}
1 - R_{\gamma}^2 & = \frac{Y^T(I - \Proj_{\gamma})Y}{\|\yobs\|_2^2}
\notag\\
& \leq \frac{\big(\|( I - \Proj_{\gamma}) \design_{\gammastar}
  \beta_{\gammastar}^\ast\|_2 + \|( \Proj_{\gamma'} - \Proj_{\gamma})
  X_{(\gammastar)^c} \betastar_{(\gammastar)^c}\|_2 + \|(I -
  \Proj_{\gamma}) \wnoise \|_2\big)^2} {\|\yobs\|_2^2} \notag \\
& \leq \frac{2\|(I - \Proj_{\gamma}) \design_{\gammastar}
  \beta_{\gammastar}^\ast \|_2^2 + 4 \|(\Proj_{\gamma'} -
  \Proj_{\gamma})
  \design_{(\gammastar)^c}\betastar_{(\gammastar)^c}\|_2^2 + 4 \|(I -
  \Proj_{\gamma}) \wnoise \|_2^2} {\|\yobs\|_2^2} \notag\\
& \overset{(i)}{\leq} \frac{2\|\big(I - \Proj_{\gamma}\big)  \design_{\gammastar}  \beta_{\gammastar}^\ast\|^2 + 4\tilde{L} \sigmazero^2 \log \covdim +
  3\numobs \sigmazero^2}{\|\yobs\|_2^2}\leq \frac{4\|\big(I - \Proj_{\gamma}\big)  \design_{\gammastar}  \beta_{\gammastar}^\ast\|^2}{\|\yobs\|_2^2}, 
\notag
\end{align}
where in step (i) we have used Assumption A, the fact that $\| ( I -
\Proj_{\gamma}) \wnoise\|_2^2 \leq \|\wnoise\|_2^2$, and
the last step uses $2\|\big(I - \Proj_{\gamma}\big)  \design_{\gammastar}  \beta_{\gammastar}^\ast\|^2 \geq 2 \numobs \Cm \|\betastar_{\gammastar \setminus
  \gamma}\|_2^2 \geq 4 \Ltil \sigmazero^2 \log \covdim + 3 \numobs
\sigmazero^2$ for $\CB \geq (4 \Ltil + 3) / \Cm$.

Consequently, as long as $\CB$ is large enough so that $\Cm^2 \CB \geq
64 \, \big(\alpha + \kappa + 8(L + \tilde{L}) + 2 \big)$, the
posterior probability ratio $\frac{\posterior(\gamma
  \mid\yobs)}{\posterior(\gammatil \mid \yobs)}$ under the event
$\Aevent_\numobs \cap \Cevent_\numobs \cap \Devent_\numobs$ is upper
bounded as
\begin{align}
\frac{\posterior(\gamma \mid \yobs)}{\posterior(\gammatil \mid \yobs)}
&= \covdim^{\kappa k} (1 + \hyperpara)^{k/2} \cdot \Big( 1 -\frac{\|\yobs\|_2^2 ( 1
  - R_{\gamma}^2 - ( 1 - R_{\gamma'}^2) )}{\|\yobs\|_2^2/g + \yobs^T (
  I - \Proj_{\gamma} ) \yobs} \Big)^{\numobs/2} \notag \\
& \leq \covdim^{\kappa k} (1 + \hyperpara)^{k/2} \cdot \Big(1 - \frac{\|\big(I - \Proj_{\gamma}\big)  \design_{\gammastar}  \beta_{\gammastar}^\ast\|^2 / 4
  } {4\numobs \sigmazero^2/\sstar + 4
  \|\big(I - \Proj_{\gamma}\big)  \design_{\gammastar}  \beta_{\gammastar}^\ast\|^2} \Big)^{\numobs/2} \notag \\
& \overset{(i)}{\leq} \covdim^{\kappa k} (1 + \hyperpara)^{k/2}\cdot \Big ( 1 -
\min \Big \{ \frac{1}{32},\, \frac{\Cm \CB \sstar
  \log \covdim}{32 \numobs } \Big\}\Big)^{\numobs/2}
\notag\\
&\overset{(ii)}{\leq}   \covdim^{\kappa k} (1 + \hyperpara)^{k/2} \cdot
\Big( 1 - (\alpha + \kappa+ 8(L + \tilde{L}) + 2) \sstar \log \covdim \cdot \frac{2}{\numobs} \Big)^{\numobs/2}
\notag\\
& \leq \covdim^{\kappa k} \cdot \covdim^{\alpha k} \cdot \covdim^{-\sstar (\alpha + \kappa + 8(L + \tilde{L}) + 2)} \leq \covdim^{ (\kappa + \alpha) (k - \sstar) - 8(L + \tilde{L}) - 2}.
\label{EqnForwardPosRatio}
\end{align}
where in step (i) we have used the inequality $a/(b+a)\geq \min\{1/2,
a/(2b)\}$ for any $a,b>0$ and inequality~\eqref{EqnStrongSignal}, and
step (ii) follows by our assumption on $\CB$ and the assumption
$\sstar\geq 1$ made at the beginning of this underfitted case.

Since model $\MODEL_{\gammatil}$ is overfitted, by the intermediate
result~\eqref{EqnOverfit}, we have that under the event
$\Aevent_\numobs \cap \Bevent_\numobs \cap \Cevent_\numobs$
\begin{align*}
\frac{\posterior(\gammatil \mid \yobs)}{\posterior(\gammastar \mid
  \yobs)} &\leq \covdim^{8(L + \tilde{L}) - (\alpha + \kappa)\,
  |\gammatil\setminus \gammastar|} = \covdim^{8(L + \tilde{L}) -
  (\kappa + \alpha)\, ( k + \ell -\sstar)}.
\end{align*}
Combining the last two displays, we obtain that under the event
$\Aevent_\numobs \cap \Bevent_\numobs \cap \Cevent_\numobs \cap
\Devent_\numobs$
\begin{align*}
\frac{\posterior(\gamma \mid \yobs)}{\posterior(\gammastar \mid
  \yobs)} &\leq \covdim^{- (\kappa + \alpha)\, \ell - 2} \leq
\covdim^{- 2\ell - 2},
\end{align*}
where in the last step we have used Assumption C.

%%%%%%%%%%%%%%%%%%%%%%%%%%%%%%%%%%%%%%%%%%%%%%%%%%%%%%%%%%%%%%%%%%%%%%%%%

%%%%%%%%%%%%%%%%%%%%%%%%%%%%%%%%%%%%%%%%%%%%%%%%%%%%%%%%%%%%%%%%%%%%%%%%%%

%% BIBLIOGRAPHY 
\bibliographystyle{plain}

\bibliography{MCMC}
\end{document}